\documentclass[11pt]{amsart}
\usepackage{amsmath}
\usepackage{amsxtra}
\usepackage{amscd}
\usepackage{amsthm}
\usepackage{amsfonts}
\usepackage{amssymb}
\usepackage{eucal}

\usepackage{latexsym}

\numberwithin{equation}{section}

\newtheorem{claim}[subsubsection]{Claim}
\newtheorem{cor}[subsubsection]{Corollary}
\newtheorem{lem}[subsubsection]{Lemma}
\newtheorem{prop}[subsubsection]{Proposition}

\newtheorem{conj}[subsubsection]{Conjecture}
\newtheorem{thm}[subsubsection]{Theorem}
\newtheorem{defn}[subsubsection]{Definition}
\newtheorem{rem}[subsubsection]{Remark}
\newtheorem{example}[subsubsection]{Example}

\theoremstyle{definition}
\theoremstyle{remark}

\newcommand{\thmref}[1]{Theorem~\ref{#1}}
\newcommand{\secref}[1]{Sect.~\ref{#1}}
\newcommand{\lemref}[1]{Lemma~\ref{#1}}
\newcommand{\propref}[1]{Proposition~\ref{#1}}
\newcommand{\corref}[1]{Corollary~\ref{#1}}

\newcommand{\nc}{\newcommand}
\nc{\renc}{\renewcommand}
\nc{\ssec}{\subsection}
\nc{\sssec}{\subsubsection}
\nc{\on}{\operatorname}

\nc\ol{\overline}
\nc\wt{\widetilde}
\nc\tboxtimes{\wt{\boxtimes}}
\nc{\alp}{\alpha}

\nc{\ZZ}{{\mathbb Z}}
\nc{\NN}{{\mathbb N}}
\nc{\CC}{{\mathcal C}}
\nc{\OO}{{\mathbb O}}
\renc{\SS}{{\mathbb S}}
\nc{\DD}{{\mathbb D}}
\nc{\GG}{{\mathbb G}}

\nc{\Fq}{{\mathbb F}_q}
\nc{\Fqb}{\ol{{\mathbb F}_q}}
\nc{\Ql}{\ol{{\mathbb Q}_\ell}}
\nc{\id}{\text{id}}
\nc\X{\mathcal X}

\nc{\Hom}{\on{Hom}}
\nc{\Lie}{\on{Lie}}
\nc{\Loc}{\on{Loc}}
\nc{\Pic}{\on{Pic}}
\nc{\Bun}{\on{Bun}}
\nc{\IC}{\on{IC}}
\nc{\Aut}{\on{Aut}}
\nc{\rk}{\on{rk}}
\nc{\Sh}{\on{Sh}}
\nc{\IrrSh}{\on{IrrSh}}
\nc{\Perv}{\on{Perv}}
\nc{\pos}{{\on{pos}}}
\nc{\Conv}{\on{Conv}}
\nc{\Sph}{\on{Sph}}
\nc{\Sym}{\on{Sym}}
\nc{\BunBb}{\overline{\Bun}_B}
\nc{\Buno}{\overset{o}{\Bun}}
\nc{\BunPb}{{\overline{\Bun}_P}}
\nc{\BunBM}{\overline{\Bun}_{B(M)}}
\nc{\BunPbw}{{\widetilde{\Bun}_P}}
\nc{\BunBP}{\widetilde{\Bun}_{B,P}}
\nc{\GUb}{\overline{G/U}}
\nc{\GUPb}{\overline{G/U(P)}}

\nc{\Hhom}{\underline{\on{Hom}}}
\nc\syminfty{\on{Sym}^{\infty}}
\nc\lal{\ol{\lambda}}
\nc\xl{\ol{x}}
\nc\thl{\ol{\theta}}
\nc\nul{\ol{\nu}}
\nc\mul{\ol{\mu}}
\nc\Sum\Sigma
\nc{\oX}{\overset{o}{X}{}}
\nc{\M}{{\mathcal M}}
\nc{\N}{{\mathcal N}}
\nc{\F}{{\mathcal F}}
\nc{\D}{{\mathcal D}}
\nc{\Q}{{\mathcal Q}}
\nc{\Y}{{\mathcal Y}}
\nc{\G}{{\mathcal G}}
\nc{\E}{{\mathcal E}}
\nc{\CalC}{{\mathcal C}}
\nc\Dh{\widehat{\D}}

\nc{\C}{{\mathcal C}}
\nc{\K}{{\mathcal K}}
\renewcommand{\H}{{\mathcal H}}

\nc{\T}{{\mathcal T}}
\nc{\V}{{\mathcal V}}
\renc{\P}{{\mathcal P}}
\nc{\A}{{\mathcal A}}
\nc{\B}{{\mathcal B}}
\nc{\U}{{\mathcal U}}

\nc{\Gr}{\on{Gr}}

\nc{\frn}{{\check{\mathfrak u}(P)}}
\nc{\p}{\mathfrak p}
\nc{\q}{\mathfrak q}
\nc\f{{\mathfrak f}}

\nc{\qo}{{\mathfrak q}}
\nc{\po}{{\mathfrak p}}
\nc{\s}{{\mathfrak s}}
\nc\w{\text{w}}

\nc\mathi\iota
\nc\Spec{\on{Spec}}
\nc\Mod{\on{Mod}}
\nc{\tw}{\widetilde{\mathfrak t}}
\nc{\pw}{\widetilde{\mathfrak p}}
\nc{\qw}{\widetilde{\mathfrak q}}
\nc{\jw}{\widetilde j}

\nc{\grb}{\overline{\Gr}}
\nc{\I}{\mathcal I}

\nc{\lambdach}{{\check\lambda}}
\nc{\Lambdach}{{\check\Lambda}{}}
\nc{\much}{{\check\mu}}
\nc{\omegach}{{\check\omega}}
\nc{\nuch}{{\check\nu}}
\nc{\etach}{{\check\eta}}
\nc{\alphach}{{\check\alpha}}
\nc{\betach}{{\check\beta}}
\nc{\rhoch}{{\check\rho}}

\nc{\Hb}{\overline{\H}}

\nc{\BA}{{\mathbb{A}}}
\nc{\BC}{{\mathbb{C}}}
\nc{\BG}{{\mathbb{G}}}
\nc{\BM}{{\mathbb{M}}}
\nc{\BN}{{\mathbb{N}}}
\nc{\BP}{{\mathbb{P}}}
\nc{\BR}{{\mathbb{R}}}
\nc{\BZ}{{\mathbb{Z}}}
\nc{\BS}{{\mathbb{S}}}

\nc{\CA}{{\mathcal{A}}}
\nc{\CB}{{\mathcal{B}}}
\nc{\CK}{{\mathcal{K}}}
\nc{\calD}{{\mathcal{D}}}
\nc{\CE}{{\mathcal{E}}}
\nc{\CF}{{\mathcal{F}}}
\nc{\CG}{{\mathcal{G}}}
\nc{\CH}{{\mathcal{H}}}
\nc{\CI}{{\mathcal{I}}}
\nc{\CJ}{{\mathcal{J}}}
\nc{\CL}{{\mathcal{L}}}
\nc{\CM}{{\mathcal{M}}}
\nc{\CN}{{\mathcal{N}}}
\nc{\CO}{{\mathcal{O}}}
\nc{\CP}{{\mathcal{P}}}
\nc{\CQ}{{\mathcal{Q}}}
\nc{\CR}{{\mathcal{R}}}
\nc{\CS}{{\mathcal{S}}}
\nc{\CT}{{\mathcal{T}}}
\nc{\CU}{{\mathcal{U}}}
\nc{\CV}{{\mathcal{V}}}
\nc{\CW}{{\mathcal{W}}}
\nc{\CX}{{\mathcal{X}}}
\nc{\CY}{{\mathcal{Y}}}
\nc{\CZ}{{\mathcal{Z}}}

\nc{\cM}{{\check{\mathcal M}}{}}
\nc{\csM}{{\check{\mathcal A}}{}}
\nc{\oM}{{\overset{\circ}{\mathcal M}}{}}
\nc{\obM}{{\overset{\circ}{\mathbf M}}{}}
\nc{\oCA}{{\overset{\circ}{\mathcal A}}{}}
\nc{\obA}{{\overset{\circ}{\mathbf A}}{}}
\nc{\ooM}{{\overset{\circ}{M}}{}}
\nc{\osM}{{\overset{\circ}{\mathsf M}}{}}
\nc{\vM}{{\overset{\bullet}{\mathcal M}}{}}
\nc{\nM}{{\underset{\bullet}{\mathcal M}}{}}
\nc{\oD}{{\overset{\circ}{\mathcal D}}{}}
\nc{\obD}{{\overset{\circ}{\mathbf D}}{}}
\nc{\oA}{{\overset{\circ}{\mathbb A}}{}}
\nc{\cp}{{\overset{\circ}{\mathbf p}}{}}
\nc{\oU}{{\overset{\bullet}{\mathcal U}}{}}
\nc{\oZ}{{\overset{\circ}{\mathcal Z}}{}}
\nc{\ofZ}{{\overset{\circ}{\mathfrak Z}}{}}
\nc{\oF}{{\overset{\circ}{\fF}}}

\nc{\fa}{{\mathfrak{a}}}
\nc{\fb}{{\mathfrak{b}}}
\nc{\fc}{{\mathfrak{c}}}
\nc{\fd}{{\mathfrak{d}}}
\nc{\fg}{{\mathfrak{g}}}
\nc{\fgl}{{\mathfrak{gl}}}
\nc{\fh}{{\mathfrak{h}}}
\nc{\fj}{{\mathfrak{j}}}
\nc{\fk}{{\mathfrak{k}}}
\nc{\fl}{{\mathfrak{l}}}
\nc{\fm}{{\mathfrak{m}}}
\nc{\fn}{{\mathfrak{n}}}
\nc{\fu}{{\mathfrak{u}}}
\nc{\fp}{{\mathfrak{p}}}
\nc{\fr}{{\mathfrak{r}}}
\nc{\fq}{{\mathfrak{q}}}
\nc{\fs}{{\mathfrak{s}}}
\nc{\ft}{{\mathfrak{t}}}
\nc{\fsl}{{\mathfrak{sl}}}
\nc{\hsl}{{\widehat{\mathfrak{sl}}}}
\nc{\hgl}{{\widehat{\mathfrak{gl}}}}
\nc{\hg}{{\widehat{\mathfrak{g}}}}
\nc{\chg}{{\widehat{\mathfrak{g}}}{}^\vee}
\nc{\hn}{{\widehat{\mathfrak{n}}}}
\nc{\chn}{{\widehat{\mathfrak{n}}}{}^\vee}

\nc{\fA}{{\mathfrak{A}}}
\nc{\fB}{{\mathfrak{B}}}
\nc{\fD}{{\mathfrak{D}}}
\nc{\fE}{{\mathfrak{E}}}
\nc{\fF}{{\mathfrak{F}}}
\nc{\fG}{{\mathfrak{G}}}
\nc{\fK}{{\mathfrak{K}}}
\nc{\fL}{{\mathfrak{L}}}
\nc{\fM}{{\mathfrak{M}}}
\nc{\fN}{{\mathfrak{N}}}
\nc{\fP}{{\mathfrak{P}}}
\nc{\fU}{{\mathfrak{U}}}
\nc{\fV}{{\mathfrak{V}}}
\nc{\fZ}{{\mathfrak{Z}}}

\nc{\bb}{{\mathbf{b}}}
\nc{\bc}{{\mathbf{c}}}
\nc{\bd}{{\mathbf{d}}}
\nc{\be}{{\mathbf{e}}}
\nc{\bj}{{\mathbf{j}}}
\nc{\bn}{{\mathbf{n}}}
\nc{\bp}{{\mathbf{p}}}
\nc{\bq}{{\mathbf{q}}}
\nc{\bu}{{\mathbf{u}}}
\nc{\bv}{{\mathbf{v}}}
\nc{\bx}{{\mathbf{x}}}
\nc{\bs}{{\mathbf{s}}}
\nc{\by}{{\mathbf{y}}}
\nc{\bw}{{\mathbf{w}}}
\nc{\bA}{{\mathbf{A}}}
\nc{\bK}{{\mathbf{K}}}
\nc{\bB}{{\mathbf{B}}}
\nc{\bO}{{\mathbf{O}}}
\nc{\bC}{{\mathbf{C}}}
\nc{\bD}{{\mathbf{D}}}
\nc{\bH}{{\mathbf{H}}}
\nc{\bM}{{\mathbf{M}}}
\nc{\bN}{{\mathbf{N}}}
\nc{\bV}{{\mathbf{V}}}
\nc{\bW}{{\mathbf{W}}}
\nc{\bX}{{\mathbf{X}}}
\nc{\bZ}{{\mathbf{Z}}}
\nc{\bS}{{\mathbf{S}}}

\nc{\sA}{{\mathsf{A}}}
\nc{\sB}{{\mathsf{B}}}
\nc{\sC}{{\mathsf{C}}}
\nc{\sD}{{\mathsf{D}}}
\nc{\sF}{{\mathsf{F}}}
\nc{\sG}{{\mathsf{G}}}
\nc{\sH}{{\mathsf{H}}}
\nc{\sK}{{\mathsf{K}}}
\nc{\sM}{{\mathsf{M}}}
\nc{\sO}{{\mathsf{O}}}
\nc{\sQ}{{\mathsf{Q}}}
\nc{\sP}{{\mathsf{P}}}
\nc{\sS}{{\mathsf{S}}}
\nc{\sZ}{{\mathsf{Z}}}
\nc{\sfp}{{\mathsf{p}}}
\nc{\sr}{{\mathsf{r}}}
\nc{\sg}{{\mathsf{g}}}
\nc{\sff}{{\mathsf{f}}}
\nc{\sfb}{{\mathsf{b}}}
\nc{\sfc}{{\mathsf{c}}}
\nc{\sd}{{\mathsf{d}}}
\nc{\sz}{{\mathsf{z}}}

\nc{\BK}{{\bar{K}}}

\nc{\tA}{{\widetilde{\mathbf{A}}}}
\nc{\tB}{{\widetilde{\mathcal{B}}}}
\nc{\tg}{{\widetilde{\mathfrak{g}}}}
\nc{\tG}{{\widetilde{G}}}
\nc{\TM}{{\widetilde{\mathbb{M}}}{}}
\nc{\tO}{{\widetilde{\mathsf{O}}}{}}
\nc{\tU}{{\widetilde{\mathfrak{U}}}{}}
\nc{\TZ}{{\tilde{Z}}}
\nc{\tx}{{\tilde{x}}}
\nc{\tbv}{{\tilde{\bv}}}
\nc{\tfP}{{\widetilde{\mathfrak{P}}}{}}
\nc{\tz}{{\tilde{\zeta}}}
\nc{\tmu}{{\tilde{\mu}}}

\nc{\urho}{\underline{\rho}}
\nc{\uB}{\underline{B}}
\nc{\uC}{{\underline{\mathbb{C}}}}
\nc{\ui}{\underline{i}}
\nc{\uj}{\underline{j}}
\nc{\ofP}{{\overline{\mathfrak{P}}}}
\nc{\oB}{{\overline{\mathcal{B}}}}
\nc{\og}{{\overline{\mathfrak{g}}}}
\nc{\oI}{{\overline{I}}}

\nc{\eps}{\varepsilon}
\nc{\hrho}{{\hat{\rho}}}

\nc{\one}{{\mathbf{1}}}
\nc{\two}{{\mathbf{t}}}

\nc{\Rep}{{\mathop{\operatorname{\rm Rep}}}}
\nc{\Tot}{{\mathop{\operatorname{\rm Tot}}}}
\nc{\Ker}{{\mathop{\operatorname{\rm Ker}}}}
\nc{\Hilb}{{\mathop{\operatorname{\rm Hilb}}}}
\nc{\End}{{\mathop{\operatorname{\rm End}}}}
\nc{\Ext}{{\mathop{\operatorname{\rm Ext}}}}
\nc{\CHom}{{\mathop{\operatorname{{\mathcal{H}}\it om}}}}
\nc{\GL}{{\mathop{\operatorname{\rm GL}}}}
\nc{\gr}{{\mathop{\operatorname{\rm gr}}}}
\nc{\Id}{{\mathop{\operatorname{\rm Id}}}}
\nc{\de}{{\mathop{\operatorname{\rm def}}}}
\nc{\length}{{\mathop{\operatorname{\rm length}}}}
\nc{\supp}{{\mathop{\operatorname{\rm supp}}}}

\nc{\Cliff}{{\mathsf{Cliff}}}
\nc{\Fl}{\on{Fl}}
\nc{\Fib}{{\mathsf{Fib}}}
\nc{\Coh}{{\mathsf{Coh}}}
\nc{\FCoh}{{\mathsf{FCoh}}}

\nc{\reg}{{\text{\rm reg}}}

\nc{\cplus}{{\mathbf{C}_+}}
\nc{\cminus}{{\mathbf{C}_-}}
\nc{\cthree}{{\mathbf{C}_*}}
\nc{\Qbar}{{\bar{Q}}}

\nc{\bh}{{\bar{h}}}
\nc{\bOmega}{{\overline{\Omega}}}

\nc{\seq}[1]{\stackrel{#1}{\sim}}

\nc{\chA}{\check A}
\nc{\chrho}{\check \rho}
\nc{\chT}{{\check T}}
\nc{\chM}{{\check M}}
\nc{\chG}{{\check G}}
\nc{\chH}{{\check H}}

\nc{\Stab}{{\mathop{\operatorname{\rm Stab }}}}

\nc{\domwts}{{\check\Lambda^+}}
\nc{\tdomwts}{\tilde{\check{\Lambda}}^+}
\nc{\wts}{{\check\Lambda}}
\nc{\twts}{{\tilde{\check{\Lambda}}}}

\nc{\domcowts}{\Lambda^+}
\nc{\cowts}{\Lambda}
\nc{\poscowts}{\Lambda^{\on{pos}}}
\nc{\poscoroots}{R^{\on{pos}}}

\nc{\lat}{\Lambda}
\nc{\poslat}{\Lambda^{\on{pos}}}

\nc{\blambda}{{\bar\lambda}}
\nc{\risom}{\stackrel{\sim}{\rightarrow}}
\nc{\lisom}{\stackrel{\sim}{\leftarrow}}
\nc{\hAX}{{{\hat A}_X}}
\nc{\Zhecke}{\hat{Z}}
\nc{\Vect}{\on{Vect}}

\nc{\can}{{\on{can}}}

\nc{\Z}{{{}_{c}Z}}
\nc{\stZ}{{^\star Z}}
\nc{\nuZ}{{{}_{c,\geq \nu}Z}}
\nc{\zZ}{{{}_{c,\geq 0}Z}}

\nc{\Zcan}{{{}_{c}Z_{\can}}}
\nc{\bZcan}{{{}_{c}\overline{Z}}_{\can}}
\nc{\tZcan}{{{}_{c}\widetilde{Z}}_{\can}}
\nc{\stZcan}{{{}_{c}^\star Z_{\can}}}

\nc{\bBunP}{{{}_{c}\overline{\Bun}}_P}
\nc{\bzBunP}{{{}_{c,\geq 0}\overline{\Bun}}_P}
\nc{\bnuBunP}{{{}_{c,\geq \nu}\overline{\Bun}}_P}
\nc{\tBunP}{{{}_{c}\widetilde{\Bun}}_P}
\nc{\tnuBunP}{{{}_{c,\geq \nu}\widetilde{\Bun}}_P}
\nc{\tzBunP}{{{}_{c,\geq 0}\widetilde{\Bun}}_P}

\nc{\bBunB}{{{}_{c}\overline{\Bun}}_B}

\nc{\hrG}{h^{\rightarrow}_G}
\nc{\hlG}{h^{\leftarrow}_G}
\nc{\cHG}{{{}_c\H_G}}
\nc{\cHH}{{{}_c\H_H}}

\nc{\gtimes}{\stackrel{G}{\times}}
\nc{\btimes}{\stackrel{B}{\times}}
\nc{\timesBunG}{\underset{\Bun_G}{\times}}
\nc{\timesBunM}{\underset{\Bun_M}{\times}}
\nc{\timesBunP}{\underset{\Bun_P}{\times}}

\nc{\ctheta}{{{}_c\theta}}
\nc{\ceta}{{{}_c\eta}}
\nc{\ckappa}{{{}_c\kappa}}

\nc{\sJ}{{\mathsf{J}}}
\nc{\Junk}{{\mathsf{Junk}}}

\nc{\horo}{{\on{horo}}}
\nc{\const}{{\on{const}}}
\nc{\ch}{\check}

\nc{\cCZ}{{{}_c\CZ}}

\nc{\cCW}{{{}_c\CW}}

\nc{\Spf}{\on{Spf}}

\nc{\Quot}{\on{Quot}}
\nc{\Subquot}{\on{Subquot}}
\nc{\Good}{\on{Good}}
\nc{\Bad}{\on{Bad}}

\nc{\hecke}{{{}_cH^\lambda_G}}
\nc{\phecke}{{{}^p_cH^\lambda_G}}

\nc{\pcoh}{{{}^p h}}

\nc{\holstZ}{{{}^\star Z}}

\nc{\posrts}{\ch R^\pos}

\nc{\locZast}{\CV^{\eta^\pos}}

\nc{\genZastfib}{{{}^\eta W}}
\nc{\gendeepZast}{ M^{\eta}}

\nc{\gentranshecke}{{}_{\maxtor}{}'\H^{\underline\lambda}_{M}}

\nc{\ueta}{{\underline\eta}}
\nc{\ulambda}{{\underline\lambda}}
\nc{\umu}{{\underline\mu}}
\nc{\ua}{{\underline a}}
\nc{\uc}{{\underline c}}
\nc{\ud}{{\underline d}}

\nc{\disj}{{\on{disj}}}

\nc{\Sat}{{\on{Sat}}}

\nc{\affgr}{\on{Gr}}
\nc{\bdgr}{\on{Gr}^{I}}
\nc{\onebdgr}{\on{Gr}^{(1)}}

\nc{\pshonaffgr}{\mathbf P_{G(\CO)}(\affgr_G)}
\nc{\catp}{\mathbf P}
\nc{\catq}{\mathbf Q}
\nc{\catr}{\mathbf R}
\nc{\conv}{\on{Conv}}
\nc{\fact}{\on{Fact}}

\nc{\openX}{\mathring{X}}

\nc{\maxtor}{A_0}

\nc{\fcY}{{{}_{\fc}Y}}

\nc{\fcconv}{{{}_{\fc}\hspace{-0.15em}\on{Conv}}}
\nc{\fccatq}{{{}_{\fc}\mathbf Q}}

\nc{\cconv}{{{}_{c}\hspace{-0.15em}\on{Conv}}}

\nc{\fcH}{{{}_{\fc}\CH_G}}

\nc{\heckedom}{{}_{\fc}^{\ulambda} \widetilde Z}
\nc{\zastheckedom}{{}_{\fc}^{\eta} \widetilde W}

\nc{\hrS}{h^{\rightarrow}_S}
\nc{\hlS}{h^{\leftarrow}_S}

\nc{\cZ}{{{}_{c}Z}}
\nc{\fcZ}{{{}_{\fc}Z}}
\nc{\fcaZ}{{{}_{\fc}^aZ}}

\nc{\fcZstrat}{{{}_{\fc}Z^{(\Theta^+,\Theta^\pos)}}}
\nc{\genZaststrat}{{}^\eta W^{(\Theta^+,\Theta^\pos)}}

\nc{\newbase}{{}_{\uc}M^{\oeta}}
\nc{\genZast}{{{}^{\oeta}_\uc W}}

\nc{\ZI}{Z_{I}}
\nc{\ZdistinctI}{Z_{\mr I}}

\nc{\un}{\underline n}
\nc{\um}{\underline m}

\nc{\utheta}{\underline \theta}

\nc{\oeta}{\overline{\eta}}
\nc{\stareta}{\eta^\star}
\nc{\starkappa}{\kappa^\star}
\nc{\base}{M_{C,I}^{\stareta}}
\nc{\splitbase}{M_{I,J}^{\stareta,\starkappa}}
\nc{\zast}{W^{\stareta}_{C,I}}
\nc{\zaststrat}{W^{\stareta,\fp,\Theta}_{I}}
\nc{\splitzast}{W^{\stareta,\starkappa}_{I,J}}
\nc{\splitzaststrat}{W^{\stareta,\starkappa,\fp,\Theta}_{I,J}}
\nc{\distinctIzast}{{W^{\stareta}_{I}}}

\nc{\convdistinctIzast}{{\widetilde W^{\stareta,\lambda_I}_{\mr I}}}

\nc{\naivezast}{{Z^{\mu}_{I,{P^\op}}}}
\nc{\naivezastr}{{Z^{\mu}_{I,P^\op,r}}}
\nc{\ditinctInaivezast}{{Z^{\mu}_{I,P^-}}}
\nc{\distinctInaivezastr}{{Z^{\mu}_{I,P^-,r}}}
\nc{\convdistinctInaivezast}{{\widetilde Z^{\mu,\lambda_I}_{I,P^-}}}
\nc{\convdistinctInaivezastr}{{\widetilde Z^{\mu,\lambda_I}_{I,P^-,r}}}

\nc{\transhecke}{{\affgr^{\underline\mu}_{k,F,M_S}}}
\nc{\modifiedzast}{{{} W^{\hspace{0em}\stareta,\umu}_{I,k}}}
\nc{\modifiedzastr}{{{} W^{\hspace{0em}\stareta,\umu}_{I,k,r}}}
\nc{\distinctImodifiedzast}{{{} W^{\hspace{0em}\stareta,\umu}_{I,k}}}
\nc{\distinctImodifiedzastr}{{{} W^{\hspace{0em}\stareta,\umu}_{I,k,r}}}
\nc{\convdistinctImodifiedzast}{{{}\widetilde W^{\hspace{0em}\stareta,\umu,\lambda_I}_{I}}}
\nc{\convdistinctImodifiedzastr}{{{}\widetilde W^{\hspace{0em}\stareta,\umu,\lambda_I}_{I,r}}}

\nc{\nhecke}{{{}_{(n)}\CH_G}}

\nc{\distinctIhecke}{\CH_{n,G}}
\nc{\distinctIheckedomain}{{{}_{n}\widetilde Z^\ulambda}}
\nc{\distinctIzastdomain}{{{}^\stareta_{n}\widetilde  W^\ulambda}}

\nc{\nconv}{\on{Conv}_n}

\nc{\fus}{\on{Fus}}
\nc{\product}{{\circledast}}
\nc{\mr}{\mathring}

\nc{\BH}{\mathbb H}
\nc{\vect}{\on{\mathbf{Vect}}}

\nc{\pH}{{}^p H}

\nc{\hl}{{h^\leftarrow}}
\nc{\hr}{{h^\rightarrow}}
\nc{\la}{\leftarrow}
\nc{\ra}{\rightarrow}

\nc{\wh}{\widehat}

\nc{\glue}{\gamma}

\nc{\trans}{ W^{\stareta,+}_I}

\nc{\an}{{\on{an}}}
\nc{\loc}{{\on{loc}}}

\nc{\test}{{\mathcal S}}
\nc{\testdiv}{{\mathcal D}}

\nc{\hol}{{\on{hol}}}

\nc{\Filt}{\on{Filt}}

\nc{\constr}{\on{c}}
\nc{\op}{{op}}

\nc{\cla}{{\check\lambda}}

\begin{document}

\title[Spherical varieties and Langlands duality]
{Spherical varieties and Langlands duality}

\author[Gaitsgory and Nadler]{Dennis Gaitsgory and David Nadler}
\address{Harvard University, Northwestern University}
\email{gaitsgde@math.harvard.edu, nadler@math.northwestern.edu}

\begin{abstract}

Let $G$ be a connected reductive complex algebraic group.
This paper is 
devoted to the space $Z$ of meromorphic quasimaps
from a curve into an affine spherical $G$-variety $X$. The space $Z$ may be thought
of as a finite-dimensional algebraic
model for the loop space of $X$.
The theory we develop
associates to $X$ a connected reductive complex algebraic subgroup $\check H$
of the dual group $\check G$. The construction of $\check H$ is via
Tannakian formalism: we identify
a certain tensor category $\catq(Z)$ of perverse sheaves on $Z$ with the category of finite-dimensional representations of $\check H$.
The group $\check H$ encodes many aspects of the geometry of $X$.

\end{abstract}

\maketitle

\setcounter{tocdepth}{1}

\tableofcontents


\section{Introduction}

\subsection{Overview}

Let $G$ be a connected reductive complex algebraic group,
and let $S\subset G$ be a complex algebraic subgroup.
The following conditions are known to be equivalent:

\begin{enumerate}

\item There is an open $S$-orbit in the flag variety of $G$. 

\item There are finitely many $S$-orbits in the flag variety of $G$. 

\item For any $G$-equivariant line bundle $\CL\to G/S$, the multiplicity
of any irreducible $G$-representation in $H^0(G/S,\CL)$ is zero or one. 

\item For any $G$-variety $X$ and point $x\in X$ fixed by $S$,
there are a finite number of $G$-orbits in the closure of the $G$-orbit through $x$.

\end{enumerate}
(See~\cite{BLV86} for references to proofs.)

A subgroup $S\subset G$
is said to be {\em spherical}
if any of the above equivalent conditions holds.
A $G$-variety $X$ is said to be spherical
if there is a dense $G$-orbit $\openX\subset X$ and the stabilizer of any point $x\in \openX$
is a spherical subgroup. 

Examples of spherical varieties include flag varieties, symmetric spaces, and toric varieties.
In general,
spherical varieties form a class of varieties whose geometry is combinatorially tractable.
They play a fundamental role in the geometric representation theory of complex Lie algebras
and real and $p$-adic Lie groups. One can interpret the theory developed here 
as a generalization
of the geometric Satake correspondence (developed by Lusztig~\cite{L},
Ginzburg~\cite{G}, and Mirkovic-Vilonen~\cite{MV04}) from groups to spherical varieties.
In particular, this opens up the study of spherical varieties 
to the rich methods of the geometric Langlands program (see 
the papers of Beilinson-Drinfeld~\cite{BD} and Kapustin-Witten~\cite{KW} among many other important works).

\medskip

In what follows, we give a brief synopsis of the theory developed in this paper.\footnote{This paper is the combination of two earlier
preprints~\cite{GNqmI04} and~\cite{GNqmII04}.}
Our goal is to associate to any affine spherical $G$-variety $X$ a connected reductive complex algebraic
subgroup $\check H$ of the Langlands dual group $\ch G$
along with a canonical maximal torus of $\check H$. 
To define $\check H$, we construct
its category of finite-dimensional representations in the geometry of
the space $Z$ of meromorphic quasimaps from a curve into $X$.
A meromorphic quasimap consists of a point of the curve, a $G$-bundle on the curve,
and a meromorphic section of 
the associated $X$-bundle with a pole only at the distinguished point.
The space $Z$ may be thought
of as a finite-dimensional algebraic
model for the loop space of $X$.
By construction, the category of finite-dimensional representations
of $\check H$ is equivalent to a certain tensor category $\catq(Z)$ 
of perverse sheaves on $Z$.

The category of finite-dimensional representations
of $\check G$ 
naturally acts
on sheaves on $Z$ via 
the geometric Satake correspondence.
By definition, the category $\catq(Z)$ consists of all sheaves that arise
when one acts on a certain distinguished sheaf.
One of our main technical results is that objects of $\catq(Z)$
are direct sums of the intersection cohomology sheaves
of certain subspaces of $Z$.
The definition of these subspaces is local: it refers only to the restriction of a quasimap
to the formal neighborhood of its distinguished point. To understand what this restriction may look like,
we need a parametrization of the meromorphic quasimaps on a formal disk.
Since all $G$-bundles on a formal disk are trivial,
this is the same thing as a parametrization of the formal loops in 
the open $G$-orbit $\mr X\subset X$. 
A result of Luna-Vust~\cite{LV} identifies the equivalence classes of such formal loops with
the cone $\CV(X)$ of 
$G$-invariant discrete valuations of the function field of $X$.\footnote{We thank F. Knop 
for pointing out that the parametrization appears in~\cite{LV}.
Since we need some further refinements of the parametrization,
we include a proof using their 
compactification theory.}

To identify $\catq(Z)$
with the category of finite-dimensional representations
of $\check H$,
we use Tannakian formalism:
we endow $\catq(Z)$ with a tensor product,
a fiber functor, and the necessary compatibility constraints
so that it must be equivalent to the category of representations of such a group.
Recall that a category of representations comes equipped with a forgetful functor
which assigns to a representation its underlying vector space.
A fiber functor is a functor to vector spaces with all of the same properties
as the forgetful functor.
We construct  both the tensor product and fiber functor for $\catq(Z)$ 
via nearby cycles in canonically defined families. 
For the tensor product,
we consider quasimaps with multiple distinguished points, and consider
the family resulting from allowing the points to collide.
This kind of fusion product is inspired by the factorization structures of conformal field theory.
It is worth emphasizing that the fusion product 
is not derived from any group structure on $X$, but rather
from the homotopy group structure on the loop space of $X$. 
(See the end of this overview for a brief topological field theory interpretation 
of the fusion product.)

For the fiber functor, we consider the family obtained by considering quasimaps into
the specialization of $X$ to its asymptotic cone $X_0$.
The asymptotic cone $X_0$ belongs to 
a special class of $G$-varieties closely related to flag varieties.
A subgroup $S_0\subset G$
is said to be {\em horospherical} if it contains the unipotent radical of a Borel subgroup of $G$.
A $G$-variety $X_0$ is said to be horospherical if for
each point $x\in X_0$, its stabilizer is horospherical.
When $X_0$ is an affine horospherical $G$-variety,
the main result of~\cite{GNhoro04} implies that the 
category $\catq(Z_0)$ is equivalent to the category of finite-dimensional representations
of a torus.
For an arbitrary affine spherical $G$-variety $X$,
the specialization to its asymptotic cone $X_0$ provides a corresponding specialization 
for quasimaps.
By properly interpreting the nearby cycles in this family, we obtain a functor 
which corresponds to the restriction of representations from $\ch H$ to a 
canonical maximal torus.
In particular, by forgetting the torus action,
we obtain the sought-after fiber functor.

\medskip

Finally, one can interpret our results as a categorification of the structure
theory of spherical varieties.
First, let us summarize what we know about the subgroup $\ch H\subset \ch G$ 
associated to the affine spherical $G$-variety $X$. 
As stated above, we know that $\ch H$ is connected and reductive,
and it comes equipped with a canonical maximal torus. 
We know that the irreducible representations
of $\ch H$ are indexed by a subsemigroup of the cone $\CV(X)$ of $G$-invariant
discrete valuations of the function field of $X$.
One expects this subset to be of finite index in the entire cone,
or equivalently that 
the Weyl group of $\ch H$ is the same as that associated to $X$ in the structure theory
of spherical varieties~\cite{Br90,Kjams96}. 
This Weyl group controls many aspects of the geometry of $X$.
For example, by the results of~\cite{Knop94},
the center of the ring of differential operators
on $X$ is isomorphic to the invariants of the Weyl group of $\check H$ in the ring of
polynomial functions
on a Cartan algebra of $\check H$.

In the special case when $X$ is a symmetric variety, the subgroup $\check H$ coincides
with that associated to the corresponding real group in~\cite{Nad05}. One may view this
as an instance of the general Harish Chandra framework that real groups may be studied
in a complex algebraic context.

\medskip

In the next section of the introduction below, we give a more detailed description of the contents
of this paper.
Before continuing, it is worth commenting about the technical nature of 
quasimaps.
What we have called the space $Z$ of meromorphic quasimaps is not in fact a scheme
but rather an ind-stack. To the interested but daunted reader, 
we recommend thinking of $Z$ as a sophisticated version of the loop space of $X$.
In a less complicated world, we would not need to consider $Z$
but could work directly with the loop space of $X$.
Unfortunately, we do not know how to deal with the infinite-dimensionality of this space.
Another possible approach 
would be to study
sheaves on the affine Grassmannian of $G$ which are equivariant for the loop group
of the generic stabilizer of $X$. 
But in general the orbits of such a group
are neither finite-dimensional nor finite-codimensional,
and we do not know how to make sense of sheaves on them.

Although quasimaps involve many technical challenges, 
they seem to be
a natural model to study the kind of geometry we are interested in.
For example, we do not know how to see such structures as the tensor product 
on $\catq(Z)$
without some global considerations. 
This is not surprising from the perspective of the geometric Langlands program
as a topological field theory~\cite{KW}.
Under the geometric Satake correspondence,
one can interpret the tensor product
on the category of finite-dimensional representations
of the dual group $\check G$ as coming from two distinct sources.
On the one hand, there is the usual group multiplication on $G$,
which in turn induces a group structure on the loop group of $G$.
If we consider the geometric Satake category
as morphisms between two copies of the vacuum brane on a circle,
then the convolution product realizes the composition of morphisms.
On the other hand, there is the fusion product on the geometric Satake category
coming from a two dimensional pair-of-pants with boundary circles labelled by the vacuum brane.
This is an algebraic form of the homotopy group structure on the loop space of $G$
independent of the group structure on $G$.
Now we can interpret sheaves on the loop space of 
the spherical variety $X$ as a brane with which we can label the circle.
Then we can further interpret our category of spherical sheaves $\catq(Z)$
as morphisms between this brane and the vacuum brane.
Though such morphisms can not be composed (since they have different
source and target),
they can be fused: a version of a three dimensional pair-of-pants realizes the fusion product
of such morphisms. In this framework, the global nature of the fusion product
on $\catq(Z)$ is not an artifact of our technical approach
but evidence of the topological field theory structures underlying our constructions.


\subsection{Summary}
We now turn to a more detailed description of the specific contents
of this paper.

Let $C$ be a smooth complete complex algebraic curve,
and let $\Bun_G$ be the moduli stack of $G$-bundles on $C$. 
For a finite set $I$, we write $C^I$ for the product of $I$ copies of $C$,
and for a point $c_I\in C^I$, we write $|c_I|\subset C$ for the union of the points 
$c_i\in C$, for $i\in I$.

\medskip

Let $X$ be an affine spherical $G$-variety,
and let $\openX\subset X$ be the dense $G$-orbit. Fix a point $x\in\openX$, and let $S\subset G$
be the stabilizer of $x$,
so we have $\openX\simeq G/S$.
We define the ind-stack $\ZI$ of {\em meromorphic quasimaps}
to be that classifying the data
$$
(c_I\in C^I, \CP_G\in\Bun_G,\sigma:{C}\setminus {|c_I|}\to 
\CP_G{\gtimes}X|_{{C}\setminus {|c_I|}})
$$
where $\sigma$ is a section which
factors
$$
\sigma|_{C'}:C'\to 
\CP_G{\gtimes}{\openX}|_{C'}\to
\CP_G{\gtimes}{X}|_{C'},
$$
for some open curve $C'\subset C$.
We call the subset $|c_I|\subset C_\test$ the {pole points}
of the quasimap.
We call the largest open curve $C'\subset C$
on which $\sigma$ factors
$$
\sigma|_{C'}:C'\to 
\CP_G{\gtimes}{\openX}|_{C'}\to
\CP_G{\gtimes}{X}|_{C'}
$$
the {nondegeneracy locus} of the quasimap,
and we call its complement $C\setminus C'$ the {degeneracy locus}
of the quasimap.

Over the nondegeneracy locus $C'\subset C$, 
the section $\sigma$ defines a reduction of the $G$-bundle $\CP_G$
to an $S$-bundle $\CP_S$. We refer to
$\CP_S$ as the {generic $S$-bundle}
associated to the quasimap. 
We refer to the $\pi_0(S)$-bundle
induced from $\CP_S$ as the {generic $\pi_0(S)$-bundle}
associated to the quasimap. 
We call the quasimap {untwisted}
if its associated generic $\pi_0(S)$-bundle
is trivial. 

We write $'Z_I\subset Z_I$ for the ind-closed substack of untwisted meromorphic
quasimaps.
When $I$ has one element, we write $Z$ in place of $Z_{I}$,
and $'Z$ in place of $'Z_I$. 

\medskip

Let $\Sh(Z)$ be the bounded constructible derived category of sheaves of $\BC$-modules
on $Z$,
and
let $\catp(Z)$ be the full abelian subcategory of perverse sheaves 
on $Z$. We have two natural operations on $\Sh(Z)$ 
coming from {Hecke correspondences}.

First, modifications of the generic $S$-bundle of a quasimap provide
correspondences which act on $\Sh(Z)$. We refer to these as {generic Hecke
correspondences}, and to their collection
as the {generic Hecke action}.
We define
a {\em Hecke equivariant} perverse sheaf on $Z$
to be an 
object of $\catp(Z)$
equipped with a collection of isomorphisms for
smooth generic Hecke correspondences.
We write $\catp_\CH(Z)$ for the abelian category of Hecke equivariant
perverse sheaves on $Z$.

Second, modifications of the $G$-bundle of a quasimap at its pole point
also provide correspondences which act on $\Sh(Z)$. 
We refer to these as {meromorphic Hecke
correspondences}, and to their collection
as the {meromorphic Hecke action}.
To describe the types of meromorphic Hecke correspondences
which we consider, we introduce some more notation.
Let $\CK=\BC((t))$ be the field of formal Laurent series,
and let $\CO=\BC[[t]]$ be the ring of formal power series.
Let $G(\CK)$ be the $\CK$-valued points of $G$,
and let $G(\CO)$ be its $\CO$-valued points.
The affine Grassmannian $\affgr_G$ is an ind-scheme whose set of
$\BC$-points is naturally the quotient $G(\CK)/G(\CO)$.
The types of modifications of a $G$-bundle at a point
are given by the $G(\CO)$-orbits
in $\affgr_G$.
The geometric Satake correspondence~\cite{MV04}
states that the category $\pshonaffgr$ of $G(\CO)$-equivariant perverse sheaves
on $\affgr_G$ is a tensor category equivalent
to the category of finite-dimensional representations of the dual group $\ch G$.
In this paper, we consider the meromorphic Hecke correspondences 
whose integral kernels are given by objects of $\pshonaffgr$.

\medskip

Our aim is to study a certain semisimple category $\catq(Z)$ of 
Hecke equivariant perverse sheaves on $Z$.
To construct $\catq(Z)$,
we define the {\em untwisted basic stratum}
$'Z^0\subset Z$ to consist of those untwisted quasimaps of the form
$$
(c\in C, \CP_G\in\Bun_G,\sigma:{C}\to 
\CP_G{\gtimes}\openX).
$$
In other words, $'Z^0\subset Z$
consists of those quasimaps
whose nondegeneracy locus is all of $C$, and whose 
associated generic $\pi_0(S)$-bundle is trivial.
Note that if the nondegeneracy locus of  a quasimap is all of $C$,
then its associated generic $S$-bundle is defined over all of $C$.
Thus the
product of the curve $C$ with the moduli stack $\Bun_{S^0}$
for the identity component $S^0\subset S$
is naturally a $\pi_0(S)$-torsor over the basic stratum $'Z^0\subset Z$.

The intersection cohomology sheaf $'\IC_Z^{0}$ of the closure  
$'\ol {Z}^0\subset Z$  of the untwisted basic stratum
is naturally a Hecke equivariant perverse sheaf.
By acting on it by meromorphic Hecke correspondences
and taking perverse cohomology sheaves, 
we obtain a semisimple functor
$$
\conv:\pshonaffgr\to\catp_\CH(Z)
$$
which we call { convolution}.
We define $\catq(Z)$ to be the strict full subcategory of $\catp_\CH(Z)$
whose objects are isomorphic to subquotients of those Hecke equivariant
perverse sheaves arising via convolution.

We prove that the irreducible objects of $\catq(Z)$ are
isomorphic to the intersection cohomology sheaves
with constant Hecke equivariant coefficients of certain
substacks of $Z$.
To identify the substacks which may occur, 
let $\CV(X)$ be the cone
of $G$-invariant discrete valuations of the function field of $X$.
Thanks to Luna-Vust~\cite{LV},
we 
have a 
canonical parametrization
$$
G(\CO)\backslash \openX(\CK)\simeq \CV(X)
$$
which is invariant under automorphisms of the ring $\CO$.
Since all $G$-bundles on the formal disk are trivial, 
this bijection may also be thought of as parameterizing meromorphic quasimaps on the formal disk.
For a valuation $\theta\in\CV(X)$, 
we write $\openX(\CK)^\theta\subset\openX(\CK)$ for the formal loops of type $\theta$.
We define the {untwisted local stratum} $'Z^\theta\subset Z$ to consist of untwisted quasimaps 
of the form
$$
(c\in C, \CP_G\in\Bun_G,\sigma:{C\setminus c}\to 
\CP_G{\gtimes}\openX|_{C\setminus c})
$$
which when restricted to the formal neighborhood $D_c$ of the pole point $c\in C$
are represented by an element of $\openX(\CK)^\theta$
for any trivialization of the restriction of $\CP_G$ to ${D_c}$ 
and any identification of $D_c$ with the abstract formal disk.

Using direct geometric methods,
we prove the following about the convolution.

\begin{thm}
Every irreducible object of $\catq(Z)$ is isomorphic to the intersection cohomology
sheaf of an untwisted local stratum with constant Hecke equivariant coefficients.
\end{thm}

The results of this paper and those of~\cite{GNhoro04}
restrict which untwisted local strata support irreducible
objects of $\catq(Z)$ to those indexed by 
a certain finite-index subsemigroup of $\CV(X)$.
One expects that all of the strata
indexed by this subsemigroup in fact support objects of $\catq(Z)$.

The main result of this paper is the following.

\begin{thm}
The category $\catq(Z)$ is naturally a tensor category
equivalent to the category of finite-dimensional representations of a connected reductive
complex algebraic subgroup $\ch H\subset \ch G$.
\end{thm}
 
To prove this, we use Tannakian formalism.
We first construct a {fusion product} on $\catq(Z)$ making it a tensor category.
To do this, we consider the analogous category $\catq(Z_I)$
for the ind-stack $Z_I$ of meromorphic quasimaps with more than a single pole point.
We show that there is a canonical equivalence
$$
\glue_I:\catq(Z)^{\otimes I}\risom\catq(Z_I).
$$
The tensor product is then defined by allowing the pole points of $Z_I$ to collide
and taking nearby cycles.
The convolution functor is naturally a tensor functor
with respect 
to the fusion product on $\catq(Z)$ and the tensor product on $\pshonaffgr$.
Under Tannakian formalism, it corresponds to the restriction
of representations from $\ch G$.

To apply Tannakian formalism, we need a fiber functor on $\catq(Z)$.
This is an exact faithful tensor functor from $\catq(Z)$ to the category
of finite-dimensional vector spaces.
What we construct is an exact faithful tensor functor to a category of
vector spaces graded by a lattice. 
By forgetting the grading, we obtain the desired fiber functor.
Under Tannakian formalism, 
the graded fiber functor corresponds to the restriction
of representations from $\ch H$ to a maximal torus.

To describe some of the geometry involved in constructing the graded fiber functor,
we consider the case when $X$ is an affine horospherical $G$-variety.
From the results of~\cite{GNhoro04},
we deduce that in this case
the category $\catq(Z)$
is equivalent to a category of finite-dimensional vector spaces
graded by a lattice.
In other words, it is equivalent to the category
of finite-dimensional representations of a complex algebraic torus.

Now for a general affine spherical $G$-variety $X$, 
to construct the graded fiber functor on $\catq(Z)$,
we work with a family $\CX$ of $G$-varieties whose
general fiber is canonically isomorphic to $X$ and whose special fiber
is an affine horospherical
$G$-variety $X_0$.
The family $\CX$ is defined by filtering the ring of regular functions
on $X$ according to the action of $G$. The ring of regular functions
on $X_0$ is the associated graded of this filtration.
From the family $\CX$,
we obtain a family $\CZ$ of ind-stacks
whose
general fiber is canonically isomorphic to $Z$ and whose special fiber
is the ind-stack $Z_0$ of meromorphic quasimaps into $X_0$.
The nearby cycles in the family $\CZ$ provide a functor
$$
\psi:\catq(Z)\to\catp_\CH(Z_0).
$$
Unfortunately applying $\psi$ to objects of $\catq(Z)$ does not necessarily
produce objects of $\catq(Z_0)$. Instead, it may produce a complicated object
whose composition series contains an object of $\catq(Z_0)$
but also other ``noise". The presence of this noise reflects the fact that by working with $Z$,
rather than some local object, we have introduced global complications.
To deal with this, we work in a certain quotient category
of $\catp_\CH(Z_0)$ where we ignore the noise.

We call an object of $\catp_\CH(Z_0)$ a {\em bad} sheaf if 
each of its simple constituents is not an object of $\catq(Z_0)$.
We denote by $\Bad(Z_0)$ the full Serre subcategory 
of $\catp_\CH(Z_0)$
whose objects are
bad sheaves.
We denote by
$\Quot(Z_0)$ the quotient
of the category $\catp_\CH(Z_0)$
by the Serre subcategory $\Bad(Z_0)$.
We denote by $\Subquot(Z_0)$
the full image of 
$\catq(Z_0)$ under the natural projection to $\Quot(Z_0)$.
The natural projection provides an equivalence
$$
\catq(Z_0)\risom \Subquot(Z_0).
$$ 
The main technical result needed to proceed is the following.

\begin{thm}\label{thmgood}
For any object $\CP\in\catq(Z)$, the image of the 
nearby cycles $\psi(\CP)$ in the quotient category $\Quot(Z_0)$
belongs to the full subcategory $\Subquot(Z_0)$.
Thus the nearby cycles lift to a functor
$$
\Psi:\catq(Z)\to\catq(Z_0).
$$
\end{thm}

We take the functor $\Psi$ to be our graded fiber functor.
We must show that $\Psi$ is an exact faithful tensor functor.
The only part of this which is not easily deduced from previous results
is that it is a tensor functor. To see this, we must consider its interaction
with the fusion product. This involves showing that the nearby cycles
in a certain family over a base of dimension greater than one are 
independent of any choices. To see this, 
we exploit the fact that we know the assertion for objects
which come from convolution. This allows us to reinterpret
the fusion product as the fiber of a middle-extension, rather
than as nearby cycles. The characterizing properties of the middle-extension
imply the independence of the choices involved in working
over a base of dimension greater than one.

What remains to establish in order to invoke Tannakian formalism
comes for free. In~\cite{Nad05}, a brief Tannakian dictionary was collected
which contains the results we need. For example, at this point, the rigidity
of $\catq(Z)$ is automatic. Furthermore,
the assertions that the associated group $\check H$
are connected and reductive are immediate from our explicit understanding
of $\catq(Z)$.

\subsection{Acknowledgments}
We would like to thank D. Ben-Zvi for his interest and comments
on the exposition.
The second author was partially supported by NSF grant DMS-0600909.
Both authors were partially supported by DARPA grant HR0011-04-1-0031.

\newpage


{\Large \part{Main constructions}}

\vspace*{5mm}

In this part, we will carry out a series of constructions that lead to the
definition of the group $\ch H$. Most of the technical assertions will
be stated without proofs; the latter will be given in Parts II, III and IV.

The structure of this part is as follows:

\medskip

In \secref{sect quasimaps}, we introduce the space $Z$ 
of {quasimaps} from a curve $C$ into an affine variety $X$ on which $G$ acts
with an open orbit. The space $Z$ may be thought of as a finite-dimensional 
model of the space of loops $X(\CK)$. 

We then introduce the category of generic-Hecke equivariant perverse sheaves.
This may be thought of
as a model for the technically inaccessible ``category of $G(\CO)$-equivariant perverse sheaves on 
$X(\CK)$".

\medskip

In \secref{sect stratifications}, we study the space $Z$ of quasimaps under the assumption
that 
the $G$-variety $X$ is {spherical}. 

We first recall some basic facts about spherical $G$-varieties, and discuss 
the stratification of $X(\CK)$ by $G(\CO)$-orbits. We then use this to describe 
a (partial) stratification of $Z$.

At the end of the section, we specialize further to the case when $X$ is
{horospherical}, and discuss the stratification of $Z$ in this 
situation.

\medskip

In \secref{sect convolution}, we introduce the {convolution action} of the
monoidal category of spherical perverse sheaves on the affine Grassmannian
of $G$ on the category of perverse sheaves on $Z$. 

This leads to the definition of the category $\catq(Z)$ that will eventually be 
shown to be equivalent to the category of representations of a reductive subgroup 
$\ch H\subset \ch G$.

\medskip

In \secref{sect specialization}, we construct a functor from the category $\catq(Z)$
to a similar category on the space $Z_0$, where the latter is the
space of quasimaps into a canonical horospherical {\it collapse} of $X$.
The construction involves taking nearby cycles of objects of 
$\catq(Z)$, and then projecting to a suitable quotient category equivalent to $\catq(Z_0)$.

We use this functor in \secref{sect Tannakian} in order to construct a fiber
functor on $\catq(Z)$. The main idea is that the corresponding category
$\catq(Z_0)$ is equivalent to that of representations of a torus, and hence admits
a canonical forgetful functor to vector spaces.

\medskip

In \secref{sect Fusion}, we endow $\catq(Z)$ with a tensor structure
using the idea of fusing pole points.

\medskip

Finally, in \secref{sect Tannakian}, we show that the existing structures
on $\catq(Z)$ define a tensor equivalence between it and 
the category of representations of a reductive subgroup 
$\ch H\subset \ch G$. 

We identify the maximal torus of $\ch H$ and state a conjecture describing
the cone formed by its dominant weights.

\vspace*{10mm}

\section{Quasimaps}   \label{sect quasimaps}

\subsection{Definition}

Let $C$ be a smooth complete complex algebraic curve.
For a scheme $\test$, we write $C_\test$ for the product $\test\times C$, and
for an $\test$-point $c\in C(\test)$, we also write $c\subset C_\test$ for its
graph.
Let $\Bun_G$ be the moduli stack of $G$-bundles on $C$. It represents the functor 
which assigns to a scheme $\test$ the category of $G$-bundles on $C_\test$.
For a finite set $I$, we write $C^I$ for the product of $I$ copies of $C$.
For a scheme $\test$, and a point $c_I\in C^I(\test)$, 
we write $|c_I|\subset C_\test$ for the union of the graphs
$c_i\subset C_\test$.

Let $X$ be an affine $G$-variety
which we assume to have a dense $G$-orbit $\openX\subset X$. Fix a point $x\in\openX$, and let $S\subset G$
be the stabilizer of $x$,
so we have $\openX\simeq G/S$.
We define the ind-stack $\ZI$ of {\em meromorphic quasimaps}
to be that representing the functor which assigns to a scheme $\test$ the category of data
$$
(c_I\in C^I(\test), \CP_G\in\Bun_G(\test),\sigma:{C_\test}\setminus {|c_I|}\to 
\CP_G{\gtimes}X|_{{C_\test}\setminus {|c_I|}})
$$
where $\sigma$ is a section which
factors
$$
\sigma|_{C'_\test}:C'_\test\to 
\CP_G{\gtimes}{\openX}|_{C'_\test}\to
\CP_G{\gtimes}{X}|_{C'_\test},
$$
for some open subscheme $C'_\test\subset C_\test$
which is the complement $C'_\test=C_\test\setminus\testdiv$ 
of a subscheme $\testdiv\subset C_\test$ 
which is finite and flat over $\test$. 
We call the subscheme $|c_I|\subset C_\test$ the {pole points}
of the quasimap.
We call the largest subscheme $C'_\test\subset C_\test$
on which $\sigma$ factors
$$
\sigma|_{C'_\test}:C'_\test\to 
\CP_G{\gtimes}{\openX}|_{C'_\test}\to
\CP_G{\gtimes}{X}|_{C'_\test}
$$
the {nondegeneracy locus} of the quasimap,
and we call its complement $C_\test\setminus C'_\test$ the {degeneracy locus}
of the quasimap.

Over the nondegeneracy locus $C'_\test\subset C_\test$, 
the section $\sigma$ defines a reduction of the $G$-bundle $\CP_G$
to an $S$-bundle $\CP_S$. We refer to
$\CP_S$ as the {generic $S$-bundle}
associated to the quasimap. 
We refer to the $\pi_0(S)$-bundle 
induced from $\CP_S$ as the {generic $\pi_0(S)$-bundle}
associated to the quasimap. 
We call the quasimap {untwisted}
if for every geometric point $s\in \test$, the restriction of the associated generic $\pi_0(S)$-bundle
to $\{s\}\times C\subset C_\test$
is trivial. 

When $I$ has one element, we write $Z$ in place of $Z_{I}$.
We write $'Z_I\subset Z_I$ for the ind-substack of untwisted meromorphic
quasimaps. 

The following lemma will be proved in \secref{sect generic}.

\begin{lem}   \label{prime closed}
The ind-substack $'Z_I\subset Z_I$ is closed.
\end{lem}

We shall also have frequent use for the following notations.
In general, for an ind-stack of quasimaps,
we add the mark $'$ to the notation
to signify the ind-closed substack
of untwisted quasimaps.
We write $\mr C^I\subset C^I$ for the variety of those $c_I\in C^I$ such that $c_i$ is distinct
from $c_j$, for distinct $i,j\in I$.
We write $Z_{\mr I}\subset Z_I$ for the ind-open substack
of meromorphic quasimaps with distinct pole points, or in other words,
the fiber product of $Z_I$ and $C^I$ over $\mr C^I$.
In general, for an ind-stack over $C^I$, we add the mark $\mr{ }$ to the notation
to signify the fiber product with $\mr C^I\subset C^I$.

In the remainder of this paper, when defining various schemes and stacks, we 
often present the
moduli problem for geometric points and leave it to the reader to extend it
to an arbitrary base.


\subsection{Generic-Hecke equivariance}

In what follows, we write $C^{(n)}$ for the $n$th symmetric power of $C$.
For a point $c_{(n)}\in C^{(n)}$, we write $|c_{(n)}|\subset C$ for its support in $C$,
and we write $\|c_{(n)}\|\subset C$ for the finite (not necessarily reduced) subscheme it defines.

We define the ind-stack $\H_{Z_I,(n)}$ 
of {\em generic-Hecke modifications} to be that classifying data
$$
(c_I, \CP^1_G,\CP^2_G,\sigma_1,\sigma_2;
c_{(n)}, \alpha)
$$
where 
$(c_I,\CP_G^i,\sigma_i)\in Z_I$,
$c_{(n)}\in C^{(n)},$  
with  $|c_{(n)}|$ 
disjoint from the degeneracy locus of $(c_I,\CP_G^i,\sigma_i)$, 
and $\alpha$ is an isomorphism of $G$-bundles
$$
\alpha:\CP^1_G|_{C\setminus |c_{(n)}|}\risom\CP^2_G|_{C\setminus |c_{(n)}|}
$$
such that the following diagram commutes
$$
\begin{array}{ccc}
C\setminus (|c_{I}|\cup |c_{(n)}|)
& \stackrel{\sigma_1}{\to}  &
\CP^1_G\overset{G}{\times} X|_{C\setminus (|c_{I}|\cup |c_{(n)}|)}\\
\downarrow &  & \downarrow \alpha \\
C\setminus (|c_{I}|\cup |c_{(n)}|) & 
\stackrel{\sigma_2}{\to} &
\CP^2_G\overset{G}{\times} X|_{C\setminus (|c_{I}|\cup |c_{(n)}|)}
\end{array}
$$
where the left vertical map is the identity,
We call a generic-Hecke modification 
{\em trivial} if the isomorphism $\alpha$ extends to an isomorphism
over all of $C$.
We have the natural projections
$$
Z_I \stackrel{h^\la}{\leftarrow}\H_{Z_I,(n)} \stackrel{h^\ra}{\to} Z_I,
$$
and the natural projection
$$
\H_{Z_I,(n)}\to C^I\times C^{(n)}.
$$

We define a {\em smooth generic-Hecke correspondence} to be any stack $Y$ equipped
with smooth maps
$$
Z_I\stackrel{h^\la_Y}{\leftarrow} Y\stackrel{h^\ra_Y}{\rightarrow} Z_I
$$
such that for some $n$, there exists a map 
$$
Y\to \H_{Z_I,(n)}
$$ 
such that the following diagram commutes
$$
\begin{array}{ccccc}
Z_I  &  \stackrel{h^\la_Y}{\leftarrow} & Y & \stackrel{h^\ra_Y}{\to} & Z_I  \\
\downarrow & & \downarrow & & \downarrow \\
Z_I & \stackrel{h^\la}{\leftarrow} & \H_{Z_I,(n)} & \stackrel{h^\ra}{\to} & Z_I
\end{array}
$$
where the outer vertical maps are the identity.
We call a smooth generic-Hecke correspondence $Y$ {\em trivial} if the image of 
such a map
$$
Y\to \H_{Z_I,(n)}
$$ 
consists of trivial generic-Hecke modifications.
We define a morphism of smooth generic-Hecke correspondences to be
a map 
$$
p:Y_1\to Y_2
$$ 
such that 
the following diagram commutes
$$
\begin{array}{ccccc}
Z_I  &  \stackrel{h^{\la}_{Y_1}}{\leftarrow} & Y_1 & \stackrel{h^{\ra}_{Y_1}}{\to} & Z_I  \\
\downarrow & & \downarrow & & \downarrow \\
Z_I  &  \stackrel{h^{\la}_{Y_2}}{\leftarrow} & Y_2 & \stackrel{h^{\ra}_{Y_2}}{\to} & Z_I 
\end{array}
$$
where the outer vertical maps are the identity.

\subsection{Generic-Hecke equivariant sheaves}

We define a {\em generic-Hecke equivariant} perverse sheaf  on $Z_I$ to be an
object $\CF\in\catp(Z_I)$ of the category of perverse sheaves
on $Z_I$
equipped with isomorphisms
$$
I_Y:h^{\la*}_Y(\CF)\risom h^{\ra*}_Y(\CF),
$$
for every smooth generic-Hecke correspondence $Y$, satisfying the following conditions.
First, for any morphism of smooth generic-Hecke correspondences 
$$p:Y_1\to Y_2,$$
we require that
$$
I_{Y_1}=p^*(I_{Y_2}).
$$
Second, consider any stack Y,
and 
for $i\in \BZ/3\BZ$, 
consider smooth generic-Hecke correspondences $Y_i$, 
and maps
$$
p_i:Y_{}\to Y_i
$$
such that the compositions
$$
h^\la_{Y_i}\circ p_i \qquad h^\ra_{Y_i}\circ p_i
$$ 
are smooth
and the following diagrams commute
$$
\begin{array}{ccc}
Y & \stackrel{p_{i+1}}{\to} & Y_{i+1} \\
p_{i}\downarrow & & \downarrow h^\la_{Y_{i+1}} \\
Y_i  & \stackrel{h^\ra_{Y_i}}{\to} & Z_I.
\end{array}
$$
Then we require that the composite isomorphism
$$
p_3^*(I_{Y_3})\circ 
p_2^*(I_{Y_2})\circ
p_1^*(I_{Y_1})
$$
be the identity morphism.
Finally, for any trivial smooth generic-Hecke correspondence $Y$,
we require that the isomorphism $I_Y$ be the identity morphism.


Observe that a Hecke equivariant structure is determined
by its values on substacks of the 
ind-stack 
$\H_{Z_I,(1)}$ of
generic-Hecke modifications at a single point.
 
\medskip

Perverse sheaves on $Z_I$, equipped with a generic-Hecke equivariant
structure naturally form a category, which we will denote by $\catp_\CH(Z_I)$.

Observe that the pre-images in $\H_{Z_I,(n)}$ of the untwisted locus $'Z_I$ under the
projections $h^\la$ and $h^\ra$ coincide.
Hence we also may introduce the category $\catp_\CH('Z_I)$.

We will denote by $\catp_\CH(Z_{\mr I})$ and $\catp_\CH('Z_{\mr I})$
the versions of the above categories attached to the locus of distinct pole points
$Z_{\mr I}$.


\section{Stratifications}   \label{sect stratifications}

\subsection{Spherical varieties}

From now on we will assume that $X$ is a {\it spherical} $G$-variety.

\begin{defn}
$X$ is said to be {\em spherical} if a Borel subgroup of $G$
acts on $X$ with a dense orbit.
\end{defn}

For the remainder of this paper, with the exception of Sects. 
\ref{strat loops} and \ref{param of orbits} 
we will assume that $X$ is affine. Then the above
definition can be rephrased in terms of the action of $G$ on the ring 
of regular functions $\BC[X]$. As a representation of $G$, it decomposes
into a direct sum of isotypic components
$$
\BC[X]\simeq\sum_{\ch\lambda\in\domwts_G}\BC[X]_{\lambda}.
$$

\begin{defn}
$\BC[X]$ is said to be {\em simple} if each $\BC[X]_{\lambda}$
is an irreducible $G$-representation.
\end{defn}

In other words, $X$ is simple if the multiplicity of each irreducible
$V^\lambda$ in $\BC[X]$ is not greater than $1$.

\begin{prop}\label{algchar1} \cite[Theorem 1]{Pmsb86}  
$X$ is spherical if and only if $\BC[X]$ is simple.
\end{prop}

We recall below some well-known structure theory of spherical varieties.
The material here is largely borrowed from \cite{BLV86}.

\subsubsection{The associated tori}
Consider the subset $\domwts_X\subset\wts_T$ of dominant weights $\lambda\in\domwts_G$
such that $\BC[X]_{\lambda}$ is nonzero.
One shows that this is in fact a sub-semigroup. Let
$\wts_X\subset\wts_T$ be the lattice generated by $\domwts_X$.
Consider the torus
$$A:= \Spec(\BC[\wts_X]).$$
Thus the coweight lattice $\cowts_A$ of $A$ is the dual of $\wts_X$.
We define the subsemigroup $\poscowts_X\subset\cowts_A$ to 
consist of those coweights which are nonnegative on $\domwts_X\subset\wts_A$.

\begin{lem}\label{lemstrconv}
The semigroup $\poscowts_X\subset\cowts_A$ is strictly convex.
\end{lem}

\begin{proof}
By definition, the semigroup $\domwts_X$ is of full rank in the lattice $\wts_X$.
\end{proof}

Consider the saturation $\Sat(\wts_X)\subset \wts_T$ of the lattice 
$\wts_X\subset \wts_T$, and the torus
$$A_0:=\Spec(\BC[\Sat(\wts_X)]).$$

We let $F$ denote the kernel
$$1\to F\to \maxtor\to A\to 1,$$ 
or equivalently the cokernel
$$0\to\cowts_{A_0}\to\cowts_A\to F\to 0.$$

\begin{rem}
The tori $A$ and $A_0$ can be associated to any (i.e., not necessarily
affine) spherical variety, and they only depend on the homogeneous
space $\openX=G/S$. This more general theory will be reviewed in
\secref{param of orbits}.
\end{rem}

\subsubsection{The associated parabolic}

Choose a Borel subgroup $B^\op\subset G$, and let $\openX^+\subset \openX\subset X$ 
be the corresponding open $B^\op$-orbit. Let $P^\op\subset G$ be the parabolic, 
consisting of elements  $g\in G$ such that $g\cdot \openX^+\subset \openX^+$, and 
let $U^{op}\subset P^{op}$ be its unipotent radical.


The choice of $B^\op$ defines highest weight lines
$\fl^\lambda\subset \BC[X]_\lambda$. 

\begin{lem} \hfill

\smallskip

\noindent (1) For $\lambda\in \wts_X^+$, the line $\fl^\lambda\subset V^\lambda$ is
$P^\op$-stable.

\smallskip

\noindent (2) The open subvariety $\openX^{+}\subset X$ is the locus of non-vanishing 
of the lines $\fl^\lambda$.
\end{lem}

Thus we see that $A$ is naturally a quotient of $P^\op$. Moreover, the 
projection $P^{op}\to A$ factors as surjections
$$
P^{op}\twoheadrightarrow A_0\twoheadrightarrow A,
$$
where the first arrow has a connected kernel. Furthermore, we have
$$
\BC[X]^{U^\op}\simeq \BC[\wts_X^+].
$$ 
Setting 
$$\overline{A}:=\Spec(\BC[\wts_X^+]),
$$ 
we have a Cartesian diagram
\begin{equation} \label{proj to A}
\CD
\openX^+  @>>> X  \\
@V{p}VV   @V{\overline{p}}VV  \\
A @>>> \overline{A},
\endCD
\end{equation}
where the vertical maps are $P^\op$-equivariant, and the horizontal ones
are open embeddings. 

\begin{lem}
The action of $U^\op$ on $\openX^+$ is free, i.e., $\openX^+$ is a principal
$U^\op$-bundle over $A$.
\end{lem}

\subsubsection{The finite groups}

By construction, we have a map of stacks
$\openX^+/P^{op}\to \openX/G$, which induces a map on the level
of the corresponding fundamental groups
$$F\simeq \pi_1(\openX^+/P^{op})\simeq \pi_1(\openX/G)\simeq \pi_0(S).$$
Since $\openX^+\to \openX$ is an open embedding, we obtain that
the above map is surjective. In particular, we deduce that $\pi_0(S)$
is abelian.

\subsection{Stratification of loops}   \label{strat loops}

As usual, let $\CK=\BC((t))$ be the field of formal Laurent series, and let $\CO=\BC[[t]]$
be the ring of formal power series.
Let $G(\CK)$ be the group of $\CK$-valued points of $G$,
and let $G(\CO)$ be the group of $\CO$-valued points of $G$.
The affine Grassmannian $\affgr_G$ is an ind-scheme whose set of
$\BC$-points is naturally the quotient $G(\CK)/G(\CO)$.

The following result, which essentially follows from \cite{LV},
will be proved in Sects.~\ref{param of orbits} and \ref{proof of orbits}.

\begin{thm} \label{thm count orbits}
For a subgroup $S\subset G$ the following conditions are equivalent:

\smallskip

\noindent{\em (1)}
The quotient $G/S$ is a spherical $G$-variety.

\smallskip

\noindent{\em (2)} The group $S(\CK)$ acts on $\affgr_G$ 
with countably many orbits.

\smallskip

\noindent{\em (3)}
The group $G(\CO)$ acts on $(G/S)(\CK)$
with countably many orbits.
\end{thm}

In what follows for $S\subset G$ satisfying the equivalent conditions 
of the theorem, we will denote by $\CV(G/S)$ the set of $G(\CO)$-orbits
on $(G/S)(\CK)$. (We will see during the proof of Theorem~\ref{thmparamval} below in \secref{param of orbits} that $\CV(G/S)$
coincides with the cone of $G$-invariant discrete valuations of the function field of $X$.)

\begin{rem}\label{remab}
Condition (3) in the Theorem tautologically implies condition (2).
The inverse implication is less evident, since not every $\CK$-point
of $G/S$ lifts to a $\CK$-point of $G$, because $S$ may be
disconnected.
\end{rem}

\subsection{Description of orbits}
For the remainder of this section we let $X$ be an {\it affine, spherical}
$G$-variety. Recall that $\openX\subset X$ denotes an open $G$-orbit,
which can be identified with the quotient $G/S$. The next result that 
we are going to state identifies the set of
$G(\CO)$ orbits on $\openX(\CK)$ with a subset of the lattice $\cowts_A$.

\medskip

Note that $X(\CK)$ is naturally (the set of $\BC$-points) of an ind-affine ind-scheme.
It contains an open subscheme $X(\CK)\setminus (X\setminus\openX)(\CK)$, whose set
of $\BC$-points identifies with the set of $\CK$-points of 
$\openX$.\footnote{Note that since $\openX$ may not be affine, the set of $\CK$-points of $\openX$
does not a priori have an ind-scheme structure. Even when $\openX$ is affine,
$\openX(\CK)$ is {\it not} isomorphic to $X(\CK)\setminus(X\setminus\openX)(\CK)$.}

\medskip

Similarly, $\overline{A}(\CK)$ is naturally an ind-affine ind-scheme,
which contains an open sub-scheme $\overline{A}(\CK)\setminus(\overline{A}\setminus A)(\CK)$.
The set of $A(\CO)$-orbits on $\overline{A}(\CK)\setminus(\overline{A}\setminus A)(\CK)$ identifies
naturally with $\cowts_A$.\footnote{An orbit corresponding to $\lambda_1\in \cowts_A$ is
contined in the closure of the orbit corresponding to $\lambda_2$ if and only if
$\cla_1-\cla_2\in \poscowts_X$.}

\medskip

Let $\bO\subset X(\CK)\setminus(X\setminus\openX)(\CK)$ be a $G(\CO)$-orbit. 
We can view it as a (set of $\BC$-points) of a scheme of infinite type.
Let $\bO^+\subset \bO$ be the open subscheme
$$\bO\cap \left(X(\CK)\setminus(X\setminus\openX^+)(\CK)\right).$$
It is non-empty, since $G/P^\op$ is proper, and hence
$\CO$- and $\CK$-points of $G/P^\op$ are in bijection.
We have a Cartesian diagram
$$
\CD
\bO^+  @>>> \bO \\
@V{p}VV    @V{\overline{p}}VV  \\
\overline{A}(\CK)\setminus(\overline{A}\setminus A)(\CK)   @>>>  A(\CK),
\endCD
$$
induced by \eqref{proj to A}.

Since $\bO^+$ is irreducible, its image under $p$ is contained
in the closure of a single $A(\CO)$-orbit on $\overline{A}(\CK)\setminus(\overline{A}\setminus A)(\CK)$.
Thus, we obtain a map
$$\upsilon: \CV(G/S)\simeq \left(X(\CK)\setminus(X\setminus\openX)(\CK)\right)/G(\CO)\to
\left(\overline{A}(\CK)\setminus(\overline{A}\setminus A)(\CK)\right)/A(\CO)\simeq \cowts_A.$$

The following reformulation of a result of \cite{LV} 
will be proved in Sects. \ref{param of orbits} and \ref{proof of orbits} using their compactification
theory of spherical varieties.

\begin{thm} \label{thmparamval}  \hfill

\smallskip

\noindent(1) The map $\upsilon$ defines a bijection between $\CV(G/S)$ and a 
finitely generated saturated subsemigroup of full rank in $\cowts_A$.  

\smallskip

\noindent(2) An orbit
$\bO$ is contained in $X(\CO)$ if and only if 
$\upsilon(\bO)\in \poscowts_X$.
\end{thm}

\begin{rem}
As was remarked above, the lattice $\cowts_A$ is naturally associated
to the homogeneous space $\openX=G/S$, i.e., the additional data of
the affine variety $X$ containing $\openX$ is, in fact, redundant. As
will become clear in \secref{proof of orbits}, the map $\upsilon$ also
depends only on $\openX$.

By contrast, the subsemigroup $\poscowts_X\subset \cowts_A$ does
depend on $X$.
\end{rem}

Let $\Aut(\CO)$ be the group-scheme
of automorphisms of $\CO$. It naturally acts on $X(\CK)$.

\begin{cor}  \label{corparam}
The orbits of $G(\CO)$ on $X(\CK)\setminus(X\setminus\openX)(\CK)$ are stable 
under the $\Aut(\CO)$-action.
\end{cor}

\subsection{Stratification of quasimaps}\label{sectstratqm}

We content ourselves here with defining the strata of the ind-stack
$Z_I$ of meromorphic quasimaps into an affine spherical $G$-variety $X$
which play a role in what follows. 
The interested reader will be able to use the results of the previous sections
to extend the definitions given here and give a complete stratification of $Z_I$.
Recall that we write $\openX\simeq G/S\subset X$
for the dense $G$-orbit in $X$. 

To a quasimap 
$
(c_I,\CP_G,\sigma)
\in Z_I$,
and a point $c\in C$, we may associate an element 
$\bar\sigma_c\in \CV(G/S)$ 
as follows.
First, if we fix 
an isomorphism of the restriction $\CP_G|\CO_c$ with the trivial bundle $\CO_c\times G$,
then the restriction $\sigma|\CK_c$
may be thought of as a loop $\sigma_c\in \openX(\CK_c)$.
Another choice of trivialization of $\CP_G|\CO_c$ will lead to another loop
in the $G(\CO_c)$-orbit in $\openX(\CK_c)$ through $\sigma_c$.
Thus
we have constructed a well-defined element $\bar\sigma_c\in  
\openX(\CK_c)/G(\CO_c)$.
Now, if we  fix 
an identification of the completed local ring $\CO_c$ with the power series ring $\CO$,
then we obtain an element $\bar\sigma_c\in \openX(\CK)/G(\CO)$.
By Corollary~\ref{corparam}, the $G(\CO)$-orbits in $\openX(\CK)$ 
are invariant under the action
of $\Aut(\CO)$, so this element is well-defined.

We are now ready to define what we call the {\em local strata}
 of the
ind-stack $Z_I$ of meromorphic quasimaps.
For a partition $\fp$ of the set $I$, and a labelling $\Theta:\fp\to\CV(G/S)$,
we say that a quasimap
$
(c_I,\CP_G,\sigma)
\in Z_I
$
is of type $(\fp,\Theta)$ if the following conditions hold.
First, the section $\sigma$ factors
$$
\sigma:{C\setminus |c_I|}\to 
\CP_G{\gtimes}\openX|_{C\setminus |c_I|}\to
\CP_G{\gtimes}X|_{C\setminus |c_I|}.
$$
Second, the coincidences among the pole points $|c_I|\subset C$ are given by the partition $\fp$.
And third, the $\CV(G/S)$-valued labelling of $\fp$
associated to the quasimap 
is given by $\Theta$.
We define the {local stratum}
$$
Z_I^{(\fp,\Theta)}\subset {}Z_{I}
$$
to consist of those quasimaps of type $(\fp,\Theta)$.
When $\fp$ is the complete partition of $I$ into singleton parts,
we write $Z_{\mr I}^\Theta$ in place of $Z_I^{(\fp,\Theta)}$.
When $I$ has a single element, the partition is vacuous, the labelling $\Theta$ is a single
element $\theta\in\domcowts_A$, and we write $Z^\theta$
in place of $Z^{(\fp,\Theta)}$.

Note that for $\Theta=0$, the corresponding stratum $Z_I^{(\fp,0)}$ is
isomorphic to $\mr C^\fp\times \Bun_S$, where $\mr C^\fp\subset C^I$
is the locally closed subset determined by the partition $\fp$.

\medskip

In what follows, we shall be more interested in the {untwisted local strata}
$$
'Z_I^{(\fp,\Theta)}\subset {}'Z_{I}
$$
consisting of those untwisted quasimaps of type $(\fp,\Theta)$.
As above, 
when $\fp$ is the complete partition of $I$ into singleton parts,
we write $'Z_{\mr I}^\Theta$ in place of $'Z_I^{(\fp,\Theta)}$.
When $I$ has a single element, the partition is vacuous, the labelling $\Theta$ is a single
element $\theta\in\domcowts_A$, and we write $'Z^\theta$
in place of $'Z^{(\fp,\Theta)}$.

We will denote by $\IC^{\fp,\Theta}_I$, respectively $'\IC^{\fp,\Theta}_I$, $\IC^\Theta_{\mr I}$,
$'\IC^\Theta_{\mr I}$, $\IC^\theta$, $'\IC^\theta$ the intersection cohomology
sheaves on the closures of the respective locally closed subvarieties of the
space of quasimaps.


\subsection{Relation to Hecke equivariance}

Fix $c _I\in C^I$, and let $\fp$ be the partition of $I$ 
describing the coincidences among the points $|c_I|\subset C$.
Consider the ind-stack of {\em based} untwisted meromorphic quasimaps
$$
'Z_{c_I}={}'Z_I\underset{C^I}{\times} \{c_I\},
$$
and the corresponding based untwisted local strata
$$
'Z_{c_I}^{(\fp,\Theta)}={}'Z_I^{(\fp,\Theta)}\underset{C^I}{\times} \{c_I\}.
$$
Observe that the generic-Hecke modifications 
$$
Z_I \stackrel{h^\la}{\leftarrow}\H_{Z_I,(n)} \stackrel{h^\ra}{\to} Z_I,
$$
preserve the based untwisted local strata
$'Z_{c_I}^{(\fp,\Theta)}$.
In other words, we obtain diagrams
$$
\begin{array}{ccccc}
'Z_{c_I}^{(\fp,\Theta)} & \stackrel{h^\la}{\leftarrow} & \H_{'Z_{c_I}^{(\fp,\Theta)},(n)}
& \stackrel{h^\ra}{\to}  & 'Z_{c_I}^{(\fp,\Theta)}\\
\downarrow & & \downarrow & & \downarrow \\
Z_I & \stackrel{h^\la}{\leftarrow} & \H_{Z_I,(n)} & \stackrel{h^\ra}{\to}  & Z_I\\
\end{array}
$$
in which each square is Cartesian.

In particular, it makes sense to introduce the categories $\catp_\CH(Z_I^{\fp,\Theta})$,
$\catp_\CH('Z_I^{\fp,\Theta})$, 
$\catp_\CH(Z_{\mr I}^{\Theta})$,
$\catp_\CH('Z_{\mr I}^{\Theta})$,
$\catp_\CH(Z_{c_I}^{\fp,\Theta})$,
$\catp_\CH('Z_{c_I}^{\fp,\Theta})$,  etc. 
And the following are naturally generic-Hecke equivariant  objects:
$$\IC^{\fp,\Theta}_I\in \catp_\CH(Z_I), \,\,\,
{}'\IC^{\fp,\Theta}_I\in \catp_\CH('Z_I),\,\,\,
\IC^\Theta_{\mr I}\in \catp_\CH(Z_{\mr I}), \,\,\,
'\IC^\Theta_{\mr I}\in \catp_\CH({}'Z_{\mr I}), \,\,\, \mbox{etc.}$$


\medskip

We call two geometric points of $'Z_{c_I}^{(\fp,\Theta)}$ {\em generic-Hecke equivalent}
if they are equivalent under the equivalence relation on the geometric
points of $'Z_{c_I}^{(\fp,\Theta)}$ generated by the generic-Hecke correspondences.
One of the reasons for working with untwisted quasimaps is the following result,
which will be proved in \secref{sect generic}. We note that 
a similar statement without the 
restriction to untwisted quasimaps would not be true.

\begin{prop}\label{prophecketrans}
All geometric points of the based untwisted local stratum 
$'Z_{c_I}^{(\fp,\Theta)}\subset Z_I$
are generic-Hecke equivalent.
\end{prop}

This proposition is an analogue of the statement that a group $\sH$
acts transitively on a space $\sZ$. 
The next proposition is an analogue of the following toy situation.
If a group $\sH$ acts transitively on a space $\sZ$, then an $\sH$-equivariant local
system on $\sZ$ is the same thing as a representation of the component
group $\pi_0(\sH_z)$ of the stabilizer $\sH_z\subset\sH$ 
of a point $\sz\in \sZ$. If the local system is constant
as an ordinary local system, then the representation must factor through
the map $\pi_0(\sH_z)\to\pi_0(\sH_{[z]})$ where the subgroup $\sH_{[z]}\subset\sH$
is that preserving the connected component of $\sz\in \sZ$.
Thus if $H_{[z]}$ is connected, then any $\sH$-equivariant local system on $\sZ$ whose
underlying ordinary local system is constant is itself isomorphic to the standard
constant $\sH$-equivariant local system of the same rank.

\begin{prop}\label{propheq}
All generic-Hecke equivariant local systems
on the untwisted local stratum $'Z_I^{(\fp,\Theta)}\subset Z_I$,
whose underlying ordinary local system is constant,
are isomorphic to the standard constant generic-Hecke equivariant local system
of the same rank.
\end{prop}

The proof will also be given in \secref{sect generic}.

\subsection{The horospherical case} 

\subsubsection{Horospherical varieties}

Recall that a $G$-variety $X$ is {\it horospherical} if for each point $x\in X$,
its stabilizer $S_x\subset G$ contains the unipotent radical of a Borel subgroup of $G$.

We have the following well-known characterization when $X$ is affine. Let 
$$\BC[X]\simeq\sum_{\ch\lambda\in\wts_T}\BC[X]_{\lambda}$$
be its ring of function broken into $G$-isotypic components.

\begin{lem}  \label{algchar2} 
$X$ is horospherical if and only if $\BC[X]_\lambda\cdot 
\BC[X]_\mu\subset \BC[X]_{\lambda+\mu}$,
for all $\lambda,\mu\in\domwts_G$.
\end{lem}

Until the end of this subsection, we will assume that $X$ is a spherical
and horospherical affine variety. From \lemref{algchar2} we obtain that
the torus $A$ acts naturally on $X$, commuting with the action of $G$.
The following is easy to verify.

\begin{lem}
For $X$ horospherical, the map
$F\to \pi_0(S)$ is an isomorphism.
\end{lem}

\subsubsection{Loops into horospherical varieties}

For $X$ horospherical the subsemigroup $\CV(G/S)$ coincides
with the entire $\cowts_A$.  The corresponding bijection
$$\CV(G/S)\risom \cowts_A$$
can be explicitly described as follows.
Consider the unit $G(\CO)$-orbit 
$$\bO_0:=\openX(\CO)\subset X(\CK).
$$ 
Using the $A$-action of $X$ mentioned above, we can translate this
orbit by an element of $A(\CK)/A(\CO)\simeq \cowts_A$ and obtain another
$G(\CO)$-orbit. This describes the desired bijection.

\section{Convolution}     \label{sect convolution}

\subsection{The convolution diagram}

For a finite set $I$,
we define the ind-stack $\CH_{I,G}$ of  {\em meromorphic Hecke modifications}
to be that classifying data
$$
(c_I,\CP^1_{G},\CP^2_G,\alpha)
$$
where $c_I\in C^I$, 
$\CP^1_{G},\CP^2_G\in\Bun_{G},$ 
and $\alpha$ is a $G$-equivariant isomorphism 
$$
\alpha:
\CP^1_{G}|_{C\setminus |c_I|}\risom
\CP^2_G|_{C\setminus |c_I|}.
$$
We call the subset $|c_I|\subset C$ {modification points}.
We have the natural projection
$$
\CH_{I,G}\to C^I,
$$
and also the natural projections
$$
\Bun_G \stackrel{{}h^\la}{\leftarrow} \CH_{I,G} \stackrel{{}h^\ra}{\to}  \Bun_G.
$$ 

As usual, let $\mr C^I\subset C^I$
consist of those $c_I\in C^I$
such that $c_i$ is distinct from $c_j$, for distinct $i,j\in I$.
We write $\CH_{\mr I,G}$ for the fiber product
$$
\CH_{\mr I,G}=\CH_{I,G}\underset{C^I}{\times} \mr C^I.
$$
Consider the $(G(\CO)\rtimes\Aut(\CO))^I$-torsor 
$$
\widehat\Bun_{\mr I,G}\to\Bun_G
$$ 
that classifies data
$$
(c_I,\CP_G,\mu,\tau)
$$
where $c_I\in \mr C^I$, $\CP_G\in\Bun_G$, $\mu$ is an isomorphism of $G$-bundles
$$
\mu:D_{|c_I|}\times G\risom \CP_G|_{D_{|c_I|}},
$$
and $\tau$ is an isomorphism
$$
\tau:D\times |c_I|\risom D_{|c_I|}.
$$
We may identify $\CH_{\mr I,G}$ with the twisted product
$$
\CH_{\mr I,G}\simeq 
{\widehat\Bun}_{\mr I,G}\overset{(G(\CO)\rtimes\Aut(\CO))^I}{\times}(\affgr_G)^I
$$
so that $h^\ra$ corresponds to the evident projection to $\Bun_G$. 
Here as usual $\affgr_G$ denotes the affine Grassmannian of $G$.

For the ind-stack $Z_{I}$ of meromorphic quasimaps,
we have the diagram
$$
\begin{array}{ccccc}
Z_{I} & \stackrel{ h^\la}{\leftarrow} & 
\CH_{I,G}\underset{C^I\times\Bun_G}{\times}Z_{I}
& \stackrel{ h^\ra}{\to} & Z_{I} \\
\downarrow & & \downarrow & & \downarrow \\
\Bun_G & \stackrel{{} h^\la}{\leftarrow} & \CH_{I,G} & \stackrel{{} h^\ra}{\to} & \Bun_G
\end{array}
$$
in which each square is Cartesian.
Restricting the diagram to $\mr C^I\subset C^I$,
we obtain a similar diagram for the ind-stack $Z_{\mr I}$ of meromorphic
quasimaps with distinct pole points.
Consider the $(G(\CO)\rtimes\Aut(\CO))^I$-torsor
$$
{\widehat Z}_{\mr I}\to Z_{\mr I}
$$
that classifies data
$$
(z,\mu,\tau)
$$
where $z\in Z_{\mr I}$, with $G$-bundle $\CP_G\in\Bun_G$ and pole points $c_I\in\mr C^I$, 
$\mu$ is an isomorphism of $G$-bundles
$$
\mu:D_{|c_I|}\times G\risom \CP_G|_{D_{|c_I|}},
$$
and $\tau$ is an isomorphism 
$$
\tau:D\times |c_I|\risom D_{|c_I|}.
$$
For the twisted product
$$
\wt Z_{\mr I}={\widehat Z}_{\mr I}\overset{(G(\CO)\rtimes\Aut(\CO))^I}{\times}(\affgr_G)^I,
$$
we have the obvious identification
$$
\wt Z_{\mr I}
\simeq
\CH_{\mr I,G}\underset{\mr C^I\times\Bun_G(C)}{\times}Z_{\mr I}.
$$
Thus we have a diagram
$$
Z_{\mr I}  
\stackrel{  h^\la}{\leftarrow} 
\widetilde {{Z}}_{\mr I}
\stackrel{ h^\ra}{\to}  
{Z}_{\mr I} 
$$
in which $ h^\ra$ is the evident projection.

\subsection{Convolution of sheaves}

The geometric Satake correspondence~\cite{MV04}
states that the category $\pshonaffgr$ of $G(\CO)$-equivariant perverse sheaves
on $\affgr_G$ is a tensor category equivalent
to the category of finite-dimensional representations of the dual group $\ch G$.
We shall denote the corresponding functor by
$$V\in \on{Rep}(\ch G)\mapsto \CA^V\in \pshonaffgr.$$
More generally, for a labelling $V_I:I\to \on{Rep}(\ch G)$, we write 
$\CA_G^{V_I}$ for the corresponding object of $\pshonaffgr^{\otimes I}$.

For objects  $\CP\in\catp(Z_{\mr I})$, and $\CA\in\pshonaffgr^{\otimes I}$,
we may form the twisted product
$$
\CP\tboxtimes\CA\in
\catp(\widetilde {{Z}}_{\mr I})
$$
with respect to the projection $ h^\ra:\wt Z_{\mr I}\to Z_{\mr I}$.

We define the functor 
$$
H^I_G :\pshonaffgr^{\otimes I} \times \catp(Z_{\mr I})  \to  \catp(Z_{\mr I}) 
$$
by the formula
$$H^I_G (\CA,\CP) = \bigoplus_{k} \pH^k(h^\la_!(\CP\tboxtimes\CA)).$$
In what follows, we shall only be interested in applying $H^I_G(\cdot,\cdot)$ 
to the intersection cohomology sheaf $'\IC_{\mr I}^0$ of the closure 
$'\ol{Z}^0_{\mr I}\subset {}'Z_{\mr I}$
of the untwisted basic stratum. 
In this way,
we obtain a functor
$$\conv_I: \pshonaffgr^{\otimes I}\to\catp({}'Z_{\mr I}),$$
given by
$$\conv_I(\CA)=H^I_G(\CA,'\IC_{\mr I}^0).$$

When $I$ has a single element,
we write $\conv$ in place of $\conv_I$.

\medskip

Since objects of the category $\pshonaffgr^{\otimes I}$ are direct sums
of intersection cohomology sheaves, by the decomposition theorem, 
every object obtained as $\conv_I(\CA)$ is semi-simple.

The following result will be established in \secref{proof of conv}.

\begin{thm} \label{thmidentifyconv}
For any object $\CA\in\pshonaffgr^{\otimes I}$,
each irreducible summand of
the convolution $\conv_I(\CA)$ is isomorphic to the intersection cohomology sheaf 
of the closure of a connected component 
of an untwisted local stratum $'{Z}_{\mr I}^{\Theta}\subset{} 'Z_{\mr I}$
with constant coefficients.
\end{thm}

\subsubsection{Image category}

Since generic-Hecke modifications commute with meromorphic Hecke
modifications in the natural sense, the convolution functor descends to
a well-defined functor
$$\pshonaffgr^{\otimes I} \times \catp_\CH(Z_{\mr I})  \to  \catp_\CH(Z_{\mr I}).$$

Since the intersection cohomology
sheaf $'\IC_{\mr I}^0$ of the closure 
$'\ol{Z}^0_{\mr I}\subset {}'Z_{\mr I}$
of the untwisted basic stratum is naturally generic-Hecke equivariant,
we obtain a functor
$$
\conv_I: \pshonaffgr^{\otimes I}\to\catp_\CH({}'Z_{\mr I}).
$$

\begin{defn}
We define $\catq(Z_{\mr I})$ to be the full subcategory of $\catp_\CH({}'Z_{\mr I})$
whose objects are isomorphic to direct summands of perverse sheaves
that belong to the image of the above functor.
\end{defn}

When $I$ has a single element, we write $\catq(Z)$ in place of $\catq(Z_{\mr I})$.

\medskip

Putting together Theorem~\ref{thmidentifyconv} and Proposition~\ref{propheq},
we arrive at the following description of the image category.

\begin{thm}
Every irreducible object of $\catq(Z_{\mr I})$ is isomorphic to the intersection cohomology
sheaf of an untwisted local stratum $'{Z}_{\mr I}^{\Theta}\subset{} 'Z_{\mr I}$.
\end{thm}

\subsubsection{Varying the set $I$}

Our present goal is to define an equivalence
$$
\glue_I:\catq(Z)^{\otimes I}\risom\catq(Z_{\mr I}).
$$

Since the categories in questions are semi-simple, it suffices to
specify what this functor does on irreducible objects. We set
$$\glue_I(\otimes_{i\in I} {}'\IC^{\theta_i})={}'\IC^\Theta_{\mr I},\mbox{ where } 
\Theta(i)=\theta_i.$$
Here all of the sheaves have the tautological generic-Hecke equivariant structure.

Along with \thmref{thmidentifyconv}, in \secref{proof of conv}, we will prove
the following:

\begin{cor}\label{corconvglue}
For a finite set $I$, there is a canonical isomorphism
$$
\conv_I\simeq\glue_I\circ\conv^{\otimes I}:\pshonaffgr^{\otimes I}\to\catq(Z_{\mr I}).
$$
\end{cor}

\section{Specialization}      \label{sect specialization}

\subsection{Specialization of affine varieties}

Let $\CB$ be the quotient of the Lie algebra 
of the universal Cartan $T$ by the Lie algebra
of the center $Z(G)$.
Observe that its ring of regular function $\BC[\CB]$ may be represented as a polynomial
algebra with generators $t^\alpha$, 
for simple roots $\alpha\in\ch\Delta_G$, 
and relations $t^\alpha t^\beta=t^{\alpha+\beta}$,
for $\alpha,\beta\in\ch\Delta_G$. 
The quotient torus $T/Z(G)$ 
acts on $\CB$ by linear transformations with a dense orbit isomorphic to $T/Z(G)$.

Let $X$ be an affine variety, acted on by $G$, and consider the decomposition
of the ring $\BC[X]$ into isotypic components.

$$ \BC[X]\simeq \sum_{\lambda\in \domwts_G} \BC[X]_{\lambda}.$$ 
As we have seen in Proposition~\ref{algchar2}, multiplication in $\BC[X]$ respects this grading
if and only if $X$ is horospherical.
In general, multiplication in $\BC[X]$ respects the filtration by the subspaces
$$
\BC[X]_{\leq\mu}=\sum_{\lambda\leq\mu} \BC[X]_{\lambda}.
$$
More precisely, we have
$$
\BC[X]_{\leq\mu}\cdot \BC[X]_{\leq\nu}\subset \BC[X]_{\leq\mu+\nu},
$$
and 
$1\in\BC[X]_{\leq 0}.$

We define the $\BC$-algebra $\Filt[X]$ to be the direct sum 
$$
\Filt[X]=\sum_{\lambda\in\domwts_G} \sum_{\alpha\in\posrts_G}
\BC[X]_{\lambda} t^{\lambda+\alpha},
$$
with multiplication 
$$
f_\lambda t^{\lambda+\alpha} \cdot g_\mu t^{\mu+\beta}= 
f_\lambda g_\mu t^{\lambda+\mu+\alpha+\beta}.
$$
In the isotypic decomposition 
$$
f_\lambda g_\mu=\sum_{\nu\in\domwts_G} h_\nu,
$$ 
a term $h_\nu$ is possibly nonzero only if 
$$\lambda+\mu-\nu\in\posrts_G,
$$
and so only if 
$$\lambda+\mu+\alpha+\beta=\nu+\gamma,
\mbox{ for some $\gamma\in\posrts_G$.}
$$ 
Thus the multiplication in $\Filt[X]$ is well-defined.

The group $G\times T$ naturally acts on $\Filt[X]$, and 
we have a $T$-equivariant inclusion of $\BC$-algebras
$$
i:\BC[\CB]\to \Filt[X]
$$
given by the formula
$$
t^\alpha\mapsto 1t^\alpha, \mbox{ for }\alpha\in\ch\Delta_G.
$$

Let $\CX=\Spec(\Filt[X])$ and 
$\Delta=\Spec(i):\CX\to\CB$.

\begin{thm}[\cite{Pmsb86}]\label{thmdeform}
The projection $\Delta:\CX\to\CB$ is surjective and flat. 
The fiber $\CX_v=\Delta^{-1}(v)$, for regular $v\in \CB$, 
is canonically $G$-isomorphic to $X$,
and the zero fiber $\CX_0=\Delta^{-1}(0)$ is a horospherical $G$-variety.
 \end{thm}

We call the zero fiber $\CX_0$ 
the horospherical variety
associated to $X$, and usually denote it by $X_0$.

\subsection{Basic properties in spherical case}

Now assume that $X$ is an affine spherical $G$-variety.

\begin{prop}\label{openborb}
The fiber $\CX_v=\Delta^{-1}(v)$, for $v\in\CB$, is a spherical $G$-variety.
\end{prop}

\begin{proof}
As a representation of $G$, the ring of regular functions $\BC[\CX_v]$, 
is isomorphic to $\BC[X]$.
Therefore it is simple and so $\CX_v$ is spherical
by Proposition~\ref{algchar1}.
\end{proof}

By construction, the tori $A_v$, associated to $\CX_v$, depend only
on $\BC[\CX_v]$ as $G$-modules. Hence, the family $v\mapsto A_v$
is constant with the fiber $A$, attached to initial spherical $G$-variety $X$.

\medskip

For a Borel subgroup $B^\op\subset G$,
let $\mr\CX^+_v\subset\CX_v$, for $v\in\CB$, be the dense $B^\op$-orbit
provided by Proposition~\ref{openborb}. Evidently, there exists an
open sub-variety $\mr\CX^+\subset\CX$, which specializes to 
$\mr\CX^+_v$ for each $v\in\CB$. 

The parabolics $P^\op_v\supset B^\op$ depend only on the sublattices
$\wts_{\CX_v}\subset \wts_T$. Since the above family of lattices is constant,
so is the family of parabolics.

\begin{prop}\label{pbtriv}
There is a $P^\op$-equivariant isomorphism $\mr\CX^+\simeq\openX^+_0\times\CB$ 
such that the map $\Delta:\mr\CX^+\to\CB$ coincides with the projection onto the 
second factor of the product.
\end{prop}

\begin{proof}
This follows from the fact that the action of $U^\op$ on each $\mr\CX^+_v$ is
free and $\mr\CX^+_v/U^\op\simeq A$.
\end{proof}

Let $\mr\CX_v\subset\CX_v$, for $v\in\CB$, be the dense $G$-orbit
provided by Proposition~\ref{openborb}. Evidently, there exists an
open subvariety $\mr\CX\subset\CX$ that specializes to $\mr\CX_v$
for each $v$.

\begin{cor}\label{smoothfamily}
The open subset $\mr\CX\subset \CX$ is smooth and the map
$\Delta:\mr\CX\to\CB$ is smooth. 
\end{cor}

\begin{proof}
For $ x\in \mr\CX$, we may choose $B^\op$ so that $x\in\mr\CX^+$, and
then apply the previous proposition. 
\end{proof}


\subsection{Family of quasimaps}

We define the family $\CZ_I\to\CB$ of {meromorphic quasimaps}
into the family $\CX\to\CB$
to be the ind-stack 
classifying the data
$$
(v\in\CB,c_I\in C^I, \CP_G\in\Bun_G,\sigma:{C}\setminus {|c_I|}\to 
\CP_G{\gtimes}\CX_v|_{{C}\setminus {|c_I|}})
$$
where $\sigma$ is a section which
factors
$$
\sigma|_{C'}:C'_\test\to 
\CP_G{\gtimes}{\mr\CX_v}|_{C'}\to
\CP_G{\gtimes}{\CX_v}|_{C'},
$$
for some open curve $C'\subset C$.

We write $\CZ_{v,I}$ for the fiber of $\CZ_I$ over a point $v\in\CB$, and
usually denote the zero fiber  by $Z_{0,I}$.

\begin{prop}
The fiber $\CZ_{v,I}$, for regular $v\in \CB$, 
is canonically isomorphic to $Z_I$.
 \end{prop}

\begin{proof}
Immediate from Theorem~\ref{thmdeform}.
\end{proof}

\subsection{Specialization of sheaves}

In this section, we study the nearby cycles in the family $\CZ_I\to\CB$.
The quotient torus $T/Z(G)$ 
acts transitively on the regular elements of $\CB$.
For any line $\ell\subset\CB$ such that $\ell\setminus \{0\}$
consists of regular elements, we have the nearby cycles functor
$$
\psi_I:\catp(Z_{I})\to \catp(Z_{0,I}).
$$
Thanks to the action of $T/Z(G)$, the nearby cycles functors for different lines
are canonically isomorphic. When $I$ consists of one element, we shall write
$\psi_I$ in place of $\psi$.

In \secref{proof of hw spec} we will prove the following:

\begin{prop} \label{hw spec}  \hfill

\smallskip

\noindent{(1)} For $\theta\in \CV(G/S)$ the perverse sheaf $'\IC^\theta_{0}$ appears
in the Jordan-H\"older series of $\Psi({}'\IC^\theta)$ with multiplicity one. 

\smallskip

\noindent{(2)}
If for some $\theta'\in \cowts_A$, the perverse sheaf $'\IC^{\theta'}_{0}$ appears
in the Jordan-H\"older series of $\Psi({}'\IC^\theta)$, then $\theta-\theta'\in \poscowts_X$.

\end{prop}

\subsubsection{Specialization and generic-Hecke equivariance}

The following result will be established in \secref{sect generic}:

\begin{prop}   \label{Levi equiv}
The composition
$$\catp_\CH(Z_{I})\hookrightarrow \catp(Z_{I})
\overset{\psi_I}\to \catp(Z_{0,I})$$
factors canonically through a functor
$$
\catp_\CH(Z_{I})\to \catp_\CH(Z_{0,I}).
$$
\end{prop}

\subsubsection{Specialization and convolution}

Let $\CA$ be an object of $\pshonaffgr^{\otimes I}$
and $\CP\in \catp_\CH(Z_{\mr I})$.

\begin{prop}\label{propconvdeform}
There exists a functorial isomorphism
$$
\psi_I(H^I_G(\CA,\CP)))\simeq 
H^I_G(\CA,\psi_I(\CP))
$$
in the category $\catp_\CH(Z_{0,\mr I})$.
\end{prop}

\begin{proof}

For the family $\CZ_{\mr I}\to\CB$, 
we may form the diagram
$$
\begin{array}{ccccc}
\CZ_{\mr I} & \stackrel{\hl}{\leftarrow} & 
\CH_{I,G}\underset{C^I\times\Bun_G}{\times}\CZ_{\mr I}
& \stackrel{\hr}{\to} & \CZ_{\mr I} \\
\downarrow & & \downarrow & & \downarrow \\
\Bun_G & \stackrel{{}\hl}{\leftarrow} & \CH_{I,G} & \stackrel{{}\hr}{\to} & \Bun_G
\end{array}
$$
in which each square is Cartesian.
We have canonical isomorphisms
$$
\begin{array}{ccll}
\psi_I(H^I_G (\CA,\CP)) & = & \bigoplus_{k} \psi_I(\pH^k(h^\la_!(\CP\tboxtimes\CA)))
& \\
& \simeq & \bigoplus_{k} \pH^k(\psi_I(h^\la_!(\CP\tboxtimes\CA)))
& \mbox{($\psi_I$ is exact in the perverse t-structure)} \\
& \simeq & \bigoplus_{k} \pH^k(h^\la_!(\psi_I(\CP\tboxtimes\CA)))
& \mbox{($\psi_I$ commutes with proper pushforward)} \\
& \simeq & \bigoplus_{k} \pH^k(h^\la_!(\psi_I(\CP)\tboxtimes\CA))
& \mbox{($\psi_I$ commutes with twisted product)} \\
& = & H^I_G(\CA,\psi_I(\CP))
& \\
\end{array}
$$
\end{proof}

\subsubsection{Bad perverse sheaves}

By Corollary~\ref{corconvglue} and Theorem~\ref{thmhoro}, 
the irreducible objects of $\catq(Z_{0,\mr I})$
are isomorphic to the intersection
cohomology sheaves
with constant generic-Hecke equivariant coefficients
of untwisted local strata
$$'Z_{0,\mr I}^\Theta\subset Z_{0,\mr I}, \mbox{ for } \Theta: I\to\cowts_{A_0}.$$

We call an object of $\catp_\CH(Z_{0,I})$ a {\em bad} perverse sheaf if 
each of its simple constituents is not an object of the category $\catq(Z_{0,\mr I})$.
We denote by $\Bad(Z_{0,I})$ the full Serre subcategory 
of $\catp_\CH(Z_{0,I})$
whose objects are bad perverse sheaves. 

Thus, a simple object of $\Bad(Z_{0,I})$ is either the intersection cohomology
of a substack of $Z_{0, I}$ which is {\em not} an untwisted local stratum,
or is the intersection cohomology of an untwisted local stratum with coefficients in 
a {\em nontrivial} irreducible generic-Hecke equivariant local system.

We denote by
$\Quot(Z_{0,I})$ the quotient
of the category $\catp(Z_{0,I})$
by the Serre subcategory $\Bad(Z_{0,I})$.

We denote by $\Subquot(Z_{0,I})$
the full image of 
$\catq(Z_{0,\mr I})$ under the natural projection to $\Quot(Z_{0,I})$.
Since $\catq(Z_{0,\mr I})$ is semisimple, and simple objects go to simple objects under the projection
$\catp_\CH(Z_{0,I})\to \Quot(Z_{0,I})$,
the projection $\catq(Z_{0,\mr I})\to \Subquot(Z_{0,I})$ is clearly
an equivalence.

\medskip

In \secref{proof of special} we will prove the following:

\begin{thm}  \hfill  \label{thmconv}

\smallskip

\noindent{(1)}
In the quotient category $\Quot(Z_{0,I})$, we have an isomorphism
$$\psi_I('\IC^{0}_{\mr I})\simeq
{}'\IC_{0,\mr I}^{0}.$$

\smallskip

\noindent{(2)}
The above isomorphism induces an isomorphism
$$H^I_G(\CA,\psi_I('\IC^{0}_{\mr I}))\simeq H^I_G(\CA,{}'\IC_{0,\mr I}^{0})\in \Quot(Z_{0,I}),$$
functorial in $\CA\in \pshonaffgr^{\otimes I}$

\end{thm}

As a corollary, we obtain:

\begin{cor} \label{corgood}
The image of the functor
$$\catp_\CH(Z_{I})\overset{\psi_I}\to \catp_\CH(Z_{0,I})\to \Quot(Z_{0,I})$$
belongs to the full subcategory $\Subquot(Z_{0,I})$.
\end{cor}

Thus, by Corollary \ref{corgood}, we may define the functor
$$
\Psi_I:\catq(Z_{\mr I})\to\catq(Z_{0,\mr I})
$$
by taking the nearby cycles $\psi_I$, passing to the quotient category,
and then lifting back. When $I$ has a single element, we write $\Psi$
in place of $\Psi_I$.

\medskip

Recall now the functor $\gamma^I:\catq(Z)^{\otimes I}\to \catq(Z_{\mr I})$.
In \secref{proof of special} we will also prove the following:

\begin{prop}\label{pfactfiber}
There is a canonical isomorphism
$$
\Psi_I\circ\glue_I\simeq\glue_I\circ\Psi^{\otimes I}:\catq(Z)^{\otimes I}\to\catq(Z_{0,\mr I}).
$$
such that the following diagram commutes
$$
\begin{array}{ccc}
\Psi_{I}\circ\gamma_I\circ\conv^{\otimes I} & 
\simeq & \gamma_I\circ\Psi^{\otimes I}\circ\conv^{\otimes I} \\
\downarrow & & \downarrow \\
\Psi_{I}\circ\conv_{I}\circ\gamma_{I} & 
 & \gamma_I\circ\conv_0^{\otimes I} \\ 
\downarrow & & \downarrow \\
\conv_{0,I}\circ\gamma_{I} & 
= 
& \conv_{0,I}\circ\gamma_{I} \\ 
\end{array}
$$
where the vertical arrows are previously defined isomorphisms.
\end{prop}

\section{Fusion}   \label{sect Fusion}

\subsection{Universal local acyclicity}   \label{sect ULA}

The constructions in this section will be based on the following general
principle.

Let $\CY$ be a scheme (or stack), mapping to a smooth base scheme
$B$ (in practice we will take $B$ to be $C^I$ for various finite sets $I$).
Let $\CF$ be an object of $\Sh(\CY)$, and let us assume that $\CF$
is universally locaclly acyclic (ULA) over $B$ (we refer the reader to
\cite{BG02} and \cite{Ga04} for a review of the ULA property).

The following properties follow from the assumption:

\begin{lem}  \label{ULA lem}
Let $\CF'$ be a subquotient of $\pH^{n}(\CF)$ for some $n$.

\smallskip

\noindent{(1)} Let $B^0\overset{\jmath} \hookrightarrow B$ be an open subvariety.
Then $\CF'\simeq \jmath_{!*}\jmath^*(\CF')$.

\smallskip

\noindent{(2)}
Let $B'\overset{\imath}\hookrightarrow B$ be a smooth 
locally closed subvariety of codimension $k$. Then 
$\pH^{k'}(\imath^*(\CF'))=0$ for $k'\neq k$, and
$\pH^{k}(\imath^*(\CF'))$ is ULA over $B'$.

\end{lem}

Combining the two assertions of the above lemma, we obtain that
it $\CY'\overset{\imath}\hookrightarrow \CY$ is a smooth locally closed 
subvariety of $\CY$ and $\CY''\overset{\imath'}\hookrightarrow \CY'$
is a smooth closed subvariety of $\CY'$ with the complement 
$\CY^{',0}\overset{\jmath'}\hookrightarrow \CY'$ and $\CF'$ is
as above, then
$$(\imath'\circ \imath)^*(\CF')[-\on{codim}(Y',Y)]\simeq
\imath'{}^*\circ\jmath'_{!*}\bigl(\jmath'{}^*\circ \imath^*(\CF')[-\on{codim}(Y',Y)]\bigr).$$

Finally, let us recall that if $p:\CY\to \CY'$ is a proper map (where $\CY'$ is another
scheme or stack over $B$), then the push-forward $p_!(\CF)\simeq p_*(\CF)\in \Sh(\CY')$
is also ULA over $B$.

\subsection{Fusion for the affine Grassmannian}

Recall that $\affgr_G$ denotes the affine Grassmannian of $G$,
and $\pshonaffgr$ the category of $G(\CO)$-equivariant perverse 
sheaves on $\affgr_G$. We begin by recalling standard results leading to
the fusion product on $\pshonaffgr$.

Recall that $\Aut(\CO)$ denotes the group-scheme of automorphisms of $\CO$.
It naturally acts on $G(\CK)$, $G(\CO)$, and $\affgr_G$.
Let $\catp_{G(\CO)\rtimes\Aut(\CO)}(\affgr_G)$ be the category of
$G(\CO)\rtimes\Aut(\CO)$-equivariant perverse sheaves on $\affgr_G$,
and let $\catp_{\CS}(\affgr_G)$ be the category of perverse sheaves on $\affgr_G$
constructible with respect to the $G(\CO)$-orbits.

\begin{lem}\label{lemforget}
The forgetful functors are equivalences
$$
\catp_{G(\CO)\rtimes\Aut(\CO)}(\affgr_G)\risom\pshonaffgr
\risom\catp_{\CS}(\affgr_G).
$$
\end{lem}

Fix a smooth complex algebraic curve $C$.
For a finite set $I$, let $\bdgr_G$ be the Beilinson-Drinfeld Grassmannian
over $C$. It classifies data
$$
(c_I\in C^I,\CP_G\in\Bun_G(C),\alpha:C\setminus|c_I|\times G\risom\CP_G|_{C\setminus|c_I|})
$$
where $\alpha$ is an isomorphism of $G$-bundles. 
When $I=\{1,\ldots,n\}$, we write $\affgr^{(n)}_G$ in place of $\bdgr_G$.

We have the natural projection
$$
\bdgr_G\to C^I.
$$
For a subvariety $U\subset C^I$, we write $\bdgr_G|_U$ for the fiber product
$$
\bdgr_G|_U=\bdgr_G\underset{C^I}{\times} U
$$

Let $\fs:I\to J$ be a surjection of finite sets.
We have the corresponding closed embedding
$\delta_{\fs}:C^{J}\to C^I$
which in turn induces
a closed embedding
$$
\Delta_{\fs}: \affgr_G^{J}\to \bdgr_G.
$$

\begin{lem}  \hfill   \label{lembdfact}

\smallskip

\noindent{(1)}
The above map indices an isomorphism
$$\affgr_G^{J}\simeq\bdgr_G|_{C^J}.$$

\smallskip

\noindent{(2)}
We have the identification
$$\bdgr_G|_{\mr C^I}\simeq (\onebdgr_G)^I|_{\mr C^I}.$$
\end{lem}

Let $\hat C\to C$ be the $\Aut(\CO)$-torsor of formal parameters.
It classifies data
$$
(c\in C,\tau:D\risom D_c)
$$
where $D$ denotes the formal disk,
$D_c$ denotes the formal neighborhood
of $c\in C$, and $\tau$ is an isomorphism of formal disks.

\begin{lem}\label{lemglobaffgr}
We have the identification
$$
\onebdgr_G\simeq \hat C\overset{\Aut(\CO)}{\times}\affgr_G.
$$
\end{lem}

By Lemmas~\ref{lemforget} and~\ref{lemglobaffgr},
we have the fully faithful functor
$$\rho:\pshonaffgr\to\catp(\onebdgr_G),$$
given by
$$\rho(\CA)=\BC_{C}\wt\boxtimes\CA[1],$$
and, the corresponding fully faithful functor
$$\rho_I:\pshonaffgr^{\otimes I}\to\catp(\bdgr_G|_{\mr C^I}).$$

\medskip

Consider the inclusion 
$$j:\bdgr_G|_{\mr C^I}\to\bdgr_G,$$ 
and the corresponding middle-extension functor
$$j_{!*}:\catp(\bdgr_G|_{\mr C^I})\to\catp(\bdgr_G).$$ 

The following fundamental fact is directly implied by \cite{MV04} (see also \cite{Ga04}).

\begin{prop}  \label{ULA gr}
For any $\CA\in \pshonaffgr^{\otimes I}$ the perverse sheaf
$$j_{!*}(\rho_I(\CA))$$
is ULA with respect to $C^I$.
\end{prop}

For a surjection $\fs:I\to J$ of finite sets, consider the functor
$$\pshonaffgr^{\otimes I}\to \catp(\affgr^J_G),$$
given by
$$\CA\mapsto \Delta_\fs^*( j_{!*}(\rho_{I}(\CA))){[|J|-|I|]}.$$

It is easy to see that for any $\CA\in \pshonaffgr^{\otimes I}$
the resulting object of $\catp(\affgr^J_G)$ belongs to the
image of the fully faithful functor $\rho_J$. Hence, we obtain
a well-defined functor
$$\product_{G,\fs}:\pshonaffgr^{\otimes I}\to\pshonaffgr^{\otimes J},$$
characterized by the property that
$$\rho_J(\product_{G,\fs}(\CA))\simeq\Delta_\fs^*( j_{!*}(\rho_{ I} (\CA))){[|J|-|I|]}.$$

From \lemref{ULA lem} it follows that for surjections $\fs:I\to J$ and $\fr:J\to K$,
there is a canonical isomorphism
$$\product_{\fr\circ\fs}\simeq\product_\fr\circ\product_\fs,$$
providing associativity and commutativity constraints for
the above operation.
This endows $\pshonaffgr$ with a structure of symmetric
monoidal category. 

\subsection{Fusion product on $\catq(Z)$}

We shall now define a (symmetric) monoidal product on $\catq(Z)$ in parallel with the fusion product
of the previous subsection.

For the inclusion 
$$
j:Z_{\mr I}\to Z_I,
$$ 
we have the middle-extension
$$
j_{!*}:\catp_\CH(Z_{\mr I})\to\catp_\CH(Z_{I}).
$$ 

\begin{lem}  
For any object $\CQ_I\in \catq(Z_{\mr I})$, the perverse sheaf
$j_{!*}(\CQ_I)$ is ULA with respect to $C^I$.
\end{lem}

\begin{proof}

This follows from \propref{ULA gr} and the fact that the map $h^\la$ in the definition
of convolution is proper.
\end{proof}

For a surjection $\fs:I\to J$ of finite sets,
we have the corresponding closed embedding $\delta_{\fs}:C^{I'}\to C^I$
which in turn induces
a closed embedding
$$\Delta_\fs: Z_J\to Z_I.$$ 

By the previous lemma and \secref{sect ULA}, we obtain:

\begin{cor}
The functor 
$$
\Delta_\fs^*\circ j_{!*}:\catp_\CH(Z_{\mr I})\to\catp_\CH(Z_{J})
$$  
descends to a functor
$$
\Delta_\fs^*\circ j_{!*}:\catq(Z_{\mr I})\to\catq(Z_{\mr J}).
$$
\end{cor}

Recall that we have an equivalence
$$
\glue_I:\catq(Z)^{\otimes I}\risom\catq(Z_{\mr I}).
$$

We define the {\em fusion product}
$$\product_\fs:\catq(Z)^{\otimes I}\to\catq(Z)^{\otimes J}$$
to be the composite functor
$$\product_\fs(\CP_I)=\glue_{J}^{-1}(\Delta_\fs^*( j_{!*}( \glue_I(\CP_I))))[|J|-|I|].$$

By the preceding discussion, we immediately have the following.

\begin{cor}\label{corconvfusion}
For a surjection $\fs:I\to J$ of finite sets, 
There is a canonical isomorphism
$$
\product_\fs\circ\conv^{\otimes I}\simeq\conv^{\otimes J}\circ\product_{G,\fs}:
\pshonaffgr^{\otimes I}\to\catq(Z)^{\otimes J}.
$$
\end{cor}

When $J$ has one element, so that there is only one possible
surjection $\fs:I\to J$, we write $\product$ in place of $\product_\fs$.

In \secref{proof of hw fuse} we will prove the following:

\begin{prop} \label{hw fuse} \hfill

\smallskip

\noindent{(1)}
For $\theta_1,\theta_2\in\CV(G/S)$, the fusion product
${'\IC^\theta_1}\product{{}'\IC^\theta_2}$ contains
${'\IC^{\theta_1+\theta_2}}$
as a subquotient with multiplicity $1$.

\smallskip

\noindent{(2)}
If for some $\theta_1+\theta_2\neq \theta\in \CV(G/S)$, the perverse sheaf $'\IC^{\theta}$
is a constituent of ${'\IC^\theta_1}\product{{}'\IC^\theta_2}$, then
$\theta-(\theta_1+\theta_2)\notin \poscowts_X$.
\end{prop}

\subsection{Fusion and specialization}

Our goal now is to establish the following compatibility property
between fusion and specialization functors:

\begin{thm}   \label{fib fun}
There is a canonical isomorphism
$$
\product_\fs\circ\Psi^{\otimes I}\simeq\Psi^{\otimes J}\circ\product_{\fs}:
\catq(Z)^{\otimes I}\to\catq(Z_0)^{\otimes J}
$$
such that the following diagram commutes
$$
\begin{array}{ccc}
\Psi^{\otimes J}\circ\product_\fs\circ\conv^{\otimes I} & 
\simeq & \product_\fs\circ\Psi^{\otimes I}\circ\conv^{\otimes I} \\
\downarrow & & \downarrow \\
\Psi^{\otimes J}\circ\conv^{\otimes J}\circ\product_{G,\fs} & 
 & \product_\fs\circ\conv_0^{\otimes I} \\ 
\downarrow & & \downarrow \\
\conv_0^{\otimes J}\circ
\product_{G,\fs} & 
= 
& \conv_0^{\otimes J}\circ\product_{G,\fs} \\ 
\end{array}
$$
where the vertical arrows are previously defined isomorphisms.
\end{thm}

\begin{proof}

After a diagram chase starting from Proposition~\ref{pfactfiber}, 
it remains to exhibit a canonical isomorphism
\begin{equation} \label{fuse and deform}
\Delta_\fs^*\circ j_{!*}\circ\psi_I \simeq \psi_J\circ\Delta_\fs^*\circ j_{!*}
\end{equation}
compatible with convolution.

\begin{lem}
For $\CP_I\in\catq(Z_{\mr I})$, the perverse sheaf
$\psi_I\circ j_{!*}(\CP_I)$ is ULA with respect to the projection $Z_{0,I}\to C^I$.
\end{lem}

\begin{proof}
We have to show that $\psi_I(\CP_I)\in \catp(Z_{0,\mr I})$ can be extended
to an object of $\Sh(Z_{0,I})$, which is ULA with respect to $C^I$. For that
we can replace $\CP_I$ by $H^I_G(\CA,{}'\IC_{\mr I}^0)$
for some $\CA\in \pshonaffgr^I$.

We claim that 
$$\psi_I\biggl(H^I_G(j_{!*}(\CA)), {}'\IC_{I}^0\biggr).$$
provides the required extension. Indeed, as in \propref{propconvdeform},
the above expression is isomorphic
to
$$H^I_G(j_{!*}(\CA),\psi_I({}'\IC_{I}^0)),$$
and the latter is the direct image under the proper map $h^\la$
of the perverse sheaf
$$\psi({}'\IC^0)\tboxtimes j_{!*}(\CA)$$ on
$$\CH_{I,G}\underset{\Bun_G}{\times}\CZ\subset 
\CH_{I,G}\underset{C^I\times\Bun_G}{\times}\CZ_{I},$$
which is ULA over $C^I$ by \propref{ULA gr}.

\end{proof}

From the lemma we obtain that
$$\psi_I\circ j_{!*}(\CP_I)\simeq j_{!*}\circ\psi_I(\CP_I).$$
Therefore,
$$\Delta_\fs^*\circ j_{!*}\circ\psi_I (\CP_I)\simeq
\Delta_\fs^*\circ \psi_I\circ j_{!*}(\CP_I).$$

\medskip

Note that there exists a functorial morphism
$$\Delta_\fs^*\circ \psi_I(\CQ_I) \to  \psi_J\circ\Delta_\fs^*(\CQ_I)$$
for $\CQ_I\in \Sh(Z_I)$. As in the above lemma, we show that this
morphism is an isomorphism when applied to objects of the form 
$j_{!*}(\CP_I)$, $\CP_I\in\catq(Z_{\mr I})$.
Combining, we obtain the isomorphism of \eqref{fuse and deform}.

We leave it to the reader to unwind the isomorphisms to verify the asserted compatibility.
\end{proof}

\section{Tannakian formalism}  \label{sect Tannakian}

\subsection{Geometric Satake equivalence}

Our starting point is the following fundamental result, due to Lusztig, Ginzburg,
Drinfeld, and Mirkovic-Vilonen~\cite{MV04}.

Recall that the category $\pshonaffgr$ of $G(\CO)$-equivariant perverse sheaves
on $\affgr_G$ carries a natural structure of symmetric monoidal (= tensor) category.

\begin{thm}[\cite{MV04}]\label{Satake}
The category $\pshonaffgr$ is equivalent to that of 
finite-dimensional representations of the Langlands dual group $\ch G$.
\end{thm}

This theorem is the basis for the relation between the geometry
of loop spaces of $G$ and representation theory of $\ch G$.

\subsection{The fiber functor}

By the main result of \cite{GNhoro04}, we may identify the category $\catq(Z_0)$
with the category of finite dimensional representations of the torus
$\ch A_0\subset \ch G$. The following will be established in \secref{conv horo}:

\begin{prop}  \hfill  \label{horo cat}

\smallskip

\noindent{(1)}
The equivalence $\catq(Z_0)\simeq \Rep(\ch A_0)$ is naturally equipped with a
monoidal structure, compatible with the commutativity constraints.

\smallskip

\noindent{(2)}
The resulting tensor functor
$$\Rep(\ch G)\simeq\pshonaffgr\overset{\Conv}\to \catq(Z_0)\simeq \Rep(\ch A_0)$$
is naturally isomorphic to the restriction functor under 
$\ch A_0\hookrightarrow \ch T\hookrightarrow \ch G$.
\end{prop}

Composing the functor $\Psi$ with the forgetful functor from $\Rep(\ch A_0)$ to the 
category of vector spaces, we obtain that $\catq(Z)$ is a symmetric monoidal
category, equipped with a tensor functor to the category of finite-dimensional
vector spaces.

\subsection{The associated subgroup}

By \thmref{Satake}, we have a tensor functor
$$\Rep(\ch G)\to \catq(Z),$$
and by 
\thmref{fib fun} and \propref{horo cat}, its composition with
$\catq(Z)\to \on{Vect}$ is the tautological forgetful functor 
on $\Rep(\ch G)$.

By the Tannakian dictionary 
collected in \cite[Section 9]{Nad05}, 
this implies that $\catq(Z)$ is a tensor category equivalent to the category of finite dimensional
representations of a subgroup $\ch H\subset \ch G$ such that $\ch A_0\subset\ch H$.
Under this identification, the convolution $\conv$
corresponds to the restriction of representations from $\ch G$ to $\ch H$, and the nearby cycles
$\Psi$
corresponds to the restriction of representations from $\ch H$ to $\ch A_0$.

\begin{prop} \hfill

\smallskip

\noindent(1)
The subgroup $\ch H\subset \ch G$ is connected and reductive.

\smallskip

\noindent(2)
The torus $\ch A_0\subset \ch H$ is a maximal torus.

\end{prop}

\begin{proof} 

By \propref{hw fuse}, the fusion product $({'\IC^\theta})^{\product n}$ contains
$'\IC^{n\theta}$ as a subquotient. Using \cite[Corollary 2.22]{DM82}, this shows that 
$\ch H$ is connected. The fact that $\ch H$ reductive follows from the
semi-simplicity of $\catq(Z)$.

\medskip

Let us consider again the restriction functor $\Rep(\ch H)\to \Rep(\ch A_0)$
that geometrically translates as the functor $\Psi$. By \propref{hw spec}(1),
for every irreducible object $'\IC^\theta\in \catq(Z)$, the corresponding
object $\Psi({}'\IC^\theta)$ contains a distinguished irreducible constituent, 
namely, $'\IC^\theta_0$. 

By \propref{hw spec}(2) and \propref{hw fuse}(2), the resulting collection 
of lines in vector spaces underlying irreducible representations of $\ch H$
defines a Borel subgroup in $\ch H$, and $\ch A_0$ manifestly identifies
with its Cartan quotient.

\end{proof}

Note that \propref{hw spec} implies that if $'\IC^\theta$ happens to
be an object of $\catq(Z)$, then $\theta\in \cowts_{A^0}\subset \cowts_A$.

\begin{conj}
For every $\theta\in \cowts_{A^0}\cap \CV(G/S)$, the irreducible object
$'\IC^\theta\in \catp_\CH({}'Z)$ belongs to $\catq(Z)$.
\end{conj}

This conjecture would imply
that the Weyl group of $\ch H$ is the same as that associated to $X$
in the theory of spherical varieties~\cite{Br90,Kjams96}.

\newpage


{\Large \part{Proofs--A}}

\vspace*{5mm}

Parts II, III and IV of this paper are devoted to the proofs of various
assertions stated in Part I. Here in Part II, we present proofs that do not
require the use of local models. The structure of this part is as follows:

\medskip

In \secref{sect param}, we prove the statements of \secref{sect stratifications}
regarding $G(\CO)$-orbits on $\CK$-points of spherical varieties. The main
tool will be a distinguished family of partial compactifications
of spherical varieties.

\medskip

In \secref{sect generic}, we prove miscellaneous results concerning
spaces of quasimaps and the generic-Hecke action.

\medskip

In \secref{conv horo}, we consider quasimaps into the particular
target of a horospherical $G$-variety $X$. We describe explicitly 
its stratification, and the behavior of convolution functors in this case.

\vspace*{10mm}

\section{Spherical geometry}  \label{sect param}

\subsection{Structure theory}

\subsubsection{The Levi subgroup}   \label{Levi sgr}

Recall the family $\CX$ of spherical varieties, and recall that the open $P^\op$-stable 
sub-family $\mr\CX^+$ was in fact constant. Choose a point $x\in \openX^+$ and
set 
$$
Q:=\on{Stab}_{P^\op}(x)\simeq S\cap P^\op.
$$ 
We can write
$$
\mr\CX^+\simeq P^\op/M_S\times \CB,
$$
and thus $Q$ is a constant $\CB$-family, contained in the non-constant
family $v\mapsto S_v$, where $S_v\subset G$ denotes the stabilizer of
$x\in \CX_v$, for $v\in\CB$.

By what we have seen, $Q\cap U^\op=\{1\}$ and there is a short exact sequence
$$1\to Q\to P^\op/U^\op\to A\to 1.$$

\medskip

Let us now specialize to the point $0\in \CB$. Let $S_0$ denote the corresponding 
subgroup of $G$ (note that it may be disconnected). By a result of Knop,
the normalizer of $S_0$ is a parabolic subgroup, denoted $P\subset G$, which
is opposite to $P^\op$, i.e., $P\cap P^\op$ is a Levi subgroup of both $P$ and 
$P^\op$. Moreover, $[P,P]=[S_0,S_0]\subset P$.

The above implies, in particular, that $Q\subset S_0$ is a Levi subgroup of $S_0$. 

\subsection{Parametrization of loops}\label{param of orbits}

Our aim here is to describe the equivalence classes of loops in
spherical varieties, a result due to Luna-Vust. We
give an independent proof using their compactification theory~\cite{LV} since it provides
further details we will need in what follows.


\subsubsection{Toroidal compactifications}

We collect here some results from the theory
of compactifications of spherical $G$-varieties. Our basic reference is~\cite{Knop91}.
In what follows, all varieties are assumed to be normal.

Let $S\subset G$ be a spherical subgroup.
A {\em (partial) compactification} of $G/S$ is a $G$-variety $X^c$ with a point $x\in X^c$
such that the stabilizer of $x$ is $S$ and the $G$-orbit through $x$ is dense.
A compactification $X^c$ is said to be {\em toroidal} (or {\em without color}) 
if for any Borel subgroup
$B\subset G$,
no $B$-stable divisor in $G/S$ contains a $G$-orbit of $X^c$ in its closure.
Our first aim is to recall the classification of toroidal compactifications
of $G/S$. There is a more general theory describing all compactifications
but we have no need for it.

Let $\BC(G/S)$ denote the function field of $G/S$.
For a Borel subgroup $B\subset G$,
let $\BC(G/S)^{(B)}$ denote the multiplicative group of nonzero $B$-eigenfunctions in $\BC(G/S)$,
and let $\CQ(G/S)$ denote the lattice dual to the lattice of $B$-characters occuring in $\BC(G/S)^{(B)}$. 
In what follows, let $\CV(G/S)$ be the set of $G$-invariant discrete valuations of $\BC(G/S)$.
(Our aim is to show that it coincides with the set of orbits it previously was defined to parametrize
in \secref{sect stratifications}.)
There is a canonical inclusion 
$$\CV(G/S)\subset \CQ(G/S)$$ 
which takes 
$v\in\CV(G/S)$ to the homomorphism $\chi_f \mapsto v(f)$,
where $\chi_f$ is the $B$-character of
$f\in\BC(G/S)^{(B)}$. (See~\cite{Knop91}, Corollary 1.8.)
It is known that $\CV(G/S)$ forms a simplicial cone of full rank in $\CQ(G/S)$.

We call a subset $\CC\subset \CQ(G/S)$ an {\em allowable strictly convex cone} 
if $\CC$ is a strictly convex cone generated by finitely many elements of $\CV(G/S)$.
We call a nonempty finite set $\CF$
of allowable strictly convex cones in $\CQ(G/S)$ an {\em allowable fan}
if for $\CC\in\CF$, each face of $\CC$ also belongs to $\CF$,
and for each $v\in\CV(G/S)$, there is at most one $\CC\in\CF$ with $v$ in the
interior of $\CC$.

Let $X^c$ be a compactification of $G/S$.
To each $G$-orbit $Y\subset X^c$, we associate a
cone $\CC_Y(X^c)\subset\CQ(G/S)$ as follows.
Let $\mathcal D_Y(X^c)$ be the set of $B$-stable prime divisors $D\subset X^c$ 
such that $Y\subset D$. Each $D\in\mathcal D_Y(X^c)$ defines a valuation
$v_D\in\CV(G/S)$,
and we write $\CB_Y(X^c)$ for the valuations which arise in this way.
(For the construction of $v_D$, 
see~\cite{Knop91}, Lemma 1.4 and the sentence following it.)
Now define the cone $\CC_Y(X^c)\subset\CQ(G/S)$ to be that generated by $\CB_Y(X^c)$.
To the compactification $X^c$, we assign the fan $\CF(X^c)$
which is the union of the cones
$\CC_Y(X^c)$, for all $G$-orbits $Y\subset X^c$.

\begin{thm}[\cite{Knop91}, Theorem 3.3]\label{thmobjs}
The map $X^c\mapsto \CF(X^c)$ induces a bijection between isomorphism
classes of toroidal compactifications of $G/S$ and allowable
strictly convex fans in $\CQ(G/S)$.
\end{thm}

Let $S_1\subset G$ be a second spherical subgroup,
and let $\phi:G/S\to G/S_1$ be a dominant $G$-equivariant morphism.
It is possible to similarly classify the maps between the
compactifications of $G/S$ and those of $G/S_1$.
Such a morphism $\phi$ induces a map $\phi_*:\CQ(G/S)\to\CQ(G/S_1)$
such that $\phi_*(\CV(G/S))\subset\CV(G/S_1)$.
We say that a cone $\CC\subset \CQ(G/S)$ maps under $\phi$ to a cone $\CC'\subset\CQ(G/S_1)$
if we have $\phi_*(\CC)\subset\CC'$, and
a fan $\CF$ maps to a fan $\CF'$ if each cone $\CC\in\CF$ maps
to some cone $\CC'\in\CF'$.

\begin{thm}[\cite{Knop91}, Theorem 4.1]\label{thmmorphs}
Let $X^c$ and $X^c_1$ be toroidal compactifications of $G/S$ and $G/S_1$ respectively.
Then a dominant $G$-equivariant morphism $\phi:G/S\to G/S_1$ extends
to a $G$-equivariant morphism $\phi^c:X^c\to X^c_1$ if and only if $\CF(X^c)$ maps under $\phi$ to $\CF(X^c_1)$.
\end{thm}

Observe that when $\CV(G/S)$ itself is strictly convex, Theorem~\ref{thmobjs} implies
the existence of a canonical compactification corresponding to the fan $\CV(G/S)$.
It is called the {\em wonderful compactification} of $G/S$.
In this case, since the fan $\CV(G/S)$ contains a unique open cone, the wonderful compactification
contains a unique closed $G$-orbit.
In general, 
a compactification $X^c$ is said to be {\em simple} if there is a unique closed $G$-orbit in $X^c$.

\begin{thm}[\cite{BrPau87}, 5.3, Corollaire]\label{thmwonder}
The following are equivalent:
\begin{enumerate}
\item There exists a simple complete toroidal compactification 
of $G/S$.
\item The quotient $N_G(S)/S$ is finite.
\item The cone $\CV(G/S)$ is strictly convex.
\end{enumerate}
\end{thm}

The structure of the normalizer $N_G(S)\subset G$ is well-understood.

\begin{thm}[\cite{BrPau87}, 5.2, Corollaire]\label{thmnormal}
The quotient $N_G(S)/S$ is diagonalizable.
\end{thm}

Our next aim is to recall a fundamental result in the local structure theory of toroidal compactifications.
We recall some of preliminary material.
Fix a Borel subgroup $B^\op\subset G$ such that $B^\op S$ is open in $G$.
Let $P^\op\subset G$ be the parabolic subgroup of all elements $p\in G$
such that $pB^\op S=B^\op S$, and let $U^\op$ be the unipotent radical of $P^\op$.
The complement $G\setminus B^\op S$ is a union of divisors, and so
we may choose a function $f\in \BC[G]$ such that $G\setminus B^\op S$ 
is the set-theoretic zero locus of $f$.
The differential $df$ at the identity $1\in G$ defines an element
in the coadjoint representation of $G$.
The centralizer $M\subset G$ of this element is a Levi factor of $P^\op$. 
The quotient $M/M\cap S$ is a torus which we denote by $A$.

Via the embedding $A\subset G/S$, we obtain a map $\BC(G/S)^{(B^\op)}\to\BC(A)$, 
defined by
$\chi_f\mapsto f|_{A}$, which
induces
an isomorphism
$$
\CQ(G/S)\simeq \cowts_{A}.
$$
We write 
$
\domcowts_{G/S}\subset\cowts_A
$ 
for the image of $\CV(G/S)\subset\CQ(G/S)$ 
under this identification. By abuse of notation, we sometimes write $\domcowts_A$ in place
of $\domcowts_{G/S}$. When $S$ is the stabilizer of a point $x$ in the dense $G$-orbit $\openX\subset X$
of a spherical $G$-variety $X$, so that we have $\openX\simeq G/S$,
we also sometimes write $\domcowts_X$ in place of $\domcowts_{G/S}$.

Now for the moment, assume that the quotient $N_G(S)/S$ is finite.
Let $X^c$ be the wonderful compactification of $G/S$, and let $x\in X^c$
be a point with stabilizer $S$. 
Identify the $M$-orbit through $x\in X^c$
with the torus $A$.
Let $Y^+\subset X^c$ be the toric compactification of $A$
characterized by the property that for $\lambda\in\cowts_{A}$, we have 
$$\lim_{t\to 0} \lambda(t)\in Y^+
\mbox{ if and only if } \lambda\in\domcowts_{G/S}.
$$ 
In the classification of toric varieties, $Y^+$ corresponds to the cone $\domcowts_{G/S}$.

\begin{thm}[\cite{BLV86}, Th\'eor\`eme 3.5]\label{thmlocal}
The action map $U^\op\times Y^+\to X^c$ is an open embedding
and its image contains an open nonempty subset of each $G$-orbit in $X^c$.
\end{thm}

Now for an arbitrary spherical subgroup $S\subset G$ for which
the quotient $N_G(S)/S$ is not necessarily finite, consider
the spherical subgroup $S_1\subset G$ generated by $S$ and
the connected component of the normalizer $N_G(S)$. 
The natural morphism $\phi:G/S\to G/S_1$ induces
a map $\phi_*:\cowts_{A}\to\cowts_{A_1}$
such that $\phi_*(\domcowts_{G/S})\subset\domcowts_{G/S_1}$.
The following is easy to check.

\begin{lem}
The cone $\domcowts_{G/S_1}\subset \cowts_{A_1}$ is strictly convex, and the cone
$\domcowts_{G/S}\subset \cowts_{A}$ is its inverse image under $\phi_*$.
\end{lem}

By the first assertion of the lemma and 
Theorem~\ref{thmwonder}, $N_G(S_1)/S_1$ is finite, and $G/S_1$
admits a wonderful compactification $X^c_1$. Let $X^c$ be any toroidal
compactification of $G/S$ such that $\phi$ extends to a $G$-equivariant morphism
$\phi^c:X^c\to X^c_1$.
Let $Y^+_1\subset X^c_1$ be the toric compactification of $A_1$
constructed above, and let $Y\subset X^c$ be its inverse image under $\phi^c$.
Applying Theorem~\ref{thmlocal} to the pair $Y^+_1\subset X^c_1$, 
we conclude that the assertion of Theorem~\ref{thmlocal}
holds equally well
for the pair $Y^+\subset X^c$.


\subsubsection{Statement of parametrization}

Recall from the previous section
that we may identify $\CV(G/S)$ with a subset
$
\domcowts_{G/S}\subset\cowts_A.
$
Via the embedding $A\subset G/S$,
we may consider $\cowts_{A}$, and so also $\domcowts_{G/S}$,
as a subset of $(G/S)(\CK)$.

\begin{thm}\label{thmparam}
Each element of $(G/S)(\CK)$ contains a unique element of $\domcowts_{G/S}$
in the $G(\CO)$-orbit through it.
\end{thm}

The following well-known cases of the theorem are worth pointing out.

\begin{example}
Let $G$ be the product $H\times H$, and 
let $S$ be the diagonal copy of $H$.
The theorem gives the Cartan decomposition
$$
H(\CO)\backslash H(\CK)/H(\CO)\risom\domcowts_{H}
$$
where $\domcowts_H$ is the semigroup of dominant coweights of $H$.
\end{example}

\begin{example}
Let $S$ be the unipotent radical $U$ of a Borel subgroup of $G$.
The theorem gives the Iwasawa decomposition
$$
G(\CO)\backslash G(\CK)/U(\CK)\risom\cowts_G
$$
where $\cowts_G$ is the lattice of coweights of $G$.
\end{example}

\noindent
{\it Proof of Theorem~\ref{thmparam}.}
The proof 
is a simple application of the theory of compactifications of spherical
varieties discussed in the previous section. 

Let $S_1\subset G$
be the spherical subgroup generated by $S$ and
the connected component of the normalizer $N_G(S)$. 
Let $X^c_1$ be the wonderful compactification of $G/S_1$,
and let $X^c$ be any complete toroidal compactification of $G/S$
which maps to $X^c_1$. 

Via the embedding $G/S\to X^c$, we obtain an injection $(G/S)(\CK)\to X^c(\CK)$.
Since $X^c$ is compact, each element $\gamma\in X^c(\CK)$ extends
to an element $\bar\gamma\in X^c(\CO)$. By the slight generalization of Theorem~\ref{thmlocal}
discussed at the end of the previous section, 
the $G(\CO)$-orbit through $\bar\gamma$ contains an element $\bar\gamma'\in X^c(\CO)$
which lies in the image of the action map $U^\op\times Y^+\to X^c$.
Therefore we may consider $\bar\gamma'$ as an element of 
$U^\op(\CO)\times Y^+(\CO)$.
Thus acting by an element of $U^\op(\CO)$,
we see that the $G(\CO)$-orbit through $\gamma$ contains an
element in $Y^+(\CO)$. 
Acting by an element of $T(\CO)$,
we conclude that the $G(\CO)$-orbit through $\gamma$ contains an
element $\lambda\in\domcowts_{G/S}$.

To check that two distinct elements $\lambda,\lambda'\in\domcowts_{G/S}$
are not in the same $G(\CO)$-orbit,  it suffices to check that their
interactions with the divisor at infinity of $X^c$
may be distinguished. If neither $\lambda$ or $\lambda'$ is a $\BZ_{>0}$-multiple of the other,
then by Theorems~\ref{thmobjs} and~\ref{thmmorphs}, we may choose
the compactification $X^c$ to have the property that $\lim_{t\to0}\lambda(t)$
and $\lim_{t\to0}\lambda'(t)$ do not lie in the same $G$-orbit.
If either $\lambda$ or $\lambda'$ is a $\BZ_{>0}$-multiple of the other, then we may choose
the compactification $X^c$ so that $\lim_{t\to0}\lambda(t)=\lim_{t\to0}\lambda'(t)$ 
lies in a $G$-orbit of codimension $1$. It is easy to check that in this
case, the orders of intersection of the closures of $\lambda(t)$ and $\lambda'(t)$ 
with this $G$-orbit are distinct.
\newline


\subsection{Applications}\label{proof of orbits}

First, observe that \thmref{thmparamval} follows immediately from \thmref{thmparam}.

To establish \thmref{thm count orbits},
it remains to show that for a subgroup
$S\subset G$,
if $S(\CK)$
acts on $\affgr_G$ with countably many orbits,
then $S\subset G$ is spherical.

Any one parameter subgroup $\lambda:\BC^\times\to G$ defines a point of
$\affgr_G$ which we also denote by $\lambda$. It is easy to see that 
the $G$-orbit $F_\lambda\subset\affgr$
through $\lambda$
is a flag variety of $G$.
If $S(\CK)$
acts on $\affgr_G$ with countably many orbits,
then the number of orbits intersecting 
$F_\lambda \subset \affgr_G$ is countable. 
Therefore one of the orbits intersects $F_\lambda$ in an open set.
Let $\mu\in F_\lambda$ be a point in this open set.

We may identify the tangent space of $F_\lambda$ at the point $\mu$
with the quotient $\fg/\fp_\mu$, where $\fg$ is the Lie algebra of $G$,
and $\fp_\mu$ is the Lie algebra of the parabolic subgroup $P_\mu\subset G$
which stabilizes $\mu$.
In order for an $S(\CK)$-orbit in $\affgr_G$ to intersect $F_\lambda$
in an open set containing $\mu$, 
we must have that the Lie algebra $\fs$ of $S$ surjects onto $\fg/\fp_\mu$.
Choosing $\lambda$ regular, so that $\mu$ is regular as well,
we conclude that $\fs$ must surject onto the quotient of $\fg$
by a Borel subalgebra. This implies that $S$ has an open orbit
in the flag variety of $G$, and so it is a spherical subgroup.

\medskip

The following will be used in Sect.~\ref{sect trans locus}.

\begin{prop}\label{local trans}
For $\lambda\in \openX^+(\CK)$, and $u\in U^{op}(\CK)$, suppose $u\cdot\lambda$ is in a
$G(\CO)$-orbit in $X(\CK)$ in the closure of the orbit through $\lambda$. Then $u\in U^{op}(\CO)$,
and hence $u\cdot\lambda$ is in the $G(\CO)$-orbit through $\lambda$.
\end{prop}

\begin{proof}
Consider the loops $\lambda$ and $u\cdot\lambda$ as elements of $X^c(\CK)$.
If $u\cdot\lambda$ lies in $X^+(\CO)$ then we are done by similar arguments as in the proof
of Theorem~\ref{thmparam} and the fact that $U^{op}\cap S$ is the identity.

Suppose $u\cdot\lambda $ lies in $X^+(\CK) - X^+(\CO)$, and let $\eta\in\Lambda_A$
denote the projection $p(u\cdot\lambda)$. Then we may find
a conjugate open set $X_1^+$ such that $u\cdot\lambda $ lies in $X_1^+(\CO)$.
We also have the map $p_1:X_1^+\to Y^+$, and the projection 
$\eta_1=p_1(u\cdot\lambda)\in\Lambda_A$. 
It follows from the fact that $u\cdot\lambda\in (X^+(\CK) \cap X_1^+(\CK))- X^+(\CO)$
that $\eta_1$ must have a deeper pole than $\eta$. But then $u\cdot\lambda$ is in a deeper
stratum of $\openX(\CK)$ so can not be in the closure of that containing $\lambda$.
\end{proof}


\section{Generic-Hecke action}   \label{sect generic}

\subsection{Proof of \lemref{prime closed}}

Let $p:C_\test\to\test$ be the projection,
and let $\CP'_{\pi_0(S)}$ be the
generic $\pi_0(S)$-bundle of an $\test$-valued quasimap.
We first claim that the quasimap is untwisted if and only if
$\CP'_{\pi_0(S)}$ 
arises as the pullback via $p$ of a $\pi_0(S)$-bundle on $\test$.
Clearly, the latter condition is sufficient.
To see it is necessary, choose a faithful representation $V$ of $\pi_0(S)$, and consider
the associated local system 
$$
\CL'=\CP'_{\pi_0(S)}\overset{\pi_0(S)}{\times} V.
$$
By assumption, $\CL'$ is trivial on the subschemes  $\{s\}\times C\subset C_\test$
for every geometric point $s\in \test$. Thus it extends to a local system $\CL$
on all of $C_\test$,
and, since $C$ is connected, the adjunction morphism
$$
p^* R^0p_*\CL\to\CL
$$
is an isomorphism. We conclude that 
$\CL$, and thus $\CL'$ as well, is the pullback of the local system $R^0p_*\CL$ and the claim is proved.

Now suppose $\test$ is a dense subscheme of $\ol\test$, and we have a $\ol\test$-valued
quasimap whose restriction to $\test$ is untwisted. Let $\ol\CP_{\pi_0(S)}$ be the generic 
$\pi_0(S)$-bundle of the quasimap, and let $\CP_{\pi_0(S)}$ be that of its restriction.
By the above discussion, there is a 
$\pi_0(S)$-local system $\CR_{\pi_0(S)}$ on $\test$ such that
$$
p^*\CR_{\pi_0(S)}\simeq\CP_{\pi_0(S)}.
$$
If $\CR_{\pi_0(S)}$ did not extend to the complement $\ol\test\setminus \test$, 
then $\CP_{\pi_0(S)}$ would not extend to an open subset of $p^{-1}(\ol\test\setminus \test)$.
Since this is a contradiction to the existence of $\ol\CP_{\pi_0(S)}$, we conclude
that $\CR_{\pi_0(S)}$ does indeed extend to a $\pi_0(S)$-local system 
$\ol\CR_{\pi_0(S)}$ on all of $\ol\test$, and the above isomorphism extends
to an isomorphism
$$
p^*\ol\CR_{\pi_0(S)}\simeq\ol\CP_{\pi_0(S)}.
$$
Thus the $\ol\test$-valued quasimap is also untwisted and the lemma is proved.

\subsection{Proof of \propref{prophecketrans}}

For $k=1,2$, fix closed points $z_k\in{}'Z_{c_I}^{(\fp,\Theta)}$. We may think of $z_k$, for $k=1,2$,
in terms of data $(c_I,\CP_G^k,\CP_S^k)$,
where $\CP_G^k$ is a $G$-bundle on $C$,
and $\CP_S^k$ is an $S$-bundle on $C\setminus |c_I|$
equipped with an $S$-equivariant bundle map
$$
\CP_S^k\to\CP^k_G|_{C\setminus |c_I|}.
$$
The assertion of the proposition is that on some open curve $C'\subset C$
containing $|c_I|$, we have an isomorphism
$$
\CP_G^1\simeq\CP_G^2
$$
which restricts on $C'\setminus |c_I|$ to give an isomorphism
$$
\CP_S^1\simeq\CP_S^2.
$$

By choosing an appropriate
open curve $C_1\subset C$
containing $|c_I|$, we may assume that $\CP^2_G$ is the trivial bundle
$$
\CP^0_G=C_1\times G.
$$
Then to prove the proposition, it suffices to find a trivialization of $\CP^1_G$
on some open curve $C'\subset C_1$ containing $|c_I|$ such that the induced
isomorphism 
$$
\CP^1_G|_{C'}
\simeq
\CP^0_G|_{C'}
$$
restricts to give an isomorphism
$$
\CP_S^1|_{C'\setminus |c_I|}\simeq\CP_S^2|_{C'\setminus |c_I|}.
$$

Since the $\pi_0(S)$-bundles induced from $\CP_S^k$, for $k=1,2$, are trivial,
by choosing an appropriate
open curve $C_2\subset C_1\setminus |c_I|$, we may assume 
that we have an isomorphism
$$
\alpha:\CP_G^1|_{C_2}\risom\CP^0_G|_{C_2}
$$
which restricts to give an isomorphism
$$
\CP_S^1|_{C_2}\risom\CP_S^2|_{C_2}.
$$
Then via the standard trivialization of $\CP^0_G$ and the isomorphism $\alpha$, 
we obtain a trivialization $\tau$ of the restriction 
$
\CP_G^1|_{C_2}.
$
If $\tau$ extends across the points $|c_I|$, then we are done.
Otherwise, we are left to try to change $\tau$ so that it extends,
but so that we do not change the bundle $\CP_S^1|_{C_2}$.
More precisely, to prove the proposition, it suffices to find an open curve
$C_3\subset C_2$, and a map 
$$\Gamma:C_3\to S,$$ 
so that multiplying by $\Gamma$,
we obtain a trivialization 
$$
\Gamma\cdot\tau
\quad
\mbox{ of }
\quad
\CP_G^1|_{C_2}
$$
which does extend across the points $|c_I|$.

Now by definition,
the assumption that $z_k\in {}'Z_{c_I}^{(\fp,\Theta)}$, for $k=1,2$, implies that 
we may prove the analogue of the proposition on the formal punctured neighborhood 
$D^\times_{|c_I|}$.
In other words, there is a map 
$$
\gamma:D^\times_{|c_I|}\to S^0
$$ 
such that the trivialization
$$
\gamma\cdot(\tau|_{D^\times_{|c_I|}})
\quad
\mbox{ of }
\quad
\CP^1_G|_{D^\times_{|c_I|}}
$$
extends to a trivialization on the formal 
neighborhood $D_{|c_I|}$. 
Taking the class of $\gamma$
in the product of affine Grassmannians of $S^0$ at the points $|c_I|$, 
we may find 
an open curve $C_3\subset C_2$
and a map $\Gamma:C_3\to S^0$ such that we have an equality of classes
$$
[\Gamma]=[\gamma]
$$
in the product of affine Grassmannians.
In other words, we have an equality of maps
$$
\Gamma|_{D^\times_{|c_I|}}=\gamma\cdot\gamma_+,
$$
for some $\gamma_+\in S^0(\CO_{|c_I|})$.
For any congruence subgroup $S^0_{++}\subset S^0(\CO_{|c_I|})$,
we may find an open curve $C_4\subset C_3\cup{|c_I|}$ containing $|c_I|$
and a map $\Gamma_+:C_4\to S^0$ such that we have an equality
$$
[\Gamma_+|_{D_{|c_I|}}]=[\gamma_+]
$$
in the quotient group $S^0(\CO_{|c_I|})/S^0_{++}$.
In particular, we may take $S^0_{++}$ to be contained in the
stabilizer of the class of $\tau$ in the product of affine Grassmannians of $G$.
Thus we conclude that the map
$$
\Gamma\cdot\Gamma_+^{-1}:C_4\to S^0
$$
takes $\tau$ to a trivialization which extends across $|c_I|$.
This completes the proof of the proposition.

\subsection{Proof of \propref{propheq}}

Recall that a Hecke equivariant structure is determined
by its values on substacks of the 
ind-stack 
$\H_{Z_I,(1)}$ of
generic-Hecke modifications
 at a single point.
We may realize $\H_{Z_I,(1)}$ as the twisted product
of an open subset of $Z_I\times C$ with the affine Grassmannian $\affgr_{S^0}$.

\begin{lem}
Suppose smooth generic-Hecke modifications 
$Y_1,Y_2\subset \H_{Z_I,(1)}$
both lie in the twisted product of an open subset $U\subset Z_I\times C$ with
a single component of the affine Grassmannian $\affgr_{S^0}$.
Then there is a sequence of nonempty smooth generic-Hecke modifications
and morphisms of Hecke modifications
$$
Y_1 \leftarrow W_{1} \to W_2 \leftarrow \cdots \rightarrow W_{k-1} \leftarrow W_{k}\to Y_2.
$$
\end{lem}

\begin{proof}
We may assume $Y_2$ is the trivial modification.

Let $S^0_r$ be the maximal reductive quotient of $S^0$.
Consider the projection of affine Grassmannians $\affgr_{S^0}\to\affgr_{S^0_r}$,
and the induced projection on their twisted products with  $U$.
Let $W_r\subset \affgr_{S^0_r}$ be the largest $S^0_r(\CO)$-orbit such that its twisted product
$U\tilde\times W_r$ intersects the projection of $Y_1$.
In particular, the intersection is a nonempty open subset of the projection of $Y_1$.
We may truncate the inverse image of $W_r$ to obtain a smooth $S^0(\CO)$-invariant
subset $W\subset \affgr_{S^0}$ so that its twisted product $U\tilde\times W$ intersects $Y_1$ in a nonempty open set.
Thus
we have a diagram of nonempty smooth generic-Hecke modifications 
$$
Y_1\leftarrow Y_1\cap (U\tilde \times W) \rightarrow U\tilde\times W.
$$

Let $U_{S^0}$ be a maximal unipotent subgroup of $S^0$.
We may truncate the $U_{S^0}(\CK)$-orbit through the base point $U_{S^0}(\CK)\cdot [S^0(\CO)]\subset\affgr_{S^0}$
to obtain a smooth subset $V\subset\affgr_{S^0}$ such that
the intersection $W\cap V$ is nonempty. Let $(W\cap V)^{sm}$
denote the smooth part of the intersection.
Let $\hat U$ denote the open set $U$ equipped with large level structure
at the modification point.
Taking the product with $\hat U$,
we obtain a diagram of nonempty smooth generic-Hecke modifications
$$
U\tilde\times W\leftarrow \hat U \times (W\cap V)^{sm}\to \hat U \times V.
$$
Finally, taking $Y_2$ to be the trivial Hecke modification corresponding to the base point 
$[S^0(\CO)]\subset\affgr_{S^0}$, we have 
a diagram of nonempty smooth generic-Hecke modifications
$$
\hat U \times V \leftarrow \hat U \times \{[S^0(\CO)]\} \to Y_2.
$$
\end{proof}

By the lemma, it suffices to show that if the modifications lying in a component
of $\affgr_{S^0}$ preserve a component of  $'Z_I^{(\fp,\Theta)}$,
then the component of $\affgr_{S^0}$ is the one containing the trivial modification.
To see this, first observe that it is clearly true for the basic stratum  $'Z_I^{(\fp,0)}$
since it is isomorphic to the product $C^{(\fp)}\times\Bun_{S^0}$.
Here we write $C^{(\fp)}\subset C^I$ for the locally closed subvariety of points
$c_I\in C^I$ such that the coincidences among the points $|c_I|\subset C$
are given by the partition $\fp$.
In general, for any pair of points $z_1,z_2\in {}'Z_{I}^{(\fp,\Theta)}$
which are related by a generic-Hecke modification given by a connected
subscheme of $\affgr_{S^0}$, we may simultaneously
modify their associated $G$-bundles
at their pole points so that we obtain points $z^0_1,z^0_2\in {}'Z_{I}^{(\fp,0)}$
which are still related by the generic-Hecke modification given by the same connected
subscheme of $\affgr_{S^0}$. If the points $z_1,z_2\in {}'Z_{I}^{(\fp,\Theta)}$
are in the same connected component of $'Z_{I}^{(\fp,\Theta)}$,
then the same is true for the points $z^0_1,z^0_2\in {}'Z_{I}^{(\fp,0)}$.
Thus we are done since we have already seen the assertion is true 
for ${}'Z_{I}^{(\fp,0)}$.

\subsection{Generic-Hecke Levi equivariance}

We define the ind-stack
$\CH_{Z_I,(n)}^{(m)} $
of {generic-Hecke modifications with level structure} to be that classifying data
$$
(c_I, \CP^1_G,\CP^2_G,\sigma_1,\sigma_2;
c_{(n)}, \alpha; \beta_1,\beta_2)
$$
where 
$(c_I, \CP^1_G,\CP^2_G,\sigma_1,\sigma_2; c_{(n)}, \alpha)\in \CH_{Z_I,(n)},$
and 
$\beta_i$ is an isomophism of $G$-bundles
$$
\beta_i:\|m\cdot c_{(n)}\|\times G\risom \CP^i_G|_{\|m\cdot c_{(n)}\|}
$$
such that for $1\in G$, the following induced diagram commutes
$$
\begin{array}{ccc}
\|m\cdot c_{(n)}\|
& \stackrel{\beta_i}{\to}  &
\CP^i_G\overset{G}{\times} \openX|_{\|m\cdot c_{(n)}\|}\\
\downarrow &  & \downarrow \\
\|m\cdot c_{(n)}\|
& \stackrel{\sigma_i}{\to} &
\CP^i_G\overset{G}{\times} \openX|_{\|m\cdot c_{(n)}\|}
\end{array}
$$
where the vertical maps are the identity.

For sufficiently small modifications $\alpha$ and large order $m$, it makes sense
to ask for $\alpha$ to come from a subgroup of $S$. We take the subgroup in
question to be $Q$ (see \secref{Levi sgr}).

In particular, we may define the ind-substack 
$$\CH_{Z_I, (n),Q}^{(m)}\subset\CH_{Z_I, (n)}^{(m)}$$
of  {\em generic Levi modifications with level structure} to consist of those
modifications which come from $Q$.

We call a smooth generic-Hecke correspondence $Y$ a {\em Levi correspondence}
if its defining map factorizes
$$
Y\to \CH_{Z_I,(n),Q}^{(m)}\to \CH_{Z_I,(n)}
$$ 
for some $m$.
In analogy with the notion of Hecke equivariant perverse sheaf,
we have
the abelian category $\catp_{\CH_Q}(Z_I)$
of Levi equivariant perverse sheaves on $Z_I$.
We write $\catp_{\CH_Q,c}(Z_I)$ for the full subcategory of $\catp_{\CH_Q}(Z_I)$
whose underlying objects are constructible along the orbits of the
generic-Hecke modifications.

\begin{prop}\label{plevieq}
For $X$ horospherical,
the forgetful functor 
$$
\catp_\CH(Z_I)\to\catp_{\CH_Q,c}(Z_I)
$$ 
is an equivalence.
\end{prop}

\begin{proof}
Let $U(S)$ be the unipotent radical of $S$ so that we have $S\simeq U\rtimes Q$.
By the proof of \propref{propheq}, we see that a generic-Hecke equivariant structure
is determined by its values on Levi correspondences. In other words, the forgetful functor
$$
\catp_\CH(Z_I)\to\catp_{\CH_Q}(Z_I)
$$
is a fully faithful embedding. Thus to prove the proposition, we must check that every object of
$\catp_{\CH_Q,c}(Z_I)$ can be equipped with a generic-Hecke equivariant structure.
As explained in the proof of \propref{propheq}, it suffices to consider Hecke modifications
which are twisted products with subschemes of the affine Grassmannian $\affgr_{S^0}$.
Observe that
$\affgr_{S^0}$ is exhausted by 
the $U(S)(\CK)$-orbits through $\affgr_Q$, and each orbit may be written
as the increasing union of smooth affine spaces. Thus there is no obstruction to 
lifting the Hecke equivariant structure from the Levi subgroup $Q$.
\end{proof}


\section{Convolution action in the horospherical case}   \label{conv horo}

Throughout this section, we will only consider $X$ horospherical.

\subsection{Stratification in the horospherical case}    \label{horo strata}

When $X$ is horospherical,
we will need a complete
stratification of $Z_I$, not only what we called local strata.
We provide the definitions here and refer the reader to~\cite{GNhoro04} for more details.

For a positive coweight $\theta^\pos\in\cowts^\pos_A$,
we write $\fU(\theta^\pos)$ for 
a decomposition $\theta^\pos=\sum_m n_m\theta^\pos_m$,
for $\theta^\pos_m\in\poscowts_X\setminus\{0\}$ distinct, and $n_m$ positive integers.
We say that a $\poscowts_X$-valued divisor on $C$ is of type $\fU(\theta^\pos)$
if it is of the form $\sum_{m}\sum_{n=1}^{m_n} \theta^\pos_{m} \cdot c_{m,n}$, 
for $c_{m,n}\in C$ distinct.
For a partition $\fp$ of the set $I$, a labelling $\Theta:\fp\to\cowts_{A}$,
and a decomposition $\fU(\theta^\pos)$ of a positive coweight
$\theta^\pos\in \cowts^\pos_A$, 
we say that a quasimap
$
(c_I,\CP_G,\sigma)
\in Z_I
$
is of type $(\fp,\Theta,\fU(\theta^\pos))$ if 
the coincidences among the pole points $|c_I|\subset C$ are given by the partition $\fp$,
and
the $\domcowts_{A}$-valued divisor on $C$
associated to the quasimap 
is equal to $\Theta\cdot c_I+
\sum_{m}\sum_{n=1}^{n_m} \theta^\pos_{m} \cdot c_{m,n},$
for $c_{m,n}\in C$ distinct and disjoint from $|c_I|\subset C$.
We define the {stratum}
$$
Z_I^{(\fp,\Theta,\fU(\theta^\pos))}\subset {}Z_{I}
$$
to consist of those quasimaps of type $(\fp,\Theta,\fU(\theta^\pos))$.
When $I$ is empty, the partition $\fp$ and labelling $\Theta$ are vacuous, 
and we write $Z^{\fU(\theta^\pos)}_\emptyset$
in place of $Z^{(\fp,\Theta,\fU(\theta^\pos))}_\emptyset$.

\subsection{Adding an auxiliary $A$-bundle}
In many arguments, we will need the following generalization of $Z_I$.
Define the ind-stack $\stZ_I$ to be that classifying data
$$
(c_I\in C^I,\CP_G\in\Bun_G,\CP_{A}\in\Bun_A,\sigma: \CP_{A}|_{C\setminus |c_I|}\to 
\CP_G{\gtimes}X|_{C\setminus |c_I|})
$$
where 
$\sigma$ is an $A$-equivariant map (here we are using the fact that $A$ acts on $X$),
which factors
$$
\sigma|_{C'}:\CP_{A}|_{C'}\to 
\CP_G{\gtimes}{\openX}|_{C'}\to
\CP_G{\gtimes}{X}|_{C'}
$$
for some open curve $C'\subset C\setminus |c_I|$.
We write $\fr: Z_I \to  \stZ_I$ for the obvious induction map, and have 
a Cartesian diagram
$$
\begin{array}{ccc}
Z_I & \stackrel{\fr}{\to} & \stZ_I \\
\downarrow & & \downarrow \\
\Bun_{\langle1\rangle} & \to & \Bun_{A}. \\
\end{array}
$$

\subsubsection{Description of the strata}
We stratify $\stZ_I$ in the same way
that we stratified $Z_I$.
For data $(\fp,\Theta,\fU(\theta^\pos))$, we have the corresponding stratum
$\stZ_I^{\fp,\Theta,\fU(\theta^\pos)}$, and a Cartesian diagram
$$
\begin{array}{ccc}
Z_I^{(\fp,\Theta,\fU(\theta^\pos))} & \stackrel{\fr}{\to} & \stZ_I^{(\fp,\Theta,\fU(\theta^\pos))} \\
\downarrow & & \downarrow \\
\Bun_{\langle1\rangle} & \to & \Bun_{A}. \\
\end{array}
$$


We can describe the stratum $\stZ_I^{(\fp,\Theta,\fU(\theta^\pos)}$ more explicitly
in terms of the stack $\Bun_P$ as follows.

For a partition $\fU(\theta^\pos)$ as above, consider the corresponding partially
symmetrized power of the curve $C^{\fU(\theta^\pos)}=\prod_m C^{(n_m)}$. 
We will write 
$\prod_m\prod_{n=1}^{n_m} c_{m,n}$
for an element of $C^{\fU(\theta^\pos)}$.
Let $\mr C^{\fU(\theta^\pos)}\subset C^{\fU(\theta^\pos)}$ be the complement
to the diagonal divisor.
Similarly, for a partition $\fp$ of the set $I$ into $k$ elements let $C^\fp$
denote $C^k$, and let $\mr C^\fp\subset C^\fp$ be the complement to the diagonal
divisor.
Finally, let  $C^{(\fp,\fU(\theta^\pos))}\subset \mr C^\fp\times \mr C^{\fU(\theta^\pos)}$
be 
the complement to the diagonal divisor.

Observe that there a canonical {\it finite} map
$$\overline{j}_{(\fp,\Theta,\fU(\theta^\pos))}: C^\fp\times 
C^{\fU(\theta^\pos)} \times \stZ_{\emptyset}\to \stZ_I,$$
given by 
$$\overline{j}_{(\fp,\Theta,\fU(\theta^\pos))}(c_I,\prod_m\prod_{n=1}^{n_m} c_{m,n},(\CP_G,\CP_A,\sigma))=
(\CP_G,\CP_A(-\Theta\cdot c_I-\sum_{m}\sum_{n=1}^{n_m}\theta^\pos_m\cdot c_{m,n}),\sigma).$$
We denote by $j_{(\fp,\Theta,\fU(\theta^\pos))}$ the composition of 
$\overline{j}_{(\fp,\Theta,\fU(\theta^\pos))}$ with the open embedding
$$C^{(\fp,\fU(\theta^\pos))}\times \Bun_P\hookrightarrow
C^\fp\times C^{\fU(\theta^\pos)} \times \stZ_{\emptyset}.$$
It is easy to see that $j_{(\fp,\Theta,\fU(\theta^\pos))}$ is an isomorphism
onto $\stZ_I^{(\fp,\Theta,\fU(\theta^\pos))}$.

\subsubsection{Translational Hecke action}
One of the primary reasons for introducing $\stZ_I$ is that it comes with extra symmetries
given by modifications of the auxiliary $A$-bundle. Namely, we have a diagram
$$
\stZ_I \stackrel{{}h^\la}{\la}{}^\star\CH_{I,A}\stackrel{{}h^\ra}{\ra} \stZ_I
$$
where the Hecke ind-stack ${}^\star\CH_{I,A}$ classifies data
$
(z,\CP'_A,\alpha)
$
where $z\in \stZ_I$ is given by data $(c_I,\CP_G,\CP_A,\sigma)$,
$\CP'_A\in\Bun_A$, $\alpha$ is an isomorphism of $A$-bundles
$$
\alpha:\CP'_A|_{C\setminus |c_I|}
\risom \CP_A|_{C\setminus |c_I|},
$$
$h^\la$ is the obvious projection to $z$, and $h^\ra$ is the projection
to the data $(c_I,\CP_G,\CP'_A,\sigma\circ\alpha)$.

Suppose we fix a partition $\fp$ of the set $I$, and a labelling $\Theta':\fp\to\cowts_A$.
Then we may consider the ind-substack ${}^\star\CH^{(\fp,\Theta')}_{I,A}$ where the coincidences
among $c_I$ are prescribed by $\fp$, and the modifications by $\Theta'$. It is easy to
see that restricting the above diagram gives a diagram
$$
\stZ_I^{(\fp,\Theta,\fU(\theta^\pos))}\stackrel{{}h^\la}{\la}{}^\star\CH^{(\fp,\Theta')}_{I,A}\stackrel{{}h^\ra}{\ra} \stZ_I^{(\fp,\Theta+\Theta',\fU(\theta^\pos))}
$$
in which both projections are isomorphisms.

\subsection{Proof of \propref{horo cat}(1)} 
We will use the ind-stack $\stZ_I$ introduced in the previous section.
Observe that the construction of the fusion product for $Z_I$
extends to $\stZ_I$
in an obvious way. Furthermore, for the canonical map $\fr:Z_I \to \stZ_I$, 
we have a functorial identification
$
\fr^* \circ \product \simeq \product \circ \fr^*.
$
Finally, since $\fr$ respects strata, we conclude that we may use $\stZ_I$ 
to calculate the fusion product for $Z_I$.

The strategy of the proof is as follows. First, we will check that the monoidal
structures on the subcategories generated by the trivial object agree.
Then, we will use the translational symmetry of $\stZ_I$ to extend
this to other objects.

Let $V_i$, for $i\in I$, be a finite collection of vector spaces thought of as trivial representations
of $\ch A_0$. Then we must place a monoidal structure on the functor taking it to the
collection of objects $\IC^0_\stZ\otimes V_i$, for $i\in I$.

\begin{lem}
We have a functorial identification
$$
\product_{i\in I} (\IC^0_\stZ\otimes V_i) \simeq \IC^0_\stZ \otimes (\otimes_{i\in I} V_i)
$$
\end{lem}

\begin{proof}
Recall that the left hand side is obtained by considering the intersection cohomology
sheaf
$\IC^0_{\stZ_{\mr I}} \otimes (\otimes_{i\in I} V_i)$ of the basic stratum closure $\ol \stZ^0_{\mr I}$,
taking its middle-extension to all of $\stZ_I$, then restricting the result
to the locus where the pole points coincide. Observe that this only involves the substack
of $\stZ_I$ consisting of quasimaps without poles. The projection to the base $C^I$
becomes a trivial fibration when restricted to this substack.
Thus the fiber of the middle-extension is canonically isomorphic 
to the right hand side.
\end{proof}

Now to extend this identification to any 
finite collection of representations
of $\ch A_0$, we only need observe that the action of Hecke modifications of the $A$-bundle
at the pole points is transitive on the local strata and clearly compatible
with the fusion product.

\subsection{Proof of \propref{horo cat}(2)} 

Let us recall the main result of \cite{GNhoro04}, which gave 
the following explicit description of the convolution when $X$ is horospherical. 

To state it,  recall that in this case,
the normalizer $P\subset G$ of the stabilizer $S\subset G$ is a parabolic subgroup
with Levi factor $M\subset G$.
We write $2\check\rho_M\in\wts_T$ for the sum of the positive roots of $M$,
and for $\lambda\in\cowts_T$, we write $\langle2\ch\rho_M,\lambda\rangle\in\BZ$
for the natural pairing.
We write $q$ for the natural surjection
$\cowts_T\to\cowts_{A_0}$ 
which may be thought of as the map of weight lattices induced by the inclusion $\ch A_0\to \ch T$.
For $\mu\in\cowts_T$ and a representation $V$ of $\ch G$ we write $V(\mu)$ for the
corresponding weight space.

\begin{thm}\label{thmhoro}
When $X$ is horospherical, for any $V\in \Rep(\ch G)$,
we have
$$
\conv(\CA^V_G)\simeq
\bigoplus_{\kappa\in\cowts_{A_0}}
\bigoplus_{\stackrel
{\mu\in\cowts_{T}}
{q(\mu)=\kappa}}
{}'\IC^{\kappa}_{Z} \otimes V(\mu)
[\langle2\chrho_M,\mu\rangle].
$$
\end{thm}

By the theorem, under the identifications $\Rep(\ch G)\simeq \pshonaffgr$ and 
$\Rep(\ch A_0)\simeq \catq(Z)$,
the restriction $\Rep(\ch G)\to \Rep(\ch A_0)$ and the convolution
$\pshonaffgr\to\catq(Z)$ are canonically isomorphic.
To finish the proof of \propref{horo cat}(2), we need only confirm that the monoidal
structures on these functors agree. This follows immediately from the theorem
and the explicit description of the fusion product in the proof of \propref{horo cat}(1).

\subsection{Towards the proof of \thmref{thmconv}(2)}   \label{proof tedious}

The reader should skip this subsection and return to it when needed in
Section \ref{proof of thmconv(2)}. The main goal is to prove \propref{conv and bad sheaves}
below about the interaction of convolution and bad sheaves on $Z_I$.
It is included here due to the fact that its proof does not use the local model but does use 
$\stZ_I$.

\begin{prop}\label{conv and bad sheaves}
For a nonzero positive coweight
$\theta^\pos\in \cowts^\pos_A$, a decomposition $\fU(\theta^\pos)$, 
and a Hecke equivariant local system $\CL$ on the stratum
$$
Z_{\emptyset}^{\fU(\theta^\pos)}\subset Z_{\emptyset},
$$
we have that
$$
H^I_G(\CA,
\BC_{\mr C^I}\boxtimes
\IC^{ \fU(\theta^\pos)}_{\emptyset}(\CL))
\mbox{ is a bad sheaf.}
$$ 
\end{prop}

\begin{proof}
We will prove the stronger statement that forgetting the Hecke equivariant structure
on 
$
H^{I}_G(\CA,\BC_{\mr C^I}\boxtimes\IC^{ \fU(\theta^\pos)}_{\emptyset}(\CL))
$
gives an ordinary perverse sheaf none of whose summands is
isomorphic to a summand of a perverse sheaf which results from forgetting
the Hecke equivariant structure on an object of $\catq(Z_{I})$.

To use the added flexibility of the ind-stack $\stZ_{I}$,
we first need to confirm that the sheaf $\BC_{\mr C^I}\boxtimes\IC^{ \fU(\theta^\pos)}_{\emptyset}(\CL)$ 
arises as the pullback under $\fr$ of a sheaf on $\stZ_{\mr I}$.

\begin{lem}
For a positive coweight $\theta^\pos\in \cowts^\pos_A$, a decomposition $\fU(\theta^\pos)$, 
and a local system $\CL$ on the stratum
$$
Z_{\emptyset}^{\fU(\theta^\pos)}\subset Z_{\emptyset},
$$
there exists a local system ${}^\star\CL$ on the stratum 
$$
\stZ_{\emptyset}^{\fU(\theta^\pos)}\subset \stZ_{\emptyset},
$$ 
such that $\CL$ is a direct summand of
$
\fr^*({}^\star\CL).
$

\end{lem}

\begin{proof}
As usual, we write $P$ for the normalizer of $S$.
We give a proof in the case when $S\subset G$ is connected so that we may choose a
section $A\to P$ of the projection $P\to A$. In this case, we show that
there is a local system ${}^\star\CL$ and
an isomorphism
$$
\CL\simeq
\fr^*({}^\star\CL).
$$
In the general case, one must work with the components
of 
$
Z_{\emptyset}^{\fU(\theta^\pos)} 
$
separately. Otherwise the argument is the same and we leave the details to the reader.
 
Recall that we have a fibration $\stZ^{\fU(\theta^\pos)}_{\emptyset}\to\Bun_A$
with fiber above the trivial bundle $Z^{\fU(\theta^\pos)}_{\emptyset}$.
Recall as well that have an isomorphism
$$
j_{\fU(\theta^\pos)}:\mr C^{\fU(\theta^\pos)}\times\Bun_P
\risom
\stZ^{\fU(\theta^\pos)}_{\emptyset}.
$$
Thus, choosing a point $c_{\fU(\theta^\pos)}\in \mr C^{\fU(\theta^\pos)}$,
we may define a section
$
\Bun_A\to \stZ^{\fU(\theta^\pos)}_{\emptyset}
$
by the formula
$$
\CP_A\mapsto (c_{\fU(\theta^\pos)},\CP_A(-c_{\fU(\theta^\pos)})\overset{A}{\times}P).
$$
We conclude that the fundamental group of
$\stZ^{\fU(\theta^\pos)}_{\emptyset}$
is a product of that of the base $\Bun_A$ and that of the fiber $Z^{\fU(\theta^\pos)}_{\emptyset}$.
Therefore
we may extend any local system from the fiber to the total space.
\end{proof}

Observe that for any $\CP\in\catp(\stZ_{0,\mr I})$, we have the obvious compatibilty
$$
H^{I}_G(\CA,\fr^*\CP)\simeq\fr^* H^{I}_G(\CA,\CP).
$$
Thus by the previous lemma, 
to prove the proposition, it suffices to show that 
for a local system ${}^\star\CL$ on the horospherical stratum 
$$
\stZ_{\emptyset}^{\fU(\theta^\pos)}\subset \stZ_{\emptyset},
$$ 
we have that
$$
\fr^*(H^{I}_G(\CA,\BC_{\mr C^I}\boxtimes{}^\star\IC_{\emptyset}^{\fU(\theta^\pos)}({}^\star\CL)))
\mbox{ is a bad sheaf.}
$$
Furthermore, we may assume that ${}^\star\CL$ is irreducible,
so that 
under the isomorphism
$$
j_{\fU(\theta^\pos)}:\mr C^{\fU(\theta^\pos)}\times\Bun_P
\risom
\stZ^{\fU(\theta^\pos)}_{\emptyset},
$$
it comes from a product of local systems
$$
j_{\fU(\theta^\pos)*}  (\CL^{\fU(\theta^\pos)}\boxtimes\CL^0)\simeq
{}^\star\CL.
$$

Clearly we have the following.

\begin{lem}
For a positive coweight $\theta^\pos\in \cowts^\pos_A$, a decomposition $\fU(\theta^\pos)$,
and $\CP\in\catp(\stZ_{0,\mr I})$, and $\CC\in\catp(C^{\fU(\theta^\pos)})$,
there is a canonical isomorphism
$$
H^{I}_G(\CA,j_{\fU(\theta^\pos)*}(\CC\boxtimes\CP))
\simeq 
j_{\fU(\theta^\pos)*}(\CC\boxtimes H^{I}_G(\CA,\CP)).
$$
\end{lem}

By the lemma,
it remains to show that
$$
\fr^*
j_{\fU(\theta^\pos)*}
(\IC_{C}^{\fU(\theta^\pos)}(\CL^{\fU(\theta^\pos)})\boxtimes
H_G^{\lambda_I}(\BC_{\mr C^I}\boxtimes\IC^0_{\emptyset}(\CL^0)))
\mbox{ is a bad sheaf.}
$$
But clearly no irreducible summand of such a sheaf could be supported
on the closure of an untwisted local stratum.
This completes the proof of the proposition.
\end{proof}

\newpage


{\Large \part{The local model}\label{loc mod}}

\vspace*{5mm}

As was mentioned in the introduction, the space of quasimaps $Z_I$ is
supposed to model the wildly infinite-dimensional space $X(\CK)$. The 
price we pay is that although $Z_I$ is a more manageable 
object (i.e., it carries a well-defined category of sheaves), it is
not local in nature. 

The goal of this Part is to remedy this ``non-locality". Namely, we will 
introduce spaces $W^\eta_I$ that, on the one hand, will be local
with respect to the curve $C$, and on the other hand, will be equivalent
to $Z_I$ in the smooth topology. 

In Part IV, the machinery developed here will be applied to prove results
stated in Part I. 
The structure of this part is as follows:

\medskip

In \secref{cons loc mod}, we describe a general pattern, pointed out by
Drinfeld, which explains in many cases why a space that is built out of maps from a curve $C$
to a target space $\CY$ will have a local behavior with respect to $C$.
Roughly speaking, this happens when $\CY$ contains an open substack
$\CY^0$ isomorphic to a point $pt$. In our specific case this open substack
will be isomorphic to $pt/F$, where $F$ is a finite abelian group.

Our usual quasimaps may be thought of as maps from $C$
to the stack $X/G$. We show that if we replace $G$ by a certain subgroup
$R$ of the parabolic $P^\op$, we achieve the desired locality.

\medskip

In \secref{sect factor}, we establish a version of the factorization property
of the local model $W^\eta_I$ which is the expression of its locality with respect to
$C$.

\medskip

In \secref{sect fibers}, we show that the fibers of the natural projection from
$W^\eta_I$ to a suitable configuration space can be completely described in
terms of the affine Grassmannian of $G$.

\medskip

In \secref{sect rel quasimaps}, we show that the quasimaps space $Z_I$
and the local model $W^\eta_I$ are essentially equivalent in the smooth
topology. This allows us to reduce questions about the local behavior of
the former to those about the latter.

\vspace*{10mm}

\section{Construction of the local model}   \label{cons loc mod}

\subsection{Base ind-scheme}

The input for the construction described below is the torus $A$,
and the seimgroup of dominant weights $\domwts_X$, or alternatively
the semigroup of positive coweights $\poscowts_X\subset\cowts_A$.

Note that a rational section $\tau$ of an $A$-bundle on $C$
defines a $\cowts_A$-valued
divisor $\on{div}(\tau)$ on $C$. 
The support of $\on{div}(\tau)$ is a finite subset of $C$ which we denote by $|\tau|$.

For a finite set $I$, and $\eta\in\cowts_A$,
define the ind-scheme $C_I^{\eta}$ to be that classifying data
$$
(c_I;\CP_A,\tau),
$$
where $c_I\in C^I$, $\CP_A\in\Bun^\eta_A$, and $\tau$ is a rational section of $\CP_A$
such that the following holds.
For any $\lambda\in\domwts_X$, let $\CL^\lambda$ be the corresponding 
one-dimensional representation
of $A$, and let $\CL^\lambda$ be the line bundle induced from the
$A$-bundle $\CP_A$. Then we require that the meromorphic section of $\CL^\lambda$ associated to $\tau$
be regular on $C\setminus |c_I|$.
We write $c_I^\eta$ for a point of $C^\eta_I$.
It is the same thing as a point $c_I\in C^I$, and
a $\cowts_A$-valued divisor on $C$
of degree $\eta$ which takes values in $\poscowts_X$ on $C\setminus|c_I|.$
Note that $C^\eta_I$ is indeed an ind-scheme since $\domwts_X$ generates $\wts_A$.

For a point $c_I^\eta\in C_I^\eta$, 
we call the subset $|c_I|\subset C$ the {pole points} of $c_I^\eta$.
We call the union $|c_I|\cup|\tau|\subset C$
the {degeneracy locus} of $c_I^\eta$, and denote it by $|c_I^\eta|$. 

\subsection{Base ind-stack}

The additional input for the construction described here is 
the torus $\maxtor$, 
and the short exact sequence
$$
1\to F\to \maxtor\to A\to 1,
$$ 
or equivalently the short exact sequence
$$
0\to\cowts_{A_0}\to\cowts_A\to F\to 0.
$$

For a finite set $I$, and $\eta\in\cowts_{A_0}$,
define the ind-stack $ M^\eta_I$ to be 
the fiber product
$$
 M^\eta_I = \Bun_{A_0}\underset{\Bun_A}{\times} C^\eta_I
$$
where the map $\Bun_{A_0}\to\Bun_A$ is the induction 
$$
\CP_{A_0}\mapsto \CP_{A_0}\overset{A_0}{\times} A,
$$
and the map $ C^\eta_I\to\Bun_A$ is the obvious projection.
In other words, the stack $ M^\eta_I$ classifies data
$$
(c_I;\CP_A, \tau; \CP_{\maxtor}, \tau_0)
$$
where $(c_I;\CP_A,\tau)\in C_I^{\eta}$,
$\CP_{\maxtor}\in\Bun^{\eta}_{\maxtor}$, 
and $\tau_0$ is an $A$-equivariant isomorphism
$$
\tau_0: \CP_{\maxtor} \overset{\maxtor}{\times} A\risom \CP_A.
$$
Note that we assume $\eta\in\cowts_{A_0}$,
or else the stack $ M^\eta_I$ would be empty.

We have the natural projection 
$$
 M^\eta_I\to C_I^\eta.
$$
For a point $m_I^\eta\in M^\eta_I$, over a point $c_I^\eta\in C^\eta_I$,
we call the subset $|c_I|\subset C$ the {pole points} of $m_I^\eta$, and denote it by $|m_I|$.
We call the subset
$|c^\eta_I|\subset C$ the {degeneracy locus} of $m_I^\eta$, and denote it by $|m_I^\eta|$.

To a point $m_I^\eta\in M^\eta_I$, 
with representative 
$(c_I;\CP_A, \tau; \CP_{\maxtor}, \tau_0)$, 
we may associate an $F$-bundle 
on $C\setminus |m_I^\eta|$ via
the pullback of $\tau$ under $\tau_0|_{C\setminus |m_I^\eta|}$.
Since $F$ is finite, this is the same thing as an $F$-local system
on $C\setminus|m_I^\eta|$. 


\subsection{The general pattern}

\subsubsection{A simplified version}   \label{simple version}

Let us recall the following general construction\cite{BFG03}, Sect. 2.16. Let $\CY$ be 
an algebraic stack with an open substack $\CY^0\subset \CY$ isomorphic to 
a point $pt$.

Assume that we are given an $A$-bundle $\CP_{A,\CY}$ over $\CY$,
and its trivialization over $\CY^0$. Moreover, assume that for every 
$\lambda\in \domwts_X$ the meromorphic section of the corresponding line
bundle $\CL_{\CY}^\lambda$ is regular, and that $\CY^0$ is the locus of
non-vanishing of these sections.

For a curve $C$ we can consider the space (i.e., functor on the category of schemes)
$\on{Maps}_I(C,\CY)$ over $C^I$
that classifies maps $$(C\setminus |c_I|)\to \CY$$ 
such that all but finitely many points of $C\setminus |c_I|$ get 
mapped to $\CY^0$. The space $\on{Maps}_I(C,\CY)$ splits into connected components
according to the degree $\eta\in \poscowts_X$ of the pull-back of $\CP_{A,\CY}$.

By construction, we have a canonical map
$$\on{Maps}_I(C,\CY)^\eta\to C^\eta_I.$$
For a given map $\sigma:C\setminus |c_I|\to \CY$ the locus of $C\setminus |c_I|$ for which this map
lands in $\CY^0$ equals the complement of the degeneracy locus of the resulting
element of $C^\eta_I$.

\subsubsection{A generalization}   \label{Drinf gen}

Next let $\CY$ be an algebraic stack equipped with an $A_0$-torsor $\CP_{A_0,\CY}$.
Let $\CP_{A,\CY}$ denote the induced $A$-bundle, and assume that it is
trivialized over an open substack $\CY^0\subset \CY$, such that this trivialization
has the same properties as in Section \ref{simple version}.

In particular, over $\CY^0$, the $A_0$-torsor $\CP_{A_0,\CY}$ admits a
canonical reduction to an $F$-torsor $\CP_{F,\CY}$. We assume furthermore
that the resulting map
$$\CY^0\to pt/F$$ is an isomorphism.

\medskip

We define the space $\on{Maps}_I(C,\CY)$ in the same way as above.
By construction, we have a canonical map
$$\on{Maps}^\eta_I(C,\CY)\to M_I^\eta.$$

Note that $\on{Maps}^\eta_I(C,\CY)$ is empty unless $\eta\in \cowts_{A_0}$.

\subsubsection{The case of spherical varieties}

Let us now explain in what situation we will apply the pattern of Section
\ref{Drinf gen}.

Recall the canonical surjection 
$$P^\op\twoheadrightarrow P^\op/U^\op\simeq M\twoheadrightarrow A_0,$$ and 
let us choose a splitting $M\hookleftarrow A_0$. Let $R$ denote the preimage
of $A_0$ in $P^\op$. By construction, $R\cap S\simeq F$.

We set $\CY:=X/R$ and $\CY_0:=\openX^+/R$. The $A_0$-torsor $\CP_{A_0,\CY}$
is the pull-back of the tautological $A_0$-torsor under
$$X/R\to pt/R\to pt/A_0.$$

The lines $\fl^\lambda\subset \BC[X]$ for $\lambda\in \domwts_X$ define
the trivialization of the induced $A$-torsor over $\CY_0$ with the required properties.

\medskip

We have the resulting space 
$$W_I^\eta:= \on{Maps}^\eta_I(C,\openX^+/R).
$$
Let us describe it in terms of bundles on $C$, which will in particular imply that
$W_I^\eta$ is an ind-algebraic stack.

\medskip

Namely, $W_I^\eta$ classifies the data of
$$
(c_I;\CP_A, \tau; \CP_{\maxtor},
\tau_0;
\CP_{R},\tau_R,\sigma)
$$
where
$(c_I;\CP_A, \tau; \CP_{\maxtor}, \tau_0)\in M^\eta_I$,
$\CP_{R}\in\Bun_{R}$, 
$\tau_R$ is an $\maxtor$-equivariant isomorphism
$$
\tau_R:\CP_{R}\overset{R}{\times}\maxtor\risom \CP_{\maxtor},
$$
and $\sigma$ is a section 
$$
\sigma:C\setminus |c_I|\to\CP_{R}\overset{R}{\times} X|_{C\setminus|c_I|}.
$$
For some open curve $C'\subset C\setminus |c_I|$,
the section $\sigma$ 
is required to factor
$$
\sigma|_{C'}:C'\to 
\CP_{R}{\overset{R}{\times}}{\openX^+}|_{C'}\to
\CP_{R}{\overset{R}{\times}}{X}|_{C'},
$$
and the composition
$$
C'\stackrel{\sigma|_{C'}}\to 
\CP_{R}{\overset{R}{\times}}{\openX^+}|_{C'}
\stackrel{}{\to}
\CP_A|_{C'} 
$$
is required to coincide with the rational section
$$
\tau|_{C'}:C'\stackrel{}{\to}\CP_A|_{C'}.
$$
Note that we assume $\eta\in\cowts_{A_0}$,
or else the stack $ W^\eta_I$ would be empty.


\section{Factorization}   \label{sect factor}

\subsection{``Simple" factorization}

Assume for a moment that we are in the context of Section \ref{simple version}. 
For $\eta_1,\eta_2\in \poscowts_X$, consider the natural map
$$
C_{I_1}^{\eta_1}\times C_{I_2}^{\eta_2}\to
C_{I_1\cup I_2}^{\eta_1+\eta_2}$$ 
which we denote by
$$(c_{I_1}^{\eta_1},c_{I_2}^{\eta_2})\mapsto c_{I_1}^{\eta_1}\uplus c_{I_2}^{\eta_2}.$$
Its restriction to 
the open subscheme 
$$(C_{I_1}^{\eta_1}\times C_{I_2}^{\eta_2})_\disj \subset
C_{I_1}^{\eta_1}\times C_{I_2}^{\eta_2}$$
of points $(c_{I_1}^{\eta_1},c_{I_2}^{\eta_2})$ 
such that $|c_{I_1}^{\eta_1}|$ is disjoint from $|c_{I_2}^{\eta_2}|$
is \'etale.

\medskip

We claim
\begin{multline} \label{simple factorization}
\bigl(\on{Maps}_I(C,\CY)^{\eta_1}\times \on{Maps}_I(C,\CY)^{\eta_2}\bigr)
\underset{C_{I_1}^{\eta_1}\times C_{I_2}^{\eta_2}}\times
(C_{I_1}^{\eta_1}\times C_{I_2}^{\eta_2})_\disj \\
\simeq \on{Maps}_I(C,\CY)^{\eta_1+\eta_2}\underset{C_{I_1\cup I_2}^{\eta_1+\eta_2}}
\times (C_{I_1}^{\eta_1}\times C_{I_2}^{\eta_2})_\disj.
\end{multline}

We refer the reader to \cite{BFG03} for the proof. We shall now generalize
this to the context where the finite group is present.

\subsection{Factorization and the base ind-stack}

We have a natural map
$$
 M^{\eta_1}_{I_1}{\times} M^{\eta_2}_{I_2}\to
 M^{\eta_1+\eta_2}_{I_1\cup I_2}
$$
defined by
$$
((c_{I_1}^{\eta_1}; \CP^1_{\maxtor}, \tau^1_0),
(c_{I_2}^{\eta_2}; \CP^2_{\maxtor}, \tau^2_0))\mapsto
(c_{I_1}^{\eta_1}\uplus c_{I_2}^{\eta_2}; 
\CP^1_{\maxtor}\otimes\CP^2_{\maxtor}, \tau^1_0\otimes\tau^2_0).
$$

Consider the fiber products
$$
(M^{\eta_1}_{I_1}{\times} M^{\eta_2}_{I_2})_\disj
=
 (M^{\eta_1}_{I_1}{\times} M^{\eta_2}_{I_2})
\underset{C_{I_1}^{\eta_1}\times C_{I_2}^{\eta_2}}{\times}
(C_{I_1}^{\eta_1}\times C_{I_2}^{\eta_2})_\disj
$$
$$
( M^{\eta_1+\eta_2}_{I_1\cup I_2})_\disj
=
 M^{\eta_1+\eta_2}_{I_1\cup I_2}
\underset{C_{I_1\cup I_2}^{\eta_1+\eta_2}}{\times}
(C_{I_1}^{\eta_1}\times C_{I_2}^{\eta_2})_\disj.
$$
It is easy to see that the induced map
$$( M^{\eta_1}_{I_1}{\times} M^{\eta_2}_{I_2})_\disj
\to
( M^{\eta_1+\eta_2}_{I_1\cup I_2})_\disj
$$
is an \'etale cover.

\medskip

To glue mapping spaces, we will also need the stack 
$$
(M^{\eta_1}_{I_1}{\times} M^{\eta_2}_{I_2})^\sim_\disj
$$
over 
$
(M^{\eta_1}_{I_1}{\times} M^{\eta_2}_{I_2})_\disj
$
whose fiber over a point
$((c_{I_1}^{\eta_1}; \CP^1_{\maxtor}, \tau^1_0),(c_{I_2}^{\eta_2}; \CP^2_{\maxtor}, \tau^2_0))$
is the data of trivialization of the $F$-bundle $\CP^1_F$ over the finite scheme $|c^{\eta_2}_{I_2}|$
and of the $F$-bundle $\CP^2_F$ over the finite scheme $|c^{\eta_1}_{I_1}|$.
It is clear that the forgetful map
$$(M^{\eta_1}_{I_1}{\times} M^{\eta_2}_{I_2})^\sim_\disj\to
(M^{\eta_1}_{I_1}{\times} M^{\eta_2}_{I_2})_\disj$$
is also an \'etale cover.

\subsection{The general case} \label{sect factorgen}

Now let $\CY^0\subset \CY$ be as in \secref{Drinf gen}.

\begin{prop}      \label{propfact}
There is a canonical isomorphism
\begin{multline*}
\bigl(\on{Maps}_{I_1}(C,\CY)^{\eta_1}\times \on{Maps}_{I_2}(C,\CY)^{\eta_2}\bigr)
\underset{M^{\eta_1}_{I_1}{\times} M^{\eta_2}_{I_2}}\times 
(M^{\eta_1}_{I_1}{\times} M^{\eta_2}_{I_2})^\sim_\disj  \\
\simeq \on{Maps}_{I_1\cup I_2}(C,\CY)^{\eta_1+\eta_2}\underset{M^{\eta_1+\eta_2}_{I_1\cup I_2}}
\times (M^{\eta_1}_{I_1}{\times} M^{\eta_2}_{I_2})^\sim_\disj.
\end{multline*}
\end{prop}

\begin{proof} 
Let $f_1,f_2$ be the maps in the data of an element of the left hand side.
We define the corresponding map $f$ of the right hand side as follows.
First, break up the curve $C$ as the union of the punctured curve
$C \setminus (|c_{I_1}^{\eta_1}|\cup |c_{I_2}^{\eta_2}|)$
and the disjoint completed formal neighborhoods $D_{|c_{I_1}^{\eta_1}|}$, $D_{|c_{I_2}^{\eta_2}|}$.
Define the map $f$ on the pieces to be
$$
f_1\otimes f_2: {C \setminus (|c_{I_1}^{\eta_1}|\cup |c_{I_2}^{\eta_2}|)} \to \CY^0 = pt/F
$$ 
$$
f_1:D_{|c_{I_1}^{\eta_1}|}\to \CY
\qquad
f_2:D_{|c_{I_2}^{\eta_2}|}\to \CY
$$
Here we have written $f_1\otimes f_2$ to denote the map classifying the
$F$-bundle obtained by tensoring
the $F$-bundles $\CP^1_F, \CP^2_F$ classified by $f_1, f_2$ respectively.

Now to see that the pieces canonically glue, we use the remaining data. Namely, we have trivializations
$\tau_1,\tau_2$ of the $F$-bundles $\CP^1_F, \CP^2_F$ over the finite schemes 
$|c_{I_2}^{\eta_2}|, |c_{I_1}^{\eta_1}|$ respectively. Observe that $\tau_1,\tau_2$
canonically extend to trivializations of $\CP^1_F, \CP^2_F$ over 
$D_{|c_{I_2}^{\eta_2}|}$, $D_{|c_{I_1}^{\eta_1}|}$ respectively.
Thus $\tau_1,\tau_2$ provide isomorphisms of $F$-bundles
over the punctured completed formal neighborhoods
$$
\CP^1_F|_{D^\times _{|c_{I_1}^{\eta_1}|}} \stackrel{\tau_2}{\simeq}
(\CP^1_F\otimes \CP^2_F)
|_{D^\times _{|c_{I_1}^{\eta_1}|} }
\qquad
\CP^2_F|_{D^\times _{|c_{I_2}^{\eta_2}|}} \stackrel{\tau_1}{\simeq }
(\CP^1_F\otimes \CP^2_F)
|_{D^\times _{|c_{I_2}^{\eta_2}|} }
$$
Thus we may identify the restricted maps 
$$
f_1|_{D^\times _{|c_{I_1}^{\eta_1}|} } \simeq (f_1\otimes f_2)|_{D^\times _{|c_{I_1}^{\eta_1}|} }
\qquad
f_2|_{D^\times _{|c_{I_2}^{\eta_2}|} } \simeq (f_1\otimes f_2)|_{D^\times _{|c_{I_2}^{\eta_2}|} }
$$
and glue them to obtain a single map $f$. Clearly this identifies the moduli problem
of the left hand side with that of the right.
\end{proof}


\subsection{Complements: open curves} \label{sect open curves}

It is useful to have generalizations of the constructions of the preceding
sections for open curves. Our primary application will be the following. In the context of
\secref{Drinf gen},
the moduli problem
$\on{Maps}_I(C,\CY)^{\eta}$ involves an $F$-bundle $\CP_F$ on the curve
$C\setminus |c_I^\eta|$. Suppose we would like to apply the factorization of \secref{sect factorgen}.
 It may turn out that
\propref{propfact} has less content than needed. Namely, given a particular
$F$-bundle $\CP_F$ on $C\setminus (|c_{I_1}^{\eta_1}|\cup |c_{I_2}^{\eta_2}|)$,
there may not be any $F$-bundles on $C\setminus|c_{I_1}^{\eta_1}|$ and $C\setminus|c_{I_2}^{\eta_2}|$
whose tensor is equal to $\CP_F$. Thus the base changes on both sides of \propref{propfact} will miss
any map living over $\CP_F$.
But if we remove an auxiliary point $c\in C$, then we may always use the added flexibility
of nontrivial $F$-monodromies around $c$ to find 
$F$-bundles on $C\setminus(c\cup |c_{I_1}^{\eta_1}|)$ and $C\setminus (c\cup |c_{I_2}^{\eta_2}|)$
whose tensor is equal to $\CP_F$. With this motivation,
we outline below how one may generalize
our constructions to the open curve
$C\setminus c$.

\medskip

Recall that for a finite set $I$, and $\eta\in\cowts_A$, a point $c_I^\eta$ of the ind-scheme $C_I^\eta$
classifies the data of a point $c_I\in C^I$, and a
$\cowts_A$-valued divisor on $C$ of degree $\eta$ which
takes values in $\poscowts_A$ on $C\setminus|c_I|$.
Thus it makes sense to consider $C_I^\eta$ for a not necessarily complete curve
such as $C\setminus c$.

\medskip

For a finite set $I$, and $\eta\in\cowts_A$,
define the ind-stack 
$ M^{\eta}_{I,1}$ to be that classifying data 
$$
(c;c_I^\eta; \CP_{\maxtor}, \tau_0)
$$
where $c\in C$, $c_I^\eta\in (C\setminus c)_I^{\eta}$,
$\CP_{\maxtor}$ is an ${\maxtor}$-bundle on $C\setminus c$,
and $\tau_0$ is an $A$-equivariant isomorphism
$$
\tau_0 : \CP_{A_0}\overset{A_0}{\times} A
\risom \CP^0_A(c_I^\eta),
$$
where $\CP_A^0(c_I^\eta)$ denotes the trivial $A$-bundle on $C\setminus c$,
twisted by the $\cowts_A$-valued divisor associated to $c_I^\eta$.
Note that even if $\eta\in\cowts_A$ is not in $\cowts_{A_0}\subset\cowts_A$, the stack 
$ M^{\eta}_{I,1}$ still makes sense and is nonempty.

For the moment,
consider the ind-stack $\ol M^\eta_{I,1}$ that classifies the same data
as $ M^\eta_{I,1}$ except that the ${\maxtor}$-bundle $\CP_{\maxtor}$
it classifies
is defined on $C$ rather than $C\setminus c$.
The following lemma confirms that this makes no difference.
Its assertion is completely local and we leave its proof to the reader.

\begin{lem}\label{lemextend}
The natural restriction $\ol M^\eta_{I,1}\to M^\eta_{I, 1}$ is an open and closed embedding.
\end{lem}

For a finite set $I$, and $\eta\in\cowts_A$,
define the open ind-substack
$$
( M^\eta_I\times C)_\disj\subset  M^\eta_I\times C
$$
to consist of pairs $(m^\eta_I,c)$ such that $c$ is disjoint from the degeneracy locus $|m_I^\eta|$.
The following is immediately implied by Lemma~\ref{lemextend}.

\begin{lem}
The natural map 
$$
( M^\eta_I\times C)_\disj\to  M^{\eta}_{I,1}
$$
is an open and closed embedding.
\end{lem}

\medskip

Now for a finite set $I$, and $\eta\in\cowts_A$,
define the 
space $\on{Maps}_{I,1}(C,\CY)^{\eta}$ to be that classifying $c\in C$, $c_I\in (C\setminus c)^I$,
and a map $C \setminus (|c_I| \cup c) \to \CY$ as described in \secref{Drinf gen}.
By construction, we have a canonical map
$$
\on{Maps}^\eta_{I,1}(C,\CY)\to M_{I,1}^\eta.
$$
We also have an obvious analogue of the factorization of \propref{propfact}.
Note that even if $\eta\in\cowts_A$ is not in $\cowts_{A_0}\subset\cowts_A$, the space 
$\on{Maps}_{I,1}(C,\CY)^{\eta}$ still makes sense and is nonempty.

\medskip

In the case where $\CY=X/R$, we set $ W^{\eta}_{I,1}:=\on{Maps}_{I,1}(C,\CY)^{\eta}$.
It is an ind-algebraic stack classifying data
$$
(c;c_I^\eta; \CP_{\maxtor},
\tau_0;
\CP_{R},\tau_R,\sigma)
$$
where $(c;c_I^\eta; \CP_{\maxtor}, \tau_0)\in M^{\eta}_{I,1}$,
$\CP_{R}$ is an ${R}$-bundle on $C\setminus c$,
$\tau_R$ is an $\maxtor$-equivariant isomorphism
$$
\tau_R:\CP_{R}\overset{R}{\times}\maxtor\risom \CP_{\maxtor},
$$
and $\sigma$ is a section 
$$
\sigma:C\setminus(c\cup |c_I|)\to\CP_{R}\overset{R}{\times} X|_{C\setminus(c\cup |c_I|)}.
$$
For some open curve $C'\subset C\setminus(c\cup|c_I|)$,
the section $\sigma$ 
is required to factor
$$
\sigma|_{C'}:C'\to 
\CP_{R}{\overset{R}{\times}}{\openX^+}|_{C'}\to
\CP_{R}{\overset{R}{\times}}{X}|_{C'},
$$
and the composition 
$$
C'\stackrel{\sigma|_{C'}}\to 
\CP_{R}{\overset{R}{\times}}{\openX^+}|_{C'}
\stackrel{}{\to}
 \CP_A^0(c_I^\eta)|_{C'} 
$$
is required to have divisor $c_I^\eta$.
As explained in \secref{sect fibers} below,
the canonical map
$$
W^\eta_{I,1}\to M_{I,1}^\eta.
$$
is ind-representable.

Finally, for a finite set $I$, and $\eta\in\cowts_A$,
define the open ind-substack
$$
( W_{I}^{\eta}\times C)_\disj\subset W_{I}^{\eta}\times C
$$
to consist of pairs such that $c$ is disjoint from the degeneracy locus $|c_I^\eta|$.
The following is immediately implied by Lemma~\ref{lemextend}.

\begin{lem}
The natural map 
$$
(W_{I}^{\eta}\times C)_\disj\to W_{I,1}^{\eta}
$$
is an open and closed embedding.
\end{lem}


\section{Description of fibers}    \label{sect fibers}

\subsection{Relation to the affine Grassmannian}

From now on we will specialize to the case where $\CY=X/R$, and so 
$\on{Maps}_I(C,\CY)^\eta=W^\eta_I$. Our present goal is to describe
the fibers of the morphism $W^\eta_I\to M^\eta_I$ in terms of
the affine Grassmannian $\affgr_G$. This description will imply,
in particular, that the above morphism is (ind)-representable even
if $C$ is not complete. First, we will do this on a point-wise level.

\medskip

Given a point $m_I^\eta\in M_I^\eta$, we have a divisor
$c_I^\eta=\Sigma\, \eta_k \cdot c_k\in C_I^\eta$, and an $F$-bundle $\CP_F$
over $C\setminus|c^\eta_I|$. Recall that by construction, $\eta_k$ is arbitrary when $c_k$
is a pole point, but $\eta_k$ is constrained to lie in $\poscowts_X$ otherwise.
By restriction, we obtain an $F$-torsor $\CP_F^k$ 
on the formal neighborhood $D_{c_k}$ of each $c_k$.

Recall that $F\simeq R\cap S$ is a subgroup of $G$. Consider the 
twisted version of the affine Grassmannian $\affgr_{G,\CP^k_F}$
that 
classifies the data of a $G$-torsor $\CP_G$ on $D_{c_k}$ and an identification
$$
\beta:\CP_G\simeq G\overset{F}\times \CP^k_{F}|_{D_{c_k}^\times}.
$$
By \cite{BL}, this is equivalent to giving $\CP_G$ over
an open subset $U\subset C$ with
$$
c_k\in U\subset C\setminus (\underset{k'\neq k}\cup c_{k'})
$$ 
and giving $\beta$ over $U\setminus c_k$.

Let $(W^\eta_I)_{m_I^\eta}$ denote the fiber of $W^\eta_I$ above $m^\eta_I$.
We claim that there is natural morphism
\begin{equation} \label{map to affgr}
(W^\eta_I)_{m_I^\eta}\to \underset{k}\Pi\, 
\affgr_{G,\CP_F^k}.
\end{equation}
Namely, consider the restriction of the data 
to each $D_{c_k}$. We obtain
a $G$-bundle $\CP_G$, endowed with a reduction to $R$, and a reduction
to $S$ over $D^\times_{c_k}$. Furthermore,
the resulting map $D^\times_{c_k}\to S\backslash G/R$ hits the open
substack 
$$
pt/F\simeq \openX^+/R\subset S\backslash G/R
$$
and the induced $F$-torsor is equal to $\CP^k_F$. Thus the restriction
of $\CP_G$ to each $D^\times_{c_k}$ is induced from $\CP^k_F$.

\subsubsection{Description of the image}    \label{descr image}

We will now show that the map \eqref{map to affgr} is an (ind)-locally
closed embedding, and describe its image.

\medskip

First, let $\affgr_{R,\CP_F^k}$ be the corresponding twisted version of
the affine Grassmannian of the group $R$. The projection $R\to A$ induces a canonical map
$$
\affgr_{R,\CP_F^k}\to \affgr_A,
$$
and we write 
$
\affgr^{\eta_k}_{R,\CP_F^k}\subset \affgr_{R,\CP^k_F}
$ for the 
preimage of the corresponding connected component $\affgr_A^{\eta_k}\subset \affgr_A$.
The inclusion $R\subset G$ induces a canonical locally closed embedding
$$
\affgr^{\eta_k}_{R,\CP_F^k}\hookrightarrow \affgr_{G,\CP_F^k}.
$$
In fact, one can check that $\affgr^{\eta_k}_{R,\CP_F^k}$ is the orbit of a connected
component of a twisted version of the group 
$U^\op(\CK_{c_k})\times A_0(\CO_{c_k})$ acting on $\affgr_{G,\CP_F^k}$.

By construction, the map \eqref{map to affgr} factors through 
$\affgr^{\eta_k}_{R,\CP_F^k}$.

\medskip

Next, let $\CG_{F^k}$ be the $G(\CO_{c_k})$-torsor over $\affgr_{G,\CP_F^k}$
that classifies triples $(\CP_G,\alpha,\beta)$, where $(\CP_G,\beta)$
is as in the definition of $\affgr_{G,\CP_F^k}$, and $\alpha$ is a trivialization
of $\CP_G$ over $D_{c_k}$.
We have a natural map
\begin{equation} \label{loops to loops}
\CG_{F^k}\to X(\CK_{c_k})\setminus(X\setminus\openX)(\CK_{c_k})
\end{equation}

For an element $\theta\in \CV(G/S)$ recall that we write 
$$\bO_\theta\subset X(\CK_{c_k})\setminus(X\setminus\openX)(\CK_{c_k})$$
for the corresponding $G(\CO_{c_k})$-orbit. Let $\overline{\bO}_{\theta}$
denote its closure.

The preimage of $\overline{\bO}_{\theta}$ under the map \eqref{loops to loops}
is a $G(\CO_{c_k})$-invariant closed subscheme of $\CG_{F^k}$. By invariance,
we have the
corresponding closed subscheme
$$\ol\affgr_{G,\CP_F^k}^{S,\theta}\subset \affgr_{G,\CP_F^k}.
$$ 
Similarly, let
$\affgr_{G,\CP_F^k}^{S,\theta}$ be the open subset of 
$\ol\affgr_{G,\CP_F^k}^{S,\theta}$ corresponding to $\bO_\theta$.

\medskip

Now, fix a labeling 
$
\Theta: k \mapsto \theta_k\in \CV(G/S)
$ 
such that $\theta_k$ lies in $\poscowts_X$ when $c_k$ is not a pole point.
Let $(W_I^{\eta,\Theta})_{m^\eta_I}$
(resp., $(\ol{W}_I^{\eta,\Theta})_{m^\eta_I}$) be the locally closed (resp., closed)
substack of $(W_I^\eta)_{m^\eta_I}$, corresponding to the condition that
for each $k$, the map $D^\times_{c_k}\to X$ (defined up to $G(\CO_{c_k})$-conjugacy)
lies in $\bO_{\theta_k}$ (resp., $\ol\bO_{\theta_k}$). 

Each piece $(\ol{W}_I^{\eta,\Theta})_{m^\eta_I}$ is an algebraic stack,
the entire $(W_I^{\eta})_{m^\eta_I}$ is their union, and the various pieces
$(W_I^{\eta,\Theta})_{m^\eta_I}$ define a stratification.

By construction, we immediately have the following.

\begin{lem}   \label{descr fibers}
The map \eqref{map to affgr} defines isomorphisms
$$(W_I^{\eta,\Theta})_{m^\eta_I}\simeq \underset{k}\Pi\,
\affgr^{\eta_k}_{R,\CP_F^k}\cap \affgr_{G,\CP_F^k}^{S,\theta_k}$$
$$(\ol{W}_I^{\eta,\Theta})_{m^\eta_I}\simeq \underset{k}\Pi\,
\affgr^{\eta_k}_{R,\CP_F^k}\cap \ol\affgr_{G,\CP_F^k}^{S,\theta_k}.$$
\end{lem}

\subsubsection{Description in families}

The contents of the previous two subsections can be repeated over the base
$M^\eta_I$ rather than for individual fibers. Let us spell out the relevant definitions,
amended slightly to meet our future needs. 

\medskip

We define the twisted version of
the Beilinson-Drinfeld Grassmannian $\affgr_{G,M^\eta_I}$ over $M^\eta_I$
to classify data
$(m^\eta_I;\CP_G,\beta)
$ 
where $m_I^\eta\in M^\eta_I$ with associated generic $F$-bundle $\CP_F$,
$\CP_G$ is a $G$-torsor over $C$, and $\beta$ is a reduction of $\CP_G$ to $\CP_F$
away from $|m^\eta_I|$. The fiber of $\affgr_{G,M^\eta_I}$ over a given point
$m^\eta_I$ is the product 
$$
(\affgr_{G,M^\eta_I})_{m^\eta_I}\simeq
\underset{k}\Pi\, \affgr_{G,\CP_F^k}.
$$

We similarly define $\affgr_{R,M^\eta_I}$ and its connected
component $\affgr^\eta_{R,M^\eta_I}$ which are both locally closed ind-subschemes
of $\affgr_{G,M^\eta_I}$.

\medskip

Now, for simplicity, we focus on the open locus
of $W^\eta_{I}$ where the pole points are distinct. 
Fix a labeling 
$$
\Theta:I \to\CV(G/S)
\qquad
i \mapsto \theta_i.
$$ 
We have the locally closed substack
$\affgr^{S,\Theta}_{G,M^\eta_I}$ (resp. closed substack
$\ol\affgr^{S,\Theta}_{G,M^\eta_I}$)
of $\affgr_{G,M^\eta_I}$  whose fiber over $m^\eta_I\in M^\eta_I$ with distinct pole points
is the product
$$
\underset{i\in I}\Pi\, \affgr_{G,\CP^i_F}^{S,\theta_i}
\times
\underset{k\not \in I}\Pi\, \affgr_{G,\CP^k_F}^{S,0}
\mbox{ (resp. }
\underset{i}\Pi\, \ol\affgr_{G,\CP^i_F}^{S,\theta_i}
\times 
\underset{k\not \in I}\Pi\, \ol\affgr_{G,\CP^k_F}^{S,0}
\mbox{)}.
$$

We also have the locally closed
substack $W_I^{\eta,\Theta}$ (resp., closed substack $\ol{W}_I^{\eta,\Theta}$)
corresponding to the following conditions.
For each $i\in I$, the map
$D^\times_{c_i}\to X$ (defined up to $G(\CO_{c_i})$-conjugacy)
lies in $\bO_{\theta_i}$ (resp., $\ol\bO_{\theta_i}$). 
For each $k\not \in I$, the map
$D^\times_{c_k}\to X$ (defined up to $G(\CO_{c_k})$-conjugacy)
lies in $\bO_{0}$ (resp., $\ol\bO_{0}$).

By construction, we immediately have the following.

\begin{lem}
We have canonical isomorphisms
$$W_I^{\eta,\Theta}\simeq \affgr^\eta_{R,M^\eta_I}\underset{M^\eta_I}\times
\affgr^{S,\Theta}_{G,M^\eta_I}$$
$$
\ol{W}_I^{\eta,\Theta}\simeq \affgr^\eta_{R,M^\eta_I}\underset{M^\eta_I}\times
\ol\affgr^{S,\Theta}_{G,M^\eta_I}.$$
\end{lem}


\subsection{The transverse locus}   \label{sect trans locus}

Let $\eta$ be an element of $\CV(G/S)$, and let $\CP_F$ be an
$F$-bundle on the formal punctured disc $D^\times$ with monodromy 
around the origin equal to the image of $\eta$ in $\cowts_A/\cowts_{A_0}$.

The following is a straightforward reformulation of \propref{local trans}.

\begin{prop}  \hfill  \label{trans prop}

\smallskip

\noindent(1) The intersection
$$\affgr^{\eta}_{R,\CP_F}\cap \affgr_{G,\CP_F}^{S,\eta}$$
is a point-scheme.

\smallskip

\noindent(2) 
The inclusion
$$\affgr^{\eta}_{R,\CP_F}\cap \affgr_{G,\CP_F}^{S,\eta}\hookrightarrow
\affgr^{\eta}_{R,\CP_F}\cap\ol\affgr_{G,\CP_F}^{S,\eta}
$$
is an isomorphism.
\end{prop}

Let $C_I^{\eta,+}$ be the subscheme of $C_I^\eta$
where we require that for $c^\eta_I=\Sigma \, \eta_k\cdot c_k$ such that all $\eta_k$ belong
to $\CV(G/S)$. Let $M^{\eta,+}_I$ be the corresponding substack 
of $M^\eta_I$ obtained by base change.

\medskip

Now consider the closed substack of the base change 
$$
W^{\eta,+}_I\subset W^\eta_I\underset{M^\eta_I}\times M^{\eta,+}_I
$$
corresponding at the level of fibers to
$$\underset{k}\Pi\,
\affgr^{\eta_k}_{R,\CP^k_F}\cap \ol\affgr_{G,\CP^k_F}^{S,\eta_k}
$$
in terms of the identification of \lemref{descr fibers}.

\propref{trans prop} immediately implies the following.

\begin{cor}  \label{deep stratum}
The projection $W^{\eta,+}_I\to M^{\eta,+}_I$ is an isomorphism.
\end{cor}

In what follows, we will refer to the substack $W^{\eta,+}_I\subset W^{\eta}_I$
as the {\em transverse locus}. By construction, a point $w\in W^\eta_I$
is transverse if and only if the following holds:

\begin{enumerate}
\item $w$ projects to a divisor $c^\eta_I = \Sigma \eta_k\cdot c_k$ with
each $\eta_k\in \CV(G/S)$,\
\item the associated map 
$D^\times_{c_k}\to X$ (defined up to $G(\CO_{c_k})$-conjugacy)
lies in $\bO_{\eta_k}$.
\end{enumerate}


\section{Relation to quasimaps}   \label{sect rel quasimaps}

\subsection{An intermediate stack}

We shall now introduce a stack, denoted $Z_{I,P^\op}$, that mediates
between the quasimaps space $Z_I$ and the local model $W_I$.

\medskip

Consider first the fiber product $Z_I \underset{\Bun_G}{\times} \Bun_{P^\op}$.
By definition, it classifies data 
$$(c_I\in C^I, \CP_P\in\Bun_P,\sigma:{C}\setminus {|c_I|}\to 
\CP_P{\overset{P}\times}X|_{{C}\setminus {|c_I|}})$$
where $\sigma$ is a section which factors
$$\sigma|_{C'}:C'\to 
\CP_P{\overset{P}\times}{\openX}|_{C'}\to
\CP_P{\overset{P}\times}{X}|_{C'},$$
for some open curve $C'\subset C\setminus|c_I|$.

Let $Z_{I,P^\op}$ be the open subset of the above fiber product
that corresponds to the condition that
$\sigma$ factors as 
$$\sigma|_{C'}:C'\to 
\CP_P{\overset{P}\times}{\openX^+}|_{C''}\to
\CP_P{\overset{P}\times}{X}|_{C''},$$
for some (possibly smaller) open subset $C''\subset C'$.

\subsection{Projection onto $Z_I$}

We have the obvious forgetful map $Z_{I,P^\op}\to Z_I$. Of course,
we do not claim that it is smooth. However, it will be smooth over a 
large enough open subset.

Let $M$ denote the Levi factor of $P^\op$. Define the open substack 
$\Bun_{M,r}\subset\Bun_M$ to be that for which
$$
H^1(C,\CP_{M}\overset{M}{\times}V)=0,
$$
for all $M$-modules $V$ which appear as subquotients of 
$\on{Lie}(U^\op)$. For any stack $\CQ$ mapping to $\Bun_M$, define the 
open substack $\CQ_r\subset\CQ$ 
to be the fiber product 
$$\CQ_r=\CQ\underset{\Bun_{M}}{\times}\Bun_{M,r}.$$

For $\mu\in\cowts_{M/[M,M]}$ let $\Bun_P^\mu$ be the corresponding connected
component of $\Bun_P$. 

\begin{lem} \hfill  \label{model}

\smallskip

\noindent(1) The map $\fr:\Bun_{P,r}\to Bun_G$ is smooth. 

\smallskip

\noindent(2) 
Any open substack 
of finite type $\Bun_G^{\on{fin}}\subset \Bun_G$,
is contained in the image of $\Bun^\mu_{P,r}$ for a
a sufficiently large $\mu\in\cowts_{M/[M,M]}$, 
and the restriction of $\fr$ to the inverse image
$\fr^{-1}(\Bun_G^{\on{fin}})\subset \Bun_{P,r}^\mu$
has connected fibers.
\end{lem}

These assertions are all well-known, except perhaps for the 
connectedness one, which is implied immediately by the following:

\begin{proof}
Using~\cite{BFG03}, we know there is a nonempty 
open subset $U_0\subset \Bun_G$ and coweight $\mu_0\in\cowts_G$
such that the restriction of the natural map $p:\Bun^{\mu_0}_B\to\Bun_G$ to the inverse image
$p^{-1}(U_0)\subset\Bun_B^{\mu_0}$ has connected fibers. 
Choose $\mu\in\cowts_G$ such that $U$ is in the image of the projection
$$
\Bun_G\stackrel{h^\la}{\leftarrow} \CH_G^\mu\underset{\Bun_G}{\times} U_0.
$$
Consider the diagram
$$
\begin{array}{ccccccc}
\Bun_B^{\mu+\mu_0} & \stackrel{h^\la}{\leftarrow} &
\CH_{B}^\mu\underset{\Bun_B}{\times}\Bun_B^{\mu_0}&
\subset &
\CH_{G}^\mu\underset{\Bun_G}{\times}\Bun_B^{\mu_0}
& \stackrel{h^\ra}{\to} & \Bun_B^{\mu_0}\\
\downarrow & & \downarrow & & \downarrow & & \downarrow  \\
\Bun_G & \stackrel{{}h^\la}{\leftarrow} & \CH^\mu_G & = & \CH^\mu_{G} & \stackrel{{}h^\ra}{\to} & \Bun_G.
\end{array}
$$
Since the map 
$$
\Bun_B^{\mu+\mu_0} \stackrel{h^\la}{\leftarrow} 
\CH_{B}^\mu\underset{\Bun_B}{\times}\Bun_B^{\mu_0}
$$
is a bijection,
to prove the lemma, it suffices to show that for a point $x\in U\subset\Bun_G$, 
the inverse image 
$$
(h^\la)^{-1}(p^{-1}(x))\subset
\CH_{G}^\mu\underset{\Bun_G}{\times}\Bun_B^{\mu_0}
$$
is connected.
We know that $(h^\la)^{-1}(x)\subset\CH^\mu_G$ is connected.
Furthermore, since
the rightmost square of the diagram is Cartesian,
we know that the restriction of the map
$$
\CH_{G}^\mu\underset{\Bun_G}{\times}\Bun_B^{\mu_0}
\to
\CH_G^\mu
$$
to the inverse image of $(h^\ra)^{-1}(U_0)$ has connected fibers.
By construction, we have that the intersection
$$
(h^\la)^{-1}(x)\cap (h^\ra)^{-1}(U_0) 
$$ 
is nonempty and dense in $(h^\la)^{-1}(x)$,
and so we are done.
\end{proof}

Consider the corresponding open subset
$$Z_{I,P^\op,r}^\mu\subset Z_{I,P^\op}^\mu
= Z_I\underset{\Bun_G}\times \Bun_P^\mu.
$$
Applying base change with respect to $Z_I\to \Bun_G$, we obtain that the
statements of the above lemma apply to the forgetful morphism
$$Z_{I,P^\op,r}^\mu\to Z_I.$$
Thus the space $Z_{I,P^\op}$ can be used for the local study of the space
$Z_I$.

\subsection{Passage to the local model}

One relationship between $Z_{I,P^\op}$ and $W_I$ is the following obvious one.
Consider the map $\Bun_{A^0}\to \Bun_M$ corresponding to our choice
of the splitting $A_0\hookrightarrow M$. For $\eta\in \cowts_{A_0}$,
let $\mu$ be its image in $\cowts_{M/[M,M]}$.
Then we have a Cartesian diagram
$$
\CD
W_I^\eta @>>> Z^\mu_{I,P^\op} \\
@VVV    @VVV  \\
\Bun^\eta_{A_0}  @>>>  \Bun^\mu_M.
\endCD
$$

But as we now explain, in fact $Z^\mu_{I,P^\op}$ is equivalent in the smooth topology to
$W^\eta_I$. 
Recall that $Q$ denotes the subgroup $S\cap P^\op$, so we have a short exact sequence
$$1\to Q\to M\to A\to 1.$$
Let $Q^0$ denote the connected component of $Q$, which equals the kernel of
$M\to A_0$. Our choice of the splitting $A_0\hookrightarrow M$ gives rise to a decomposition
$$Q\simeq Q^0\times F.$$ 
Because of our intended application,
we state the following for $Q^0$ though its assertion and proof use nothing
special about $Q^0$.

\begin{lem}
Let $\Bun_{Q^0}^{\on{fin}}$ be an open substack of finite type of $\Bun_{Q^0}$,
and let $\fc \subset C$ be a finite subset.
Then there exists a scheme $\CY$ with a smooth surjective map
$\CY\to \Bun_{Q^0}^{\on{fin}}$ with connected fibers, such that the 
corresponding universal
bundle on $C\times \CY$ is trivialized on a Zariski-open subset of
$C\times \CY$ containing $\fc\times \CY$.
\end{lem}

\begin{proof}
Let $B\subset Q^0$ be a Borel subgroup, let $N\subset B$ be its unipotent radical,
 and let $B\to B/N=T$ be its Cartan quotient.
As we have seen in \lemref{model}, we can replace $\Bun_{Q^0}$ by
$\Bun^\nu_B$ for a large enough coweight $\nu\in \Lambda_T$. Namely, we can choose $\nu$ large enough 
so that there is
an open $\CV\subset \Bun_B^\nu$ that maps smoothly
onto $\Bun_{Q^0}^{\on{fin}}$ with connected fibers.

Consider the canonical map $\Bun_B\to \Bun_T$, choose a scheme
$\CY_{T}$ satisfying the assertion of the lemma for $\Bun_T$, and set 
$$
\CY=\CY_{T}\underset{\Bun_T}\times \CV.
$$
Let $\CW$ be an open subset of $\CY\times C$ over which the
$T$-bundle is trivialized. We may assume that $\CW$ is affine
over $\CY$. But then the whole $B$-bundle is trivialized, since
an $N$-bundle over an affine is trivial.
\end{proof}

Observe that the lemma generalizes to families. Namely, if we choose an $\CS$-family
$\fc \subset C\times\CS$ finite over $\CS$, then there is a scheme $\CY_\CS$ over $\CS$ 
with a smooth surjective map
$\CY_\CS \to \Bun_{Q^0}^{\on{fin}}$ with connected fibers, such that the 
corresponding universal
bundle on $C\times \CY_\CS$ is trivial on a Zariski-open subset of
$C\times \CY_\CS$ containing the inverse image of $\fc$.

Note as well that the above construction of $\CY_\CS$ involves the choice
of a coweight $\nu\in \Lambda_{Q^0}$. We write $\CY^\nu_\CS$ to denote this dependence.

Consider an $\CS$-family $\fc\subset C\times\CS$ finite over $\CS$,
and let $Z_{I,P^{op},\fc}^\mu$, respectively $W^\eta_{I ,\fc}$, be the $\CS$-family of $Z_{I,P^{op}}^\mu$,
respectively $W^\eta_I$, of maps whose
degeneracy locus belongs to $\fc$.

\begin{lem}\label{base change to local model}
For $\mu=\eta+\nu$ in $\Lambda_{M/[M,M]}$,
we have an identification of $\CS$-families
$$
Z_{I,P^{op},\fc}^\mu  \underset{\Bun_{Q^0}\times \CS}{\times} \CY^\nu_\CS \simeq W^\eta_{I,\fc} 
\underset{\CS}{\times}\CY^\nu_\CS.
$$
\end{lem}

\begin{proof}
We define a map from the left hand side to the right as follows.

Let $\CW_\CS$ be the Zariski open subset of $C\times \CY_\CS$ on which the universal
$Q^0$-bundle is trivialized.
The base change of the left hand side provides a reduction of the $P^{op}$-bundle to an
$R$-bundle over $\CW_\CS$. 

On the other hand, let $\CY_\fc\subset \CY_\CS$ be the inverse image
of $\fc$, and let $\CY_\CS \setminus \CY_\fc$ be its complement.
The quasimap data of the left hand side equips the $P^{op}$-bundle with a reduction to 
$Q\simeq Q^0\times F$
over $\CY_\CS \setminus \CY_\fc$.

By construction, these two reductions are compatible over the intersection 
$\CW_\CS\cap (\CY_\CS \setminus \CY_\fc)$ in the sense that  
we may define a global $R$-bundle by gluing the $R$-bundle over $\CW_\CS$
with the $R$-bundle induced from the $F$-bundle 
over $\CY_\CS \setminus \CY_\fc$.
The quasimap data of the right hand side is simply the tautological reduction 
over $\CY_\CS \setminus \CY_\fc$ of the 
induced $R$-bundle
to the $F$-subbundle from which it was induced.
\end{proof}

We conclude that the space $W^\eta_{I}$ can be used for the local study of the space
$Z^\mu_{I,P^{op}}$.


\subsection{Behavior with respect to the transverse locus}
Consider the data $(\CP_{P^\op},\sigma)$ of a quasimap on the formal disk:
$\CP_{P^\op}$ is a $P^\op$-bundle on $D$, and $\sigma$ is a map
$$
D^\times \to \CP_{P^\op}\overset{P^\op}\times\openX^+.
$$
Alternatively, the data can be interpreted as a point of
the fiber product
$$\on{Maps}(D,pt/P^\op)\underset{\on{Maps}(D^\times,pt/P^\op)}\times
\on{Maps}(D^\times,\openX^+/P^\op).$$

On the one hand, using the projections $P^\op\to A$ and $\openX^+\to A$, we 
see that $(\CP_{P^\op},\sigma)$ gives rise to a point of
$\affgr_A$. We interpret this an element $\eta \in\cowts_A$.

On the other hand, $(\CP_{P^\op},\sigma)$ defines a $P^\op(\CO)$-orbit
$\bO'_{P^\op}$ in $\openX^+(\CK)$, and hence a $G(\CO)$-orbit
$\bO_\lambda$ in $\openX(\CK)$, for some $\lambda\in \CV(G/S)$.

By \secref{param of orbits}, the following conditions are equivalent:
\begin{enumerate}
\item $\eta=\lambda\in \cowts_A$.\
\item $\bO'_{P^\op}$ is open in $\bO_\lambda$.
\end{enumerate}

If these conditions are satisfied, we say that $(\CP_{P^\op},\sigma)$
is {\it transverse}. We say that a point of the stack $Z_{I,P^\op}$
is transverse at $c$ (resp., globally transverse), if its restriction to the 
formal disk around $c$ (resp., for every $c$) is transverse.

It is clear 
that for a given point of $z_I\in Z_I$, the locus of points
in the preimage
$$
\fr^{-1}(z_I)\subset Z_{I,P^\op}
$$
that are transverse at a given finite collection of points of $C$
is open. Moreover, $\fr^{-1}(z_I)\cap Z_{I,P^\op,r}^\mu$ is non-empty 
for $\mu$ large enough.

\subsubsection{}

Recall now that in \secref{sect trans locus} we introduced the notion
of transversality for points of $W^\eta_I$. We claim that this notion
matches with the above notion of global transversality 
under the identification of \lemref{base change to local model}:
$$
Z_{I,P^{op},\fc}^\mu  \underset{\Bun_{Q^0}\times \CS}{\times} \CY^\nu_\CS \simeq W^\eta_{I,\fc} 
\underset{\CS}{\times}\CY^\nu_\CS
$$

Indeed, suppose that the pair $(\CP_{P^\op},\sigma)$ as above is such that
$\CP_{P^\op}$ is induced from an $A_0$-bundle $\CP_{A_0}$. As in
\secref{descr image}, such data give rise to an $F$-bundle $\CP_F$ on
$D^\times$ and a point of 
$\affgr^{\eta}_{R,\CP_F}\cap \affgr_{G,\CP_F}^{S,\lambda}$.
We see that in both cases, the transversality condition is that $\eta=\lambda$.

We denote by $Z_{I,P^\op}^{\mu,+}\subset Z_{I,P^\op}^{\mu}$
the locus of transverse points. By the above, this is a closed
ind-substack.

\subsubsection{Transversality in the horospherical case}

Let us assume that $X$ is horospherical. Recall that
$\Bun_P$ is naturally a locally closed substack of $\stZ_I$ and for
any $\eta$ we have an ind-locally closed embedding
$$j_I:C^\eta_I\times \Bun^{\mu'}_P\hookrightarrow \stZ^{\mu'+\eta}_I.$$

Let $\stZ_{I,P^\op}^{\mu,\mu'}$ denote the analogue of the stack
$Z_{I,P^\op}^\mu$ when we use $\stZ_I^{\mu'}$ instead of $Z_I$. 
Note that the fiber product
$$\Bun^{\mu'}_P\underset{\stZ^{\mu'}_I}\times \stZ_{I,P^\op}^{\mu'}\subset
\Bun_P^{\mu'}\underset{\Bun_G}\times \Bun^{\mu'}_{P^\op}$$
identifies naturally with $\Bun^{\mu'}_M$.

From $j_I$ we obtain a map
$$C^\eta_I\times \Bun^\mu_M\to \stZ_{I,P^\op}^{\mu+\eta,\mu}.$$
The following is easy to check.

\begin{lem}
The above map is an isomorphism on the transverse locus
$$\stZ_{I,P^\op}^{\mu+\eta,\mu,+}\subset \stZ_{I,P^\op}^{\mu+\eta,\mu}.$$
It induces an isomorphism
$$C^\eta_I \underset{\Bun_{A}}\times \Bun^\mu_M\simeq
Z_{I,P^\op}^{\mu+\eta,+},$$
\end{lem}

The second map of the lemma induces a map
$$M^\eta_I\to W^{\eta,+}_I,$$
which is the inverse of the one of \corref{deep stratum}.

\newpage


{\Large \part{Proofs--B}}

\vspace*{5mm}

In this part we will apply the local model
$W^\eta_I$ of \secref{sect factor} to prove two main technical assertions
of Part 1, namely \thmref{thmidentifyconv} (in \secref{proof of conv}) and
\thmref{thmconv} (in \secref{proof of special}).

\vspace*{5mm}


\section{Convolution and local model}  \label{proof of conv}

\subsection{Proof of \thmref{thmidentifyconv}} 

Our proof is organized into several steps, each involving the use of more
specialized structure than the previous.
Recall that since $\pshonaffgr^{\otimes I}$ is semisimple,
by the decomposition theorem, the convolution $\conv_I(\CA)$ is semisimple as well. 
Our goal is to see
that each of its summands is the intersection cohomology sheaf with constant
coefficients
of a connected component of an untwisted local stratum ${}' Z_{\mr I}^\Theta \subset {}'Z_{\mr I}$.
 
\subsubsection{Reformulation of convolution}
We begin with the trivial observation that since $\pshonaffgr^{\otimes I}$ is semisimple,
it suffices to prove the theorem for simple objects. These are the intersection cohomology
sheaves 
$
\CA^{\lambda_I}_G 
$
of the product of $G(\CO)$-orbit closures
$$
\ol{\affgr}_G^{\lambda_I}=\prod_{i\in I}\ol{\affgr}_G^{\lambda_i}\subset(\affgr_G)^I,
\mbox{ where $\lambda_I:I\to\domcowts_G$.}
$$

This allows for the following
slight reformulation of the convolution.
By restricting the $(G(\CO)\rtimes\Aut(\CO))^I$-torsor
$$
\wh Z_{\mr I}\to Z_{\mr I}
$$
to the closure $'\ol{Z}_{\mr I}^0\subset {}Z_{\mr I}$ of the untwisted basic stratum, 
we obtain a $(G(\CO)\rtimes\Aut(\CO))^I$-torsor
$$
'\wh{\ol{Z}}^0_{\mr I}\to {}'\ol{Z}^0_{\mr I}.
$$
The twisted product
$$
'\wt{\ol{Z}}^{\lambda_I}_{\mr I}
=
{}'{\wh{\ol Z}}^0_{\mr I}\overset{(G(\CO)\rtimes\Aut(\CO))^I}{\times}\ol{\affgr}_G^{\lambda_I}
$$
fits into a diagram
$$
'Z_{\mr I}  
\stackrel{ h^\leftarrow}{\leftarrow} 
{}'\widetilde {\ol{Z}}^{\lambda_I}_{\mr I}
\stackrel{h^\rightarrow}{\to}  
{}'\ol{Z}^0_{\mr I} 
$$
in which $h^\rightarrow$ is the evident projection,
and $h^\la$ is the modification map.
In other words, the diagram results from restricting the diagram
$$
Z_{\mr I}  
\stackrel{  h^\la}{\leftarrow} 
\widetilde {{Z}}_{\mr I}
\stackrel{ h^\ra}{\to}  
{Z}_{\mr I}
$$
along the inclusion
$$
'\widetilde {\ol{Z}}^{\lambda_I}_{\mr I}
\hookrightarrow
\wt Z_{\mr I}.
$$

In the category $\catp(\wt Z_{\mr I})$,
we have a canonical isomorphism
$$
'\wt{{\IC}}_{\mr I}^{\lambda_I}\simeq 
{}'\IC_{\mr I}^0
\tboxtimes
\CA^{\lambda_I}_G 
$$
identifying 
the intersection
cohomology sheaf of $'\widetilde{\ol Z}_{\mr I}^{\lambda_I}$
with the twisted product.
Thus we have a canonical isomorphism
$$
\conv_I(\CA^{\lambda_I}_G)\simeq
\sum_k\pH^k(h^\leftarrow_{\lambda_I !}({}'\wt{{\IC}}_{\mr I}^{\lambda_I})).
$$
We will prove the theorem using this formulation of convolution.

\subsubsection{Moving the pole points}
Suppose we have established the following.

\begin{claim}\label{basedclaim}
For
any fixed configuration of pole points $c_I\in C^I$, 
the restriction of each irreducible summand of the convolution $\conv_I(\CA^{\lambda_I}_G)$
to the based ind-substack $'Z_{c_I}\subset {}'Z_{I}$ is the middle-extension of the constant sheaf on
a component of some untwisted based local stratum 
$'Z^{\Theta}_{c_I}\subset {}'Z_{c_I}$.
\end{claim}

Then we may complete the proof of \thmref{thmidentifyconv} as follows.
Observe that the closure of the untwisted basic stratum is a product
$$
'\overline{Z}_{\mr I}^0 = {}'\overline {Z}_\emptyset^0\times \mr C^I,
$$
and we have a commutative diagram
$$
\begin{array}{ccc}
'\ZdistinctI & \stackrel{ h^\la}{\leftarrow} &'\wt Z_{\mr I}  \\
\hspace{-1em}\pi \downarrow & & \hspace{1em}\downarrow  h^\ra \\
\mr C^I & \stackrel{p}{\leftarrow} & '\overline {Z}^0_\emptyset\times \mr C^I 
\end{array}
$$
where the map $p$ is the obvious projection.
We have an isomorphism
$$
h^\ra_!({}'\IC^0\tboxtimes\CA^{\lambda_I}_G)\simeq
{}'\IC^0_\emptyset\boxtimes V^{\lambda_I}_{gr}
$$
where we write $V^{\lambda_I}_{gr}$ for the graded constant sheaf on $\mr C^I$
with fiber the tensor product of the irreducible representations of $\ch G$
of highest weights $\lambda_i\in\domcowts_G$, for $i\in I$,
with canonical grading given by the action of the principal nilpotent.
Thus we have an isomorphism
$$
\pi_! h^\la_!({}'\IC^0\tboxtimes\CA^{\lambda_I}_G)
\simeq
p_!({}'\IC^0_\emptyset\boxtimes V^{\lambda_I}_{gr}).
$$

We conclude that the pushforward via $\pi$ of the convolution $\conv_I(\CA^{\lambda_I}_G)$ 
is a constant sheaf. 
Therefore the
top cohomology of the pushforward
of each of its irreducible summands is also a constant sheaf
(over the generic point of the base it is a sublocal system of a constant local system).
Coupled with Claim~\ref{basedclaim}, this establishes \thmref{thmidentifyconv}.
Thus our aim in what follows is to prove Claim~\ref{basedclaim}.


\subsubsection{Convolution for local model}
We begin by constructing a version of convolution for
the local model. 
Consider the $(G(\CO)\rtimes\Aut(\CO))^I$-torsor
$$
\wh W^{\eta}_{\mr I}\to W^{\eta}_{\mr I}
$$
that classifies data
$$
(w,\mu,\tau)
$$
where $w\in W^{\eta}_{\mr I}$, with pole points $c_I\in \mr C^I$
and $R$-bundle $\CP_{R}\in\Bun_{R}$, 
$\mu$ is an isomorphism of $G$-bundles
$$
\mu:D_{|c_I|}\times G\risom \CP_{R}\overset{R}{\times} G|_{D_{|c_I|}},
$$
and $\tau$ is an isomorphism
$$
\tau:D\times |c_I|\risom D_{|c_I|}.
$$
The twisted product
$$
\wt W^{\eta}_{\mr I}=\widehat W^{\eta}_{\mr I}
\overset{(G(\CO)\rtimes\Aut(\CO))^I}{\times}
(\affgr_G)^I
$$ 
classifies data
$$
(w,\CP_G,\alpha)
$$
where $w\in W^{\eta}_{\mr I}$, with pole points $c_I\in \mr C^I$ and
$R$-bundle $\CP_{R}\in\Bun_{R}$,
$\CP_G\in\Bun_G$, and $\alpha$ is an isomorphism of $G$-bundles
$$
\alpha:\CP_{R}\overset{R}{\times} G|_{C\setminus{|c_I|}}
\risom
\CP_G|_{C\setminus{|c_I|}}.
$$

To a point $c\in C$, and a point $\tilde w\in\wt W^{\eta}_{\mr I}$,
with pole points $c_I\in \mr C^I$, $R$-bundle $\CP_{R}\in\Bun_{R}$,
section
$$
\sigma:C\setminus|c_I|\to\CP_{R}\overset{R}{\times} X|_{C\setminus|c_I|},
$$
$G$-bundle $\CP_G\in\Bun_G$, and isomorphism 
$$
\alpha:\CP_{R}\overset{R}{\times} G|_{C\setminus{|c_I|}}
\risom
\CP_G|_{C\setminus{|c_I|}},
$$
we may associate two elements
$\bar\sigma_c^\leftarrow,\bar\sigma_c^\rightarrow\in \CV(\mr X)$ 
as follows.
For some open curve $C'\subset C\setminus|c_I|$, the section $\sigma$ factors
$$
\sigma|_{C'}:C'\to 
\CP_{R}{\overset{R}{\times}}{\mr X}|_{C'}\to
\CP_{R}{\overset{R}{\times}}{X}|_{C'}.
$$
For any trivialization of the restriction $\CP_{R}|_{D_c}$,
we obtain a projection 
$$
\CP_{R}\overset{R}{\times} \mr X|_{D^\times_c}\to \mr X.
$$
For any choice of formal coordinate at $c$, 
we define $\bar\sigma_c^\leftarrow$ to be the composition of
the resulting identification $D^\times\risom D^\times_c$ with the restriction
$$
\sigma|_{D^\times_c}:D^\times_c \to\CP_{R}\overset{R}{\times} \mr X|_{D^\times_c}\to \mr X.
$$
Similarly, for any trivialization of the restriction $\CP_{G}|_{D_c}$,
we obtain a projection 
$$
\CP_{G}\overset{G}{\times} \mr X|_{D^\times_c}\to \mr X.
$$
For any choice of formal coordinate at $c$, 
we define $\bar\sigma_c^\rightarrow$ to be the composition of
the resulting identification $D^\times\risom D^\times_c$ with the restriction of the composition
$$
\alpha\circ\sigma|_{D^\times_c}:
D^\times_c \to\CP_{R}\overset{R}{\times} \mr X|_{D^\times_c}\stackrel{}{\risom}
\CP_{G}\overset{G}{\times} \mr X|_{D^\times_c}\to \mr X.
$$
By Corollary~\ref{corparam}, the resulting elements
$\bar\sigma_c^\leftarrow,\bar\sigma_c^\rightarrow\in (X(\CK) - (X(\CK) - \openX(\CK)))/G(\CO)$
are independent of the choices.

For a labelling $\lambda_I:I\to\domcowts_G$, 
we define the opposite labelling
$\lambda_I^{{op}}:I\to\domcowts_G$ 
to be that for which $\lambda^{{op}}_i$ is the dominant coweight
in the Weyl group orbit through $-\lambda_i$, for $i\in I$.
We define the ind-substack 
$$
'\widetilde{ W}^{\eta,\lambda_I}_{\mr I}
\subset
{}'\widehat W^{\eta}_{\mr I}
\overset{(G(\CO)\rtimes\Aut(\CO))^I}{\times}
\ol{\affgr}_G^{\lambda^{{op}}_I}
$$
to consist of those points such that 
the associated element
$\bar\sigma_c^\rightarrow\in\CV(\mr X)$ is trivial for all $c\in C$,
and we define the ind-substack
$$
{}'\widetilde{\ol{ W}}^{\eta,\lambda_I}_{\mr I}
\subset
{}'\widehat W^{\eta}_{\mr I}
\overset{(G(\CO)\rtimes\Aut(\CO))^I}{\times}
\ol{\affgr}_G^{\lambda^{{op}}_I}
$$
to be its closure. 
We have the evident projection
$$
'{ W}^{\eta}_{\mr I}
\stackrel{}{\leftarrow}
{}'\widetilde{\ol W}^{\eta,\lambda_I}_{\mr I}
$$
which we denote by ${h^\leftarrow_{\lambda_I}}{}$.

Recall that $\pshonaffgr^{\otimes I}$ is semisimple.
We define the local version of the convolution functor
$$
\conv^{\eta}_I: \pshonaffgr^{\otimes I}\to\catp({}' W^\eta_{\mr I})
$$
on simple objects by the formula
$$
\conv^{\eta}_I(\CA_G^{\lambda_I})=\sum_{k}
\pH^k(h^\leftarrow_{\lambda_I !}({}'\wt{{\IC}}^{\eta,\lambda_I}_{\mr I}))
$$
where $'\wt{{\IC}}^{\eta,\lambda_I}_{\mr I}$ denotes the intersection cohomology sheaf
of $'\wt{\ol{ W}}^{\eta,\lambda_I}_{\mr I}$.


\subsubsection{Relating the convolutions}
Recall that in \secref{sect rel quasimaps},
we constructed a correspondence relating the quasimaps space
and the local model
$$
W^\eta_{I,\fc}
\leftarrow
W^\eta_{I,\fc} \underset{\CS}{\times}\CY^\nu_\CS
\simeq
Z_{I,P^{op},\fc}^\mu  \underset{\Bun_{Q^0}\times \CS}{\times} \CY^\nu_\CS 
\to
Z_{I,P^{op},\fc}^\mu 
\to
Z_{I} 
$$

Given a a labelling $\lambda_I:I\to\domcowts_G$, 
we may extend the convolution maps of the preceding two sections
to a Cartesian diagram
$$
\begin{array}{ccccccccc}
'\wt{\ol W}^\eta_{I,\fc}
& \leftarrow &
'\wt{\ol W}^\eta_{I,\fc} \underset{\CS}{\times}\CY^\nu_\CS
& \simeq &
'\wt{\ol Z}_{I,P^{op},\fc}^\mu  \underset{\Bun_{Q^0}\times \CS}{\times} \CY^\nu_\CS 
& \to & 
'\wt{\ol Z}_{I,P^{op},\fc}^\mu 
& \to &
'\wt{\ol Z}_{I} \\
 \downarrow &  & \downarrow &  & \downarrow &  & \downarrow & & \downarrow \\
'W^\eta_{I,\fc}
& \leftarrow &
'W^\eta_{I,\fc} \underset{\CS}{\times}\CY^\nu_\CS
& \simeq &
'Z_{I,P^{op},\fc}^\mu  \underset{\Bun_{Q^0}\times \CS}{\times} \CY^\nu_\CS 
& \to & 
'Z_{I,P^{op},\fc}^\mu 
& \to &
'Z_{I} 
\end{array}
$$

We will use the following properties of this diagram. They all follow from the results of
\secref{sect rel quasimaps}:

\begin{itemize}

\item For large enough parameters, the horizontal maps are all smooth.

\smallskip

\item For fixed pole points $c_I$, given any based local stratum $'Z_{c_I}^\Theta$, 
we may arrange so that the corresponding stratum of $'W_{c_I}^{\eta,\Theta}$
is transverse at the pole points.

\smallskip

\item The fibers of the rightward pointing horizontal maps are connected.

\end{itemize}

Since intersection cohomology sheaves are preserved by smooth pullback,
to calculate the convolution locally in the smooth topology,
it suffices to understand what happens for the local model.
By the last point, to calculate how the cohomology sheaves of the convolution might twist along a locus of constructibility,
again it suffices to understand what happens for the local model.


\subsubsection{Proof of Claim~\ref{basedclaim}}
By the discussion immediately preceding, to prove Claim~\ref{basedclaim},
it suffices to prove the following assertion for the local model.

\begin{prop}\label{based local prop}
For any fixed configuration of pole points $c_I\in C^I$, 
the restriction of the convolution 
to an untwisted based local stratum 
$'W^{\eta, \Theta}_{c_I}\subset {}'W^\eta_{c_I}$ which is transverse at its pole points
is constant.
\end{prop}

If we were only attempting to prove that the restriction were {\em locally constant},
we could simply appeal to the factorization pattern of the local model. 
Namely,
the factorization diagram for $W^\eta_{c_I}$ naturally extends to a factorization
diagram for its convolution. 
Furthermore, the local version of the convolution and its factorization readily generalize
to the local model
$W^\eta_{c_I,1}$
over the open curve $C\setminus c$.
Using this added freedom, we could reduce the assertion of the proposition to the following cases:
\begin{enumerate}
\item the restriction of $\conv_{i,1}^{\eta_i}(\CA_i)$ to the untwisted based local stratum
$ 'W^{\eta_i, \eta_i}_{c_i,1}\subset {}'W^{\eta_i}_{c_i,1}$
which is transverse 
at its single pole point $c_i$ and has no other degeneracies,
\item
the restriction of the basic sheaf ${}'\IC^0$ to the smooth locus of the  
local model ${}' W^{\eta'}_{\emptyset,1}$
with no pole points but otherwise arbitrary degeneracies.
\end{enumerate}

In the second case, the restriction of the intersection cohomology sheaf
to the smooth locus is constant.
In the first case, the stratum $'W^{\eta_i, \eta_i}_{c_i,1}$
reduces to a copy of the stack of $F$-local systems on $C\setminus (c_i \cup c)$
whose induced $\pi_0(S)$-local system is trivializable.
Thus each of its connected component  is isomorphic to the classifying space $pt/F$.
Since all sheaves on the classifying space are locally constant, we would be done.
Unfortunately, to factorize the local model, we must make base changes with disconnected
fibers. Thus in the course of the above argument, we could potentially
unwind local systems. 
Because we seek the stronger statement that our sheaves are in fact constant, 
we will proceed with a slightly modified approach.

\begin{proof}
First, let $\theta$ be the total degree of $\Theta$,
and let $\eta'=\eta - \theta$
be the total excess degeneracies. We base change to the open subset
$$
W^{\theta,\eta'}_{c_I}=
W^\eta_{c_I} \underset{C^\eta}{\times} (C^{\theta} \times C^{\eta'})_\disj 
$$
where we separate apart the excess degeneracies.

Next, let $C^{\theta}_{I,ex}$ be the space of pairs $(c^\theta_I,e)$
of disjoint divisors on $C $, 
with the first $c^\theta_I \in C^\theta_I$.
Define the space $W^{\theta}_{c_I,ex}$ to be that classifying the data
of a pair of disjoint divisors 
$
(c^\theta_I, e)\in C^{\theta}_{I,ex}
$
together with the usual local model data
of total degree $\theta$ on the open curve $C\setminus e$
with associated divisor $c^\theta_I$.
We have the obvious projection
$$
W^{\theta,\eta'}_{c_I}\to W^{\theta}_{c_I,ex}
$$
where we only keep track of the fibers of the local model
at the degeneracy points of the first type.
Moreover,
we clearly may extend the convolution maps so that we have a Cartesian diagram
$$
\begin{array}{ccc}
\wt{\ol W}^{\theta,\eta',\lambda_I}_{c_I} & \to &  \wt{\ol W}^{\theta,\lambda_I}_{c_I,ex}\\
\downarrow & & \downarrow \\
W^{\theta,\eta'}_{c_I} & \to & W^{\theta}_{c_I,ex}
\end{array}
$$
where the first vertical map is simply the base
change to $(C^{\theta}_{I} \times C^{\eta'})_\disj$ of the usual local convolution map.
Thus we see that to prove the proposition, it suffices to prove it for the convolution
given by the second vertical map.

Fix a point $c\in C$ disjoint from the pole points, and consider the base change of
the previous constructions to the open curve $C\setminus c$.
Let $(C\setminus c)^{\theta}_{I,ex}$ be the space of pairs $(c^\theta_I,e)$
of disjoint divisors on $C\setminus c$, 
with the first $c^\theta_I \in (C\setminus c)^\theta_I$.
Let $ W^{\theta}_{c_I,ex,1}$ be the space classifying the data
of a pair of disjoint divisors 
$
(c^\theta_I, e)\in (C\setminus c)^{\theta}_{I,ex}
$
together with the usual local model data
of total degree $\theta$ on the open curve $C\setminus (c\cup e)$
with associated divisor $c^\theta_I$.

Define the group scheme
$$
\CF^{\theta}_{ex}
\to (C\setminus c)^{\theta}_{I,ex}
$$ 
to be that  
classifying the data of a pair 
$$(c^\theta_I,e)\in (C\setminus c)^{\theta}_{I,ex},
$$
an $F$-bundle $\CP_F$ on the complement $C\setminus (c \cup e)$, and
 a trivialization of $\CP_F$ above the divisor $|c^\theta_I|\subset C\setminus (c\cup e)$. 
The group structure is given
by tensor product of $F$-bundles. 

Observe that the relative group $\CF^{\theta}_{ex}$ naturally
acts on the convolution map
$$
\begin{array}{ccc}
\wt{\ol W}^{\theta,\lambda_I}_{c_I,ex,1} \underset{(C\setminus c)^{\theta}_{I,ex}}{\times} \CF^{\theta}_{ex}
& \to &  \wt{\ol W}^{\theta,\lambda_I}_{c_I,ex,1}\\
 \downarrow & & \downarrow \\
W^{\theta}_{c_I,ex,1} \underset{(C\setminus c)^{\theta}_{I,ex}}{\times} \CF^{\theta}_{ex}
& \to &  W^{\theta}_{c_I,ex,1}
\end{array}
$$
To establish the proposition, we will prove the a priori stronger statement that the restriction
of the convolution to an untwisted transverse based local stratum is constant as an equivariant
object for $\CF^{\theta}_{ex}$.
But note that the action of $\CF^{\theta}_{triv}$ 
is free so that the stronger assertion
follows from the weaker analogue.

Let $\CF^{\theta}_{triv}\subset \CF^{\theta}_{ex}$ be the subgroup
of those $F$-bundles with trivial monodromy around any point in $C\setminus c$. Similarly, 
let $W^{\theta}_{c_I,triv, 1}\subset W^{\theta}_{c_I,ex, 1}$ be the subspace
of points whose associated generic $F$-bundle has trivial monodromy around the excess divisor
in $C\setminus c$.
Then it is easy to see that inclusion induces an isomorphism of quotient stacks
$$
W^{\theta}_{c_I,triv, 1}/\CF^{\theta}_{triv}
\risom 
W^{\theta}_{c_I,ex, 1}/\CF^{\theta}_{ex}.
$$
Thus to prove the stronger assertion, it suffices to prove it for the first quotient
of the above isomorphism. 
So we are left to see that the restriction
of the convolution to an untwisted transverse based local stratum of 
$W^{\theta}_{c_I,triv, 1}$ is constant. Observe that here the excess divisor
is playing no role. In other words, we may forget it completely to obtain a Cartesian diagram
$$
\begin{array}{ccc}
 \wt{\ol W}^{\theta}_{c_I,triv, 1} & \to &  \wt{\ol W}^{\theta,\lambda_I}_{c_I,1} \\
\downarrow & & \downarrow \\
W^{\theta}_{c_I,triv, 1} & \to &  { W}^{\theta,\lambda_I}_{c_I,1}
\end{array}
$$
with smooth horizontal maps.
Now, each connected component of an untwisted transverse based local stratum of 
$W^{\theta}_{c_I, 1}$ is nothing more than the classifying space $pt/F$. Thus we are
left to see that the $F$-action on the fiber of the convolution is trivial. But the $F$-action
factors through $F\to \pi_0(S)$ and all of our constructions can be made on the closed
subspace of untwisted quasimaps, i.e. those quasimaps with trivial associated $\pi_0(S)$-bundle.
Thus the $F$-action is trivial and we are done.
\end{proof}


\subsection{Proof of  \corref{corconvglue}} The corollary will follow from
the following refinement of \thmref{thmidentifyconv}.

\begin{prop}\label{propfactorizefibers}
For any object $\CA_I\in\pshonaffgr^{\otimes I}$, if we fix isomorphisms
$$
\conv(\CA_i)\simeq\sum_{\theta\in\domcowts_{\openX}}{}'\IC^\theta_{}(V^\theta_i),
$$
where $V^\theta_i$ are vector spaces,
then we have a canonical isomorphism
$$
\conv_I(\CA_I)\simeq\sum_{\Theta:I\to\domcowts_{\openX}}{}'\IC^\Theta_{\mr I}(\otimes_{i\in I}V^{\Theta(i)}_i).
$$
\end{prop}

\begin{proof}
Fix a labelling $\Theta:I\to\domcowts_{\openX}$, and 
a point $z\in {}'Z_{\mr I}^\Theta$. 
By our previous results,
to calculate the stalk of $\conv_I(\CA_I)$ at $z$, we may instead consider 
the stalk of $\conv_I^\eta(\CA_I)$
at a point $w\in{}' W_{\mr I}^{\eta,\Theta}$
such that $w$ is transverse at its pole points.
Then via the factorization pattern for the convolution discussed in the previous section, 
this stalk is the tensor product of the following:
\begin{enumerate}
\item the stalk of $\conv_{i,1}^{\eta_i}(\CA_i)$ at 
a point $w_i\in {}'W^{\eta_i}_{i,1}$
with a degeneracy of degree $\eta_i$
at its single pole point $c_i$ and no other degeneracies,
\item
the stalk of the basic sheaf ${}'\IC^0$ at a point $w'\in {}' W^{\eta'}_{\emptyset,1}$
with no pole points but otherwise arbitrary degeneracies.
\end{enumerate}

In the second case, we are taking the stalk of the intersection cohomology sheaf
at a smooth point, and so we may canonically identify the stalk with the trivial vector space $\BC$. 
But as in the previous proposition, 
$\conv_{i,1}^{\eta_i}(\CA_i)$ is constant along the locus appearing in the first case.
\end{proof}


\section{Specialization of sheaves}  \label{proof of special}

\subsection{Proof of \thmref{thmconv}(1)}    \label{proof of thmconv(1)} \hfill
By construction, we have a surjection of abelian groups
$
\pi_0(S_0)\to\pi_0(S).
$
We denote the kernel of this projection by $ker$.
Consider the closed substack of the basic stratum
$$
''Z^0_{0,\mr I}\subset Z^0_{0,\mr I}
$$ 
consisting 
of quasimaps whose generic $\pi_0(S_0)$-bundle is induced from a $ker$-bundle.
We refer to it
as the {\em partially twisted basic stratum}.

\begin{lem}\label{lconstr}
The restriction of the nearby cycles $\psi('\IC_{\mr I}^{0})$ 
to the basic stratum $Z_{0,\mr I}^0\subset Z_{0,I}$
is isomorphic to the constant sheaf of rank one
on the partially twisted basic stratum $''Z_{0,\mr I}^0\subset Z^0_{0,\mr I}$.
\end{lem}

\begin{proof}
Choose a point in the basic horospherical stratum $Z_{0,\mr I}^0$
and identify a neighborhood of it in the family $\CZ_{I}$ with a neighborhood of
some point in a local model $\CW^\eta_I$.
Since the family $\mr\CX^+\to\CB$ is constant, we see that the relevant neighborhood
in $\CW^\eta_I$ is also constant. Since the finite group $F$ coincides with $\pi_0(S_0)$,
the assertion follows from the definitions.
\end{proof}

Observe that the nearby cycles $\psi('\IC^{0}_{\mr I})$
are of the form
$$
\psi('\IC^{0}_{\mr I})\simeq
\BC_{\mr C^I}\boxtimes
\psi('\IC^0_{\emptyset}).
$$
The only irreducible object of $\catq(Z_{0,I})$
which is such a product is the basic sheaf $'\IC_{0,\mr I}^{0}$.
By the lemma, this occurs exactly once as an irreducible constituent
of $\psi('\IC^{0}_{\mr I})$. 

\subsection{Proof of \thmref{thmconv}(2)}    \label{proof of thmconv(2)} 
To understand the nearby cycles of arbitrary objects of $\catq(Z_I)$, we need a version
of the local model for the family $\CX\to\CB$.
The results from Section~\ref{loc mod} extend to this setting 
with only notational changes.
Thus we keep the discussion here brief and focus on what we will need specifically
for the current proof. For a finite set $I$, and $\eta\in\cowts_{A}$,
in analogy with the local model $W_{I}^\eta$ for the ind-stack $Z_I$, we write
$ \CW_{I}^\eta$ for the local model for the ind-stack $\CZ_I$.

For convenience, we work with the local model $\CW^\eta_{I,\BA^1}$ with respect to the curve $\BA^1$.
By choosing local coordinates on the curve $C$,
we have the following.

\begin{thm}\label{tlocmod}
For any point $z\in \CZ_I$, there is a finite set $K$, 
and for $k\in K$,
finite sets $J_k$,
coweights $\eta_k\in\cowts_A$, and points $w_k\in\CW^{\eta_k}_{J_k,\BA^1}$
such that the following properties hold.
\begin{enumerate}
\item
In the smooth topology, there is a neighborhood of $z\in  \CZ_I$ and a neighborhood of
$$
\prod_{k\in K}w_k \in \prod_{k\in K} W^{\eta_k}_{J_k,\BA^1}
$$ 
which are isomorphic.
\item The isomorphism identifies the restrictions of the untwisted local strata
to the neighborhoods.
\item The isomorphism identifies the restrictions of the horospherical strata
to the neighborhoods.
\item
For all $k\in K$,
the degeneracy locus of the point $w_k$ is a single point of $\BA^1$.
\item
If $z$ lies in the fiber $\CZ_{0,I}\subset\CZ_I$,
then for all $k\in K$,
the point
$
w_k
$ 
is transverse at all points of $\BA^1$.
\end{enumerate}
\end{thm}

\begin{proof}
All of the assertions are natural generalizations of the results of Section~\ref{loc mod},
except for the last.
For that, what remains to be seen is the following.
Consider the substack
of the family $\CW^\eta_I$ consisting of elements with no degenerations, that is,
elements whose defining section $\sigma$ extends to a map from all of $C$ to the
open $G$-orbit $\mr X$. We must check that this substack is smooth.
But by the results of Section~\ref{loc mod}, it is isomorphic in the smooth topology to
the analogous substack of $\CZ_I$. And this substack is smooth by 
Corollary~\ref{smoothfamily}.
\end{proof}

We will apply the theorem to prove the following constructibility result.
Recall that the stack $Z_{0,\emptyset}$ is the disjoint union of the horospherical strata
$$
Z_{0,\emptyset}^{\fU(\theta^\pos)}\subset Z_{0,\emptyset},
$$
for decompositions $\fU(\theta^\pos)$ of
positive coweights
$\theta^\pos\in \cowts^\pos_{A_0}$.

\begin{thm}\label{thmncconstr}
For any positive coweight
$\theta^\pos\in \cowts^\pos_{A_0}$,
and decomposition
$\fU(\theta^\pos)$,
the restriction of the nearby cycles $\psi('\IC_{\emptyset}^{0})$ 
to the stratum
$Z_{0,\emptyset}^{\fU(\theta^\pos)}\subset Z_{0,\emptyset}$
is locally constant.
\end{thm}

\begin{proof}
It suffices to prove an analogue of the statement for the local model
discussed above.
More precisely, for a positive coweight $\eta^\pos\in\poscowts_A$,
consider the nearby cycles functor
$$
\psi:\catp(\CW^{\eta^\pos}_{\emptyset,\BA^1})\to\catp(\CW^{\eta^\pos}_{0,\emptyset,\BA^1}),
$$
and the intersection cohomology sheaf $'\IC^{\eta^\pos,0}_\emptyset$ 
of the untwisted basic stratum
$
'\CW^{\eta^\pos,0,\BA^1}_{\emptyset}\subset\CW^{\eta^\pos}_{\emptyset,\BA^1}.
$
By Theorem~\ref{tlocmod},  it suffices 
to prove that 
the restriction of the nearby cycles $\psi('\IC_{\emptyset,0}^{\eta^\pos})$ 
to the horospherical stratum
$
\CW_{0,\emptyset,\BA^1}^{\eta^\pos,\fU_0(\eta^\pos)}\subset \CW^{\eta^\pos}_{0,\emptyset,\BA^1},
$
where $\fU_0(\eta^\pos)$ is the trivial decomposition,
is locally constant in a neighborhood of a point
$w\in \CW_{0,\emptyset,\BA^1}^{\eta^\pos,\fU_0(\eta^\pos)}$.
This assertions follows by transversality:
the action of the affine line $\BA^1$ on itself by translation
lifts to an action on the family $\CW^{\eta^\pos}_{\emptyset,\BA^1}$
which is locally transitive
on the horospherical stratum $\CW_{0,\emptyset,\BA^1}^{\eta^\pos,\fU_0(\eta^\pos)}\subset \CW^{\eta^\pos}_{0,\emptyset,\BA^1}$.
\end{proof}

Finally, we finish the proof of \thmref{thmconv}(2) as follows.
By \thmref{thmconv}(1), in the category $\catp_\CH(Z_{0,I})$,
we have a filtration
$$
\CF_1\subset\CF_2\subset \psi('\IC^{0}_{\mr I})
$$
such that there is an isomorphism
$$
\CF_2/\CF_1\simeq {}'\IC_{0,\mr I}^{0}.
$$
Furthermore, the restrictions of the subsheaf $\CF_1$ and quotient sheaf
$\psi('\IC^{0}_{\mr I})/\CF_2$ to the basic untwisted stratum 
${}'Z^0_{0,\mr I}\subset Z_{0,I}$ 
are zero.

To prove the theorem,
we must check that $H^{I}_G(\CA,\CF_1)$ and
$H^{I}_G(\CA,\psi('\IC^{ 0}_{\mr I})/\CF_2)$ are bad sheaves.
By the definition of bad sheaves, it suffices to prove that $H^{I}_G(\CA,\CR)$ is a bad sheaf
for each simple constituent $\CR$ of either sheaf.
We must deal with two possible types of simple constituents $\CR$.

First, we must 
consider the case when $\CR$ is the intersection cohomology sheaf of a component 
of the partially twisted basic stratum and the component is not untwisted.
In this case, we may argue using the results of~\cite{GNhoro04} as follows.
Since we explicitly know what the convolution looks like for the intersection cohomology
sheaf of the entire basic stratum and of the untwisted basic stratum, we know
what the convolution looks like for the intersection cohomology sheaf
of our partially twisted component. Namely, it must be a sum of intersection
cohomology sheaves of components of partially twisted local strata which themselves
are not untwisted. Thus the convolution is a bad sheaf and the assertion
in this case is verified.

Now, what remains is 
the case when the restriction of $\CR$ to the entire basic stratum
is zero.
Recall that the nearby cycles $\psi('\IC^{0}_{\mr I})$
are of the form
$$
\psi('\IC^{0}_{\mr I})\simeq
\BC_{\mr C^I}\boxtimes
\psi('\IC^0_{\emptyset}).
$$
By Theorem~\ref{thmncconstr}, 
the restriction of the nearby cycles $\psi('\IC^{ 0}_{\emptyset})$ to
each stratum 
$$
Z^{\fU(\theta^\pos)}_{0,\emptyset}\subset Z_{0,\emptyset}
$$
is locally constant,
and so each of its simple constituents $\CR$ is as well.
Since the restriction of $\CR$ to the basic stratum 
is zero, we conclude that $\CR$ must be of the form
$$
\BC_{\mr C^I}\boxtimes
\IC^{ \fU(\theta^\pos)}_{0,\emptyset}(\CL),
\mbox{ for $\theta^\pos\neq 0$,}
$$
where $\CL$ is some Hecke equivariant local system on $Z_{0,\emptyset}^{\fU(\theta^\pos)}$.
Now the theorem is an immediate consequence of \propref{conv and bad sheaves} 
of \secref{proof tedious}.

\subsection{Proof of \propref{pfactfiber}}
This is very similar to the proof of \corref{corconvglue}.
It is immediately implied by the following analogue of \propref{propfactorizefibers}.

\begin{prop}
For any object $\CQ_I\in\catq(Z)^{\otimes I}$, if we fix isomorphisms
$$
\Psi(\CQ_i)\simeq\sum_{\theta\in\cowts_{A_0}}{}'\IC^\theta_{}(V^\theta_i),
$$
where $V^\theta_i$ are vector spaces,
then we have a canonical isomorphism
$$
\Psi_I(\gamma_I(\CQ_I))\simeq\sum_{\Theta:I\to\cowts_{A_0}}{}'\IC^\Theta_{\mr I}(\otimes_{i\in I}V^{\Theta(i)}_i).
$$
\end{prop}

\begin{proof}
Fix a labelling $\Theta:I\to\domcowts_{\openX}$, and 
a point $z\in {}'Z_{\mr I}^\Theta$. 
By our previous results,
to calculate the stalk of $\Psi_I(\gamma_I(\CA_I))$ at $z$, we may instead consider 
the stalk of $\Psi_I^\eta(\gamma(\CA_I))$
at a point $w\in{}' W_{\mr I}^{\eta,\Theta}$
such that $w$ is transverse at its pole points.
Then via factorization, this stalk is the tensor product of the following:
\begin{enumerate}
\item the stalk of $\Psi_{i,1}^{\eta_i}(\CA_i)$ at 
a point $w_i\in {}'W^{\eta_i}_{i,1}$
with a degeneracy of degree $\eta_i$
at its single pole point $c_i$ and no other degeneracies,
\item
the stalk of $\Psi^{\eta'}_{\emptyset,1}({}'\IC^0)$ at a point $w'\in {}' W^{\eta'}_{\emptyset,1}$
with no pole points but otherwise arbitrary degeneracies.
\end{enumerate}

In the second case, by Corollary~\ref{smoothfamily}, we are taking the stalk of the nearby cycles of the intersection cohomology sheaf
at a smooth point of the family, 
and so we may canonically identify the stalk with the trivial vector space $\BC$. 

Thus to prove the proposition,
it remains to canonically identify the stalks of $\Psi_{i,1}^{\eta_i}(\CA_i)$ at
points $w_i,w'_i\in  {}' W_{i,1}^{\eta_i}$ satisfying the conditions of the first case.
Note that the condition is equivalent to the points being transverse everywhere.

For a finite set $I$, and $\eta\in\cowts_A$, 
define the scheme
$
\Bun_{F,I}^\eta
$
to be that classifying data
$$
(c_I^\eta;\CP_F,\alpha)
$$
where $c_I^\eta\in C^\eta_I$, $\CP_F\in\Bun_F$, and $\alpha$ is a trivialization
$$
\alpha:|c_I^\eta|\times F\to\CP_F|_{|c_I^\eta|}.
$$
We have a diagram of \'etale maps
$$
 M_{I,1}^\eta
\leftarrow
 M_{I,1}^\eta\underset{C^\eta_I}{\times}\Bun_{F,I}^\eta
\to
 M_{I,1}^\eta
$$
in which the left map is the obvious projection, and the right map is defined by
$$
((c;c_{I}^{\eta}; \CP_{\maxtor}, \tau_0),(c_I^\eta;\CP_F,\alpha))
\mapsto
(c;c_{I}^{\eta}; \CP_{\maxtor}\otimes\CP_F, \tau_0).
$$

This extends to a diagram of \'etale maps
$$
 '\CW_{I,1}^\eta
\leftarrow
 '\CW_{I,1}^\eta\underset{C^\eta_I}{\times}\Bun_{F,I}^\eta
\to
 '\CW_{I,1}^\eta
$$
where the left map is the obvious projection, and the right map is defined 
in a similar way to the factorization map.

Now since $w_i,w'_i\in  {}' W_{i,1}^{\eta_i}$ have the same degeneracies
and are transverse everywhere, they are related by the above correspondence.
Thus we conclude that the stalk of 
$\Psi_{i,1}^{\eta_i,}(\CA_i)$ is the same at
$w_i, w'_i.$ 
\end{proof}


\subsection{Proofs of Propositions \ref{hw spec} and \ref{hw fuse}}  
\label{proof of hw spec}\label{proof of hw fuse}

\subsubsection{The first assertions} 
We simultaneously establish the first assertions of both propositions,
in the process checking their compatibility.

Set $I=\{1,2\}$ and $\Theta=(\theta_1,\theta_2)$,
and consider the intersection cohomology sheaf ${}' \IC_{\mr I}^\Theta$ of the stratum closure  
$${}' \ol Z^\Theta_{\mr I}\subset {}' Z_{\mr I}.
$$
Our goal is to understand: 

\begin{enumerate}
\item
 the middle-extension of ${}' \IC_{\mr I}^\Theta$
across the divisor in $Z_I$ where the pole
points collide, and specifically its structure along the stratum ${}' Z^{\theta_1+\theta_2}$
of the divisor;
\item
the nearby cycles of ${}' \IC_{\mr I}^\Theta$ at the special fiber $Z_{0,\mr I}$ in the specializing family $\CZ_{\mr I}$, and specifically its structure along the stratum ${}' Z_0^{\Theta}$
of the special fiber.
\end{enumerate}

 To do this, we will focus on the entire specializing family $\CZ_I$ in a neighborhood of
the stratum 
$
{}' Z_0^{\theta_1+\theta_2}
$
in the collision divisor of the special fiber.

\medskip

Fix a point $z_0$ in the stratum ${}' Z_0^{\theta_1+\theta_2}$. In the \'etale topology,
there is a neighborhood of $z_0$ in the family $\CZ_I$ which is isomorphic to a neighborhood of a 
point $w_0$ in the local model $\CW^\eta_I$. Furthermore, we may arrange so that $w_0$
is transverse at its single pole point.

\medskip

Once we have removed a point from the curve $C$,
factorization provides a neighborhood of $w_0$ in the local model $\CW^\eta_I$ which is
isomorphic to a product of neighborhoods of points
\begin{enumerate}
\item $w_{0,1}$ in the local model $\CW^{\theta_1+\theta_2}_{I,1}$ such that $w_{0,1}$
has a degeneracy of degree ${\theta_1+\theta_2}$ at its single pole point, and no
other degeneracies,
\item $w'_{0,1}$ in the local model $\CW^{\eta'}_{\emptyset,1}$ with no pole points but
otherwise arbitrary degeneracies.
\end{enumerate}

In the second case, by Corollary~\ref{smoothfamily}, our calculation reduces to the middle-extension
and nearby cycles of the intersection cohomology sheaf
along the smooth locus of the family $\CW^{\eta'}_{\emptyset,1}$.

In the first case, the closure of the stratum
$ {}' W^{\theta_1+\theta_2,\Theta}_{\mr I, 1}
$
in the family
$\CW^{\theta_1+\theta_2}_{I,1}
$
lies in the globally transverse locus $\CW^{\theta_1+\theta_2, +}_{I,1}$. This is canonically isomorphic
to the constant family with fiber the untwisted transverse base ${}' M^{\theta_1+\theta_2, +}_{I,1}$.

Thus we see in this case as well
that our calculation reduces to the middle-extension
and nearby cycles of the intersection cohomology sheaf
along the smooth stack ${}' M^{\theta_1+\theta_2, +}_{I,1}$. 
This provides a canonical isomorphism of the stalks of these sheaves at
any points. By our previous results, this identification of stalks
is all that is needed to confirm the first assertions of Propositions \ref{hw spec} and \ref{hw fuse}.

\subsubsection{The second assertions}

The second assertion of \propref{hw spec} follows immediately from the description
of the strata and ind-scheme structure of $X(\CK) - (X - \mr X)(\CK)$. 

Since $\poscowts_X$
is strictly convex, to establish the scond assertion of \propref{hw fuse},
it suffices to prove the stronger statement:
if $'\IC^\theta$ is a constituent of $'\IC^{\theta_1}\product{}'\IC^{\theta_2}$,
then $(\theta_1+\theta_2)-\theta\in \poscowts_X$.

This follows immediately from the two assertions of \propref{hw spec}.

\newpage

\bibliographystyle{alpha}
\bibliography{ref}

\end{document}